\documentclass[english,reprint, superscriptaddress, secnumarabic, nofootinbib, nobibnotes]{revtex4-1}
\usepackage[T1]{fontenc}
\usepackage[latin9]{inputenc}
\usepackage{xcolor}
\usepackage{pdfcolmk}
\usepackage{array}
\usepackage{fancybox}
\usepackage{calc}
\usepackage{units}
\usepackage{textcomp}
\usepackage{multirow}
\usepackage{amsmath}
\usepackage{amsthm}
\usepackage{amssymb}
\usepackage{mathdots}
\usepackage{graphicx}
\usepackage{wasysym}
\usepackage[all]{xy}
\PassOptionsToPackage{normalem}{ulem}
\usepackage{ulem}

\makeatletter

\providecommand{\tabularnewline}{\\}
\providecolor{lyxadded}{rgb}{0,0,1}
\providecolor{lyxdeleted}{rgb}{1,0,0}

\theoremstyle{plain}
\newtheorem{thm}{\protect\theoremname}[section]
  \theoremstyle{definition}
  \newtheorem{defn}[thm]{\protect\definitionname}
\newenvironment{lyxlist}[1]
{\begin{list}{}
{\settowidth{\labelwidth}{#1}
 \setlength{\leftmargin}{\labelwidth}
 \addtolength{\leftmargin}{\labelsep}
 }}
{\end{list}}

\usepackage[all]{xy}

\usepackage{morefloats}

\usepackage{placeins}

\usepackage{lineno}

\newcommand{\id}{\mathsf{id}}
\newcommand{\di}{\mathsf{di}}
\newcommand{\qdi}{\mathsf{qdi}}
\newcommand{\qPr}{\mathsf{qPr}}
\newcommand{\sodd}{\mathsf{{s_{odd}}}}
\newcommand{\R}{\mathsf{R}}
\newcommand{\E}{\mathsf{E}}

\def\frontmatter@abstractheading{}

\makeatletter
\renewcommand{\p@subsection}{}
\renewcommand{\p@subsubsection}{}
\makeatother

\makeatother

\usepackage{babel}
  \providecommand{\definitionname}{Definition}
\providecommand{\theoremname}{Theorem}

\begin{document}

\title{Context-Content Systems of Random Variables: The Contextuality-by-Default
Theory}

\author{Ehtibar N.\ Dzhafarov}

\email{E-mail: ehtibar@purdue.edu }

\selectlanguage{english}%

\affiliation{Purdue University, USA }

\author{Janne V.\ Kujala}

\email{E-mail: jvk@iki.fi }

\selectlanguage{english}%

\affiliation{University of Jyväskylä, Finland }
\begin{abstract}
This paper provides a systematic yet accessible presentation of the
Contextuality-by-Default theory. The consideration is confined to
finite systems of categorical random variables, which allows us to
focus on the basics of the theory without using full-scale measure-theoretic
language. Contextuality-by-Default is a theory of random variables
identified by their contents and their contexts, so that two variables
have a joint distribution if and only if they share a context. Intuitively,
the content of a random variable is the entity the random variable
measures or responds to, while the context is formed by the conditions
under which these measurements or responses are obtained. A system
of random variables consists of stochastically unrelated ``bunches,''
each of which is a set of jointly distributed random variables sharing
a context. The variables that have the same content in different contexts
form ``connections'' between the bunches. A probabilistic coupling
of this system is a set of random variables obtained by imposing a
joint distribution on the stochastically unrelated bunches. A system
is considered noncontextual or contextual according to whether it
can or cannot be coupled so that the joint distributions imposed on
its connections possess a certain property (in the present version
of the theory, ``maximality''). We present a criterion of contextuality
for a special class of systems of random variables, called cyclic
systems. We also introduce a general measure of contextuality that
makes use of (quasi-)couplings whose distributions may involve negative
numbers or numbers greater than 1 in place of probabilities. 

\end{abstract}
\maketitle

\section{Introduction}

Contextuality-by-Default (CbD) is an approach to probability theory,
specifically, to the theory of random variables. CbD is not a model
of empirical phenomena, and it cannot be corroborated or falsified
by empirical data. However, it provides a sophisticated conceptual
framework in which one can describe empirical data and formulate models
that involve random variables. 

In Kolmogorovian Probability Theory (KPT) random variables are understood
as measurable functions mapping from one (domain) probability space
into another (codomain) probability space. CbD can be viewed as a
theory within the framework of KPT if the latter is understood as
allowing for multiple domain probability spaces, freely introducible
and unrelated to each other. However, CbD can also be (in fact, is
better) formulated with no reference to domain probability spaces,
with random variables understood as entities identified by their probability
distributions and their unique labels within what can be called sets
of random variables ``in existence'' or ``in play.'' 

Although one cannot deal with probability distributions without the
full-fledged measure-theoretic language, we avoid technicalities some
readers could find inhibitive by focusing in this paper on \emph{finite}
systems of \emph{categorical random variables} (those with finite
numbers of possible values). Virtually all of the content of this
paper, however, is generalizable \emph{mutatis mutandis} to arbitrary
systems of arbitrary random entities.

\subsection{A convention}

In the following we introduce sets of random variables classified
in two ways, by their \emph{contexts} and by their \emph{contents},
and we continue to speak of contexts and contents throughout the paper.
The two terms combine nicely, but they are also easily confused in
reading. For this reason, in this paper we do violence to English
grammar and write ``conteXt'' and ``conteNt'' when we use these
words as special terms.

\subsection{Two conteNts in two conteXts}

We begin with a simple example. A person randomly chosen from some
population is asked two questions, $q$ and $q'$. Say, $q=\textnormal{\textquotedblleft Do you like bees?\textquotedblright}$
and $q'=\textnormal{\textquotedblleft Do you like to smell flowers?\textquotedblright}$.
The answer to the first question (Yes or No) is a random variable
whose \emph{identity} (that which allows one to uniquely identify
it within the class of all random variables being considered) clearly
includes $q$, so it can be denoted $R_{q}$. We will refer to the
question $q$ as the \emph{conteNt} of the random variable $R_{q}$.
The second random variable then can be denoted $R_{q'},$ and its
conteNt is $q'$. The set of all random variables being considered
here consists of $R_{q}$ and $R_{q'}$, and we do not confuse them
because they have distinct conteNts: we know which of the two responses
answers which question. 

The two random variables have a joint distribution that can be presented,
because they are binary, by values of the three probabilities
\[
\begin{array}{c}
\Pr\left[R_{q}=\textnormal{Yes}\right],\quad\Pr\left[R_{q'}=\textnormal{Yes}\right],\\
\\
\Pr\left[R_{q}=\textnormal{Yes}\textnormal{ and }R_{q'}=\textnormal{Yes}\right].
\end{array}
\]
The joint distribution exists because the two responses, $R_{q}$
and $R_{q'}$, occur together in a well-defined empirical sense: the
empirical sense of ``togetherness'' of the responses here is ``to
be given by one and the same person.'' In other situations the empirical
meaning can be different, e.g., ``to be recorded in the same trial.'' 

Our example is too simple for our purposes. Let us assume therefore
that the two questions $q,q'$ are asked under varying controlled
conditions, e.g., one randomly chosen person can be asked these questions
after having watched a movie about the killer bees spreading northwards
(let us call this condition $c$), another after watching a movie
about deciphering the waggle dances of the honey bees ($c'$). Most
people would consider $q$ as one and the same question whether posed
under the condition $c$ or the condition $c'$; and the same applies
to the question $q'$. In other words, the conteNts $q$ and $q'$
of the two respective random variables would normally be considered
unchanged by the conditions $c$ and $c'$. 

However, the random variables themselves (the responses) are clearly
affected by these conditions. In particular, nothing guarantees that
the joint distribution of $\left(R_{q},R_{q'}\right)$ will be the
same under the two conditions. It is necessary therefore to include
$c$ and $c'$ in the description of the random variables representing
the responses. We will call $c$ and $c'$ \emph{conteXts} of (or
for) the corresponding random variables and present them as $R_{q}^{c},R_{q'}^{c},R_{q}^{c'},R_{q'}^{c'}$.
There are now four random variables in play, and we do not confuse
them because each of them is uniquely identified by its conteNt and
its conteXt.

\subsection{Jointly distributed versus stochastically unrelated random variables}

In each of the two conteXts, the two random variables are jointly
distributed, i.e., we have well-defined probabilities
\[
\left.\begin{array}{c}
\Pr\left[R_{q}^{c}=\textnormal{Yes}\right],\\
\\
\Pr\left[R_{q'}^{c}=\textnormal{Yes}\right],\\
\\
\Pr\left[R_{q}^{c}=\textnormal{Yes}\textnormal{ and }R_{q'}^{c}=\textnormal{Yes}\right]
\end{array}\right\} \textnormal{in conteXt }c,
\]
and
\[
\left.\begin{array}{c}
\Pr\left[R_{q}^{c'}=\textnormal{Yes}\right],\\
\\
\Pr\left[R_{q'}^{c'}=\textnormal{Yes}\right],\\
\\
\Pr\left[R_{q}^{c'}=\textnormal{Yes}\textnormal{ and }R_{q'}^{c'}=\textnormal{Yes}\right]
\end{array}\right\} \textnormal{in conteXt }c'.
\]
No joint probabilities, however, are defined between the random variables
picked from different conteXts. We cannot determine such probabilities
as 
\[
\begin{array}{c}
\Pr\left[R_{q}^{c}=\textnormal{Yes}\textnormal{ and }R_{q'}^{c'}=\textnormal{Yes}\right],\\
\\
\Pr\left[R_{q}^{c}=\textnormal{Yes}\textnormal{ and }R_{q}^{c'}=\textnormal{Yes}\right],\\
\\
\Pr\left[R_{q}^{c}=\textnormal{Yes}\textnormal{ and }R_{q}^{c'}=\textnormal{Yes}\textnormal{ and }R_{q'}^{c'}=\textnormal{Yes}\right],\\
\\
etc.
\end{array}
\]

We express this important fact by saying that any two variables recorded
in different conteXts are \emph{stochastically unrelated}. The reason
for stochastic unrelatedness is simple: no random variable in conteXt
$c$ can co-occur with any random variable in conteXt $c'$ in the
same empirical sense in which two responses co-occur within either
of these conteXts, because $c$ and $c'$ are mutually exclusive conditions.
The empirical sense of co-occurrence in our example is ``to be given
by the same person,'' and we have assumed that a randomly chosen
person is either shown one movie or another. If some respondents were
allowed to watch both movies before responding, we would have to redefine
the classification of our random variables by introducing a third
conteXt, $c''=\left(c,c'\right)$. We would then have three pairwise
mutually exclusive conteXts, $c,c',c''$, and six random variables,
$R_{q}^{c},R_{q'}^{c},R_{q}^{c'},R_{q'}^{c'},R_{q}^{c''},R_{q'}^{c''}$,
such that, e.g., $R_{q}^{c''}$ is jointly distributed with $R_{q'}^{c''}$
but not with $R_{q}^{c}$. 

In case one is tempted to consider joint probabilities involving $R_{q}^{c}$
and $R_{q}^{c'}$ simply equal to zero (because these two responses
never co-occur), this thought should be dismissed. Indeed, then all
four joint probabilities, 
\[
\begin{array}{c}
\Pr\left[R_{q}^{c}=\textnormal{Yes}\textnormal{ and }R_{q'}^{c'}=\textnormal{Yes}\right],\\
\\
\Pr\left[R_{q}^{c}=\textnormal{Yes}\textnormal{ and }R_{q'}^{c'}=\textnormal{No}\right],\\
\\
\Pr\left[R_{q}^{c}=\textnormal{No}\textnormal{ and }R_{q'}^{c'}=\textnormal{Yes}\right],\\
\\
\Pr\left[R_{q}^{c}=\textnormal{No}\textnormal{ and }R_{q'}^{c'}=\textnormal{No}\right],
\end{array}
\]
would have to be equal to zero, which is not possible as they should
sum to 1. These probabilities are not zero, they are \emph{undefined}.

\begin{figure}
\begin{centering}
\begin{tabular}{|c|c|c}
\cline{1-2} 
$R_{q}^{c}$$\begin{array}{cccc}
\\
\\
\end{array}$ & $R_{q'}^{c}$$\begin{array}{cccc}
\\
\\
\end{array}$ & $c$\tabularnewline
\cline{1-2} 
$R_{q}^{c'}$$\begin{array}{cccc}
\\
\\
\end{array}$ & $R_{q'}^{c'}$$\begin{array}{cccc}
\\
\\
\end{array}$ & $c'$\tabularnewline
\cline{1-2} 
\multicolumn{1}{c}{$q$} & \multicolumn{1}{c}{$q'$} & $\boxed{\boxed{\mathcal{A}}}$\tabularnewline
\end{tabular}
\par\end{centering}

\caption{\label{fig:system 2by2}A conteXt-conteNt matrix for system $\mathcal{A}$
in our opening example. The system consists of two bunches $R^{c}=\left(R_{q}^{c},R_{q'}^{c}\right)$,
$R^{c'}=\left(R_{q}^{c'},R_{q'}^{c'}\right)$ defined by (or defining)
the conteXts $c$ and $c'$, respectively. The notation $R^{c},R^{c'}$
reflects the fact that each bunch is a single random variable in its
own right, because its components are jointly distributed. The system
has two connections $\left(R_{q}^{c},R_{q}^{c'}\right)$, $\left(R_{q'}^{c},R_{q'}^{c'}\right)$
defined by (or defining) the conteNts $q$ and $q'$, respectively.
The connections are not random variables because their components
are stochastically unrelated. System $\mathcal{A}$ may be contextual
or noncontextual, depending on the distributions of the bunches $R^{c},R^{c'}$.}
 
\end{figure}

\subsection{Bunches and connections in conteXt-conteNt matrices }

The picture of the system consisting of our four random variables
is now complete. Let us call this system $\mathcal{A}$. It is an
example of a \emph{conteXt-conteNt (c-c)} \emph{system} of random
variables, and it can be schematically presented in the form of the
\emph{conteXt-conteNt (c-c) matrix} in Fig.~\ref{fig:system 2by2}.
All random variables in a \mbox{c-c} system are double-indexed: the
lower index indicates their conteNt, the upper index indicates their
conteXt. The random variables within each conteXt are jointly distributed,
they form what we call a \emph{bunch} (of random variables). One bunch
corresponds to one conteXt and occupies one row of the \mbox{c-c}
matrix. Any two variables that belong to different bunches are stochastically
unrelated. However, a random variable in conteXt $c$ may have a counterpart
in conteXt $c'$ that shares the same conteNt with it. In the system
$\mathcal{A}$ this is true for each of the two random variables in
$c$ (or $c'$): $R_{q}^{c}$ and $R_{q}^{c'}$ represent answers
to one and the same question, and so do $R_{q'}^{c}$ and $R_{q'}^{c'}$.
The set (in our example, the pair) of all random variables sharing
the same conteNt is called a \emph{connection }(because they bridge
stochastically unrelated bunches). One connection, in our example
$\left(R_{q}^{c},R_{q}^{c'}\right)$ or $\left(R_{q'}^{c},R_{q'}^{c'}\right)$,
corresponds to one conteNt and occupies one column of the \mbox{c-c}
matrix. 

We will see in Section \ref{sec: Cyclic-c-c-systems} that $\mathcal{A}$
is the simplest system within the class of so-called cyclic systems.
It can, of course, model more than the opening example with bees,
flowers, and movies. The variety of possible applications is great,
both within psychology and without. The conteNts $q,q'$ can be two
physical properties, e.g., spins of a particle measured at two moments
in time separated by a fixed interval. The conteXts $c,c'$ can then
be, respectively, the presence and absence of a third measurement
made prior to these two measurements. Alternatively, $c$ and $c'$
could be two orders in which the two measurements are conducted: $c$
standing for ``first $q$ then $q'$'' and $c'$ for ``first $q'$
then $q$.'' In sociology and psychology the prominently studied
analogue of the latter example is the paradigm in which two questions
are posed in two possible orders (Moore, 2001; Wang \& Busemeyer,
2013; Wang et al., 2014). One can also think of questions posed in
two different forms, in two different languages, or asked of the representatives
of two distinct populations (say, male and female). There is also
an inexhaustible variety of psychophysical applications. For instance,
$q$ and $q'$ may be visual stimuli, and $c,c'$ may be any two variants
of the conditions under which they are presented, such as the time
interval or spatial separation between them, or two versions of a
previously presented adapting stimulus. 

Fig.~\ref{fig: system general} shows a \mbox{c-c} matrix representation
of a more complex system, with three bunches and three connections. 

A generalization to arbitrary \mbox{c-c} systems should be obvious:
given a set of conteXts and a set of conteNts the cells in a \mbox{c-c}
matrix can be filled (or left empty) in all possible ways, although
constraints could be imposed to exclude matrices that are uninteresting
for contextuality analysis (e.g., empty matrices, or those with a
single random variable in a connection or in a bunch).

\begin{figure}
\begin{centering}
\begin{tabular}{|c|c|c|c}
\cline{1-3} 
$R_{1}^{1}$$\begin{array}{cc}
\\
\\
\end{array}$ & $R_{2}^{1}$$\begin{array}{cc}
\\
\\
\end{array}$ & $\cdot$$\begin{array}{cc}
\\
\\
\end{array}$ & $c_{1}$\tabularnewline
\cline{1-3} 
$R_{1}^{2}$$\begin{array}{cc}
\\
\\
\end{array}$ & $R_{2}^{2}$$\begin{array}{cc}
\\
\\
\end{array}$ & $R_{3}^{2}$$\begin{array}{cc}
\\
\\
\end{array}$ & $c_{2}$\tabularnewline
\cline{1-3} 
$R_{1}^{3}$$\begin{array}{cc}
\\
\\
\end{array}$ & $\cdot$$\begin{array}{cc}
\\
\\
\end{array}$ & $R_{3}^{3}$$\begin{array}{cc}
\\
\\
\end{array}$ & $c_{3}$\tabularnewline
\cline{1-3} 
\multicolumn{1}{c}{$q_{1}$} & \multicolumn{1}{c}{$q_{2}$} & \multicolumn{1}{c}{$q_{3}$} & $\boxed{\boxed{\mathcal{B}}}$\tabularnewline
\end{tabular}
\par\end{centering}

\caption{\label{fig: system general}A \mbox{c-c} matrix representation of
a \mbox{c-c} system $\mathcal{B}$ of random variables. The seven
random variables are grouped into three bunches (shown by the rows
of the matrix) and into three connections (shown by the columns).
The bunches are defined by (or define) three conteXts. The connections
are defined by (or define) three conteNts. The empty cells (shown
with a dot for emphasis) correspond to the cases when a given conteNt
is not represented (measured, responded to) in a given conteXt. The
variables within a bunch are jointly distributed, so we have three
random variables $R^{1}=\left(R_{1}^{1},R_{2}^{1}\right)$, $R^{2}=\left(R_{1}^{2},R_{2}^{2},R_{3}^{2}\right)$,
and $R^{3}=\left(R_{1}^{3},R_{3}^{3}\right)$. The connections $\left(R_{1}^{1},R_{1}^{2},R_{1}^{3}\right)$,
$\left(R_{2}^{1},R_{2}^{2}\right)$, and $\left(R_{3}^{2},R_{3}^{3}\right)$
are not random variables because no two random variables within a
connection are jointly distributed.}
\end{figure}

\subsection{\label{sub: ConteXts-and-conteNts}ConteXts and conteNts are non-unique
but distinct from each other}

How do we know that in our opening example the question $q$ and not
the movie $c$ determines the conteNt of the response, viewed as a
random variable? How do we know that the movie $c$ and not the question
$q$ determines the conteXt of this response? The answer is: we don't.
Some theory or tradition outside the mathematical theory of CbD tells
us what the conteXts and the conteNts in a given situation are, and
then the mathematical computations may commence. In these computations,
whatever conteXts and conteNts are given to us, they are treated as
strictly distinct entities because the respective bunches and connections
they define are fundamentally different: bunches are (multicomponent)
random variables, while connections are groups of pairwise stochastically
unrelated ones. 

It would be a completely different system if the conteNts in our opening
example were defined not just by the question asked but also by the
movie previously watched. The \mbox{c-c} matrix would then be as
shown in Fig.~\ref{fig:system_noconnect}. No conteNt in the system
$\mathcal{A}'$ occurs more than once, so there is nothing to bridge
the two bunches. It is not wrong to present the experiment with the
questions and movies in this way, it may very well be the best way
of treating this situation from the point of view of some empirical
model, but the resulting system is not interesting for our contextuality
analysis. The latter is yet to be introduced, but it should be sufficiently
clear if we say that the system $\mathcal{A}'$ is uninteresting because
contextuality pertains to how the random variables that share conteNts
differ in different conteXts.

To prevent turning this discussion into a game of superficial semantics,
it would not do to point out that $q_{1}$ and $q_{3}$ in the system
$\mathcal{A}'$ share ``part'' of their conteNt, and hence $R_{q_{1}}^{c}$
and $R_{q_{3}}^{c'}$ can be related to each other on these grounds.
If $R_{q_{1}}^{c}$ and $R_{q_{3}}^{c'}$ are members of the same
connection, then they should have the same conteNt, by definition.
A conteNts is, logically, merely a label for a connection.

\begin{figure}
\begin{centering}
\begin{tabular}{|c|c|c|c|c}
\cline{1-4} 
$R_{q_{1}}^{c}$$\begin{array}{cc}
\\
\\
\end{array}$ & $R_{q_{2}}^{c}$$\begin{array}{cc}
\\
\\
\end{array}$ & $\cdot$$\begin{array}{cc}
\\
\\
\end{array}$ & $\cdot$$\begin{array}{cc}
\\
\\
\end{array}$ & $c$\tabularnewline
\cline{1-4} 
$\cdot$$\begin{array}{cc}
\\
\\
\end{array}$ & $\cdot$$\begin{array}{cc}
\\
\\
\end{array}$ & $R_{q_{3}}^{c'}$$\begin{array}{cc}
\\
\\
\end{array}$ & $R_{q_{4}}^{c'}$$\begin{array}{cc}
\\
\\
\end{array}$ & $c'$\tabularnewline
\cline{1-4} 
\multicolumn{1}{c}{$q_{1}=\left(q,c\right)$} & \multicolumn{1}{c}{$q_{2}=\left(q',c\right)$} & \multicolumn{1}{c}{$q_{3}=\left(q,c'\right)$} & \multicolumn{1}{c}{$q_{4}=\left(q',c'\right)$} & $\boxed{\boxed{\mathcal{A}'}}$\tabularnewline
\end{tabular}
\par\end{centering}

\caption{\label{fig:system_noconnect}The same opening example as in Fig.~\ref{fig:system 2by2},
but represented by a different system, $\mathcal{A}'$. In this system
the conteXts are the same as in $\mathcal{A}$, but each new conteNt
(question) includes as its part the original conteXt in which it occurs
(the movie watched). The joint distributions within the two bunches
remain unchanged, but the system loses connections between the bunches.
Such a system is trivially noncontextual. }
\end{figure}

A symmetrical opposite of the system $\mathcal{A}'$ is to include
the questions asked into the conteXts in which they are being asked.
This creates the \mbox{c-c} matrix $\mathcal{A}''$ shown in Fig.~\ref{fig:system_nojoint}.
Since the empirical meaning of co-occurrence in our example is ``to
be given by the same person,'' representing our opening example by
the system $\mathcal{A}''$ amounts to simply ignoring the observed
joint events. One only records (and estimates probabilities of) the
individual events, as if the paired questions were asked separately
of different respondents. This would not be a reasonable way of representing
the situation (as it involves ignoring available information), but
it is logically possible.

The system $\mathcal{A}''$ becomes a reasonable representation, however,
in fact the only ``natural'' one, if the empirical procedure is
modified and the questions are indeed asked one at a time rather than
in pairs. Then the responses to questions about the bees and about
the flowers, whether they are given after having watched the same
movie or different movies, come from different respondents, and their
joint probabilities are undefined. 

The system $\mathcal{A}''$ has the same connections as $\mathcal{A}$,
but it is as uninteresting from the point of view of contextuality
analysis as the system $\mathcal{A}'$. A system without joint distributions
(i.e., one in which every bunch contains a single random variable)
is always trivially noncontextual. 

For completeness, we should also consider a radical point of view
that combines those in the systems $\mathcal{A}''$ and $\mathcal{A}'$.
It is shown in Fig.~\ref{fig:system4separate}: every conteXt and
every conteNt include information about both the question being asked
and the conditions under which it is asked. This creates four unique
conteNts in a one-to-one correspondence with four unique conteXts.

\begin{figure}
\begin{centering}
\begin{tabular}{|c|c|c}
\cline{1-2} 
$R_{q}^{c_{1}}$$\begin{array}{cccc}
\\
\\
\end{array}$ & $\cdot$$\begin{array}{cccc}
\\
\\
\end{array}$ & $c_{1}=\left(c,q\right)$\tabularnewline
\cline{1-2} 
$R_{q}^{c_{2}}$$\begin{array}{cccc}
\\
\\
\end{array}$ & $\cdot$$\begin{array}{cccc}
\\
\\
\end{array}$ & $c_{2}=\left(c',q\right)$\tabularnewline
\cline{1-2} 
$\cdot$$\begin{array}{cccc}
\\
\\
\end{array}$ & $R_{q'}^{c_{3}}$$\begin{array}{cccc}
\\
\\
\end{array}$ & $c_{3}=\left(c,q'\right)$\tabularnewline
\cline{1-2} 
$\cdot$$\begin{array}{cccc}
\\
\\
\end{array}$ & $R_{q'}^{c_{4}}$$\begin{array}{cccc}
\\
\\
\end{array}$ & $c_{4}=\left(c',q'\right)$\tabularnewline
\cline{1-2} 
\multicolumn{1}{c}{$q$} & \multicolumn{1}{c}{$q'$} & $\boxed{\boxed{\mathcal{A}''}}$\tabularnewline
\end{tabular}
\par\end{centering}

\caption{\label{fig:system_nojoint}The same opening example as in Fig.~\ref{fig:system 2by2},
but represented by a different system, $\mathcal{A}''$. In this system
the conteNts are the same as in $\mathcal{A}$, but each new conteXt
is identified by both the movie watched and the question asked. The
connections are well-defined here, but each bunch contains a single
variable. Such a system is trivially noncontextual. This system would
more aptly be used to describe the experiment in which the questions
are posed one question per respondent rather than in pairs. }
\end{figure}

\begin{figure}
\begin{centering}
\begin{tabular}{|c|c|c|c|c}
\cline{1-4} 
$R_{q_{1}}^{c_{1}}$$\begin{array}{cc}
\\
\\
\end{array}$ & $\cdot$$\begin{array}{cc}
\\
\\
\end{array}$ & $\cdot$$\begin{array}{cc}
\\
\\
\end{array}$ & $\cdot$$\begin{array}{cc}
\\
\\
\end{array}$ & $c_{1}=\left(c,q\right)$\tabularnewline
\cline{1-4} 
$\cdot$$\begin{array}{cc}
\\
\\
\end{array}$ & $R_{q_{2}}^{c_{2}}$$\begin{array}{cc}
\\
\\
\end{array}$ & $\cdot$$\begin{array}{cc}
\\
\\
\end{array}$ & $\cdot$$\begin{array}{cc}
\\
\\
\end{array}$ & $c_{2}=\left(c',q\right)$\tabularnewline
\cline{1-4} 
$\cdot$$\begin{array}{cc}
\\
\\
\end{array}$ & $\cdot$$\begin{array}{cc}
\\
\\
\end{array}$ & $R_{q_{3}}^{c_{3}}$$\begin{array}{cc}
\\
\\
\end{array}$ & $\cdot$$\begin{array}{cc}
\\
\\
\end{array}$ & $c_{3}=\left(c,q'\right)$\tabularnewline
\cline{1-4} 
$\cdot$$\begin{array}{cc}
\\
\\
\end{array}$ & $\cdot$$\begin{array}{cc}
\\
\\
\end{array}$ & $\cdot$$\begin{array}{cc}
\\
\\
\end{array}$ & $R_{q_{4}}^{c_{4}}$$\begin{array}{cc}
\\
\\
\end{array}$ & $c_{4}=\left(c',q'\right)$\tabularnewline
\cline{1-4} 
\multicolumn{1}{c}{$q_{1}=\left(q,c\right)$} & \multicolumn{1}{c}{$q_{2}=\left(q,c'\right)$} & \multicolumn{1}{c}{$q_{3}=\left(q',c\right)$} & \multicolumn{1}{c}{$q_{4}=\left(q',c'\right)$} & $\boxed{\boxed{\mathcal{A}'''}}$\tabularnewline
\end{tabular}
\par\end{centering}

\caption{\label{fig:system4separate}A combination of the modifications shown
in Figs. \ref{fig:system_nojoint} and \ref{fig:system_noconnect}:
the original conteNts are treated as part of new conteXts, and the
original conteXts are treated as part of new conteNts: the system
consists of four unrelated to each other random variables, and it
is trivially noncontextual. The order in which the $q$'s and $c$'s
are listed in the definitions of the conteXts and the conteNts is
immaterial.}
\end{figure}

\subsection{Multiple connections crossing multiple bunches}

There is another aspect of the non-uniqueness of representing an empirical
situation by a \mbox{c-c} system. As we see in both system $\mathcal{A}$
of Fig.~\ref{fig:system 2by2} and system $\mathcal{B}$ of Fig.~\ref{fig: system general},
two or more bunches may very well intersect with the same two or more
connections. In the system $\mathcal{B}$, the conteNts $q_{1}$ and
$q_{2}$ are represented by $R_{1}^{1},R_{2}^{1}$ in the conteXt
$c_{1}$ (a short way of saying ``in the bunch labeled by the conteXt
$c_{1}$''), and the same $q_{1}$ and $q_{2}$ are represented by
$R_{1}^{2},R_{2}^{2}$ in the conteXt $c_{2}$. This may make it desirable
(but by no means necessary) to introduce a new conteNt $q_{12}$ and
the corresponding connection formed by $R_{12}^{1}=\left(R_{1}^{1},R_{2}^{1}\right)$
and $R_{12}^{2}=\left(R_{1}^{2},R_{2}^{2}\right)$. In our conceptual
framework this means replacing $\mathcal{B}$ with another system,
shown in Fig.~\ref{fig: alternative1-1}. The new system has a different
set of conteNts and its contextuality analysis generally will not
coincide with that of the system $\mathcal{B}$. 

\begin{figure}
\begin{centering}
\begin{tabular}{|c|c|c|c|c}
\cline{1-4} 
$R_{1}^{1}$$\begin{array}{cc}
\\
\\
\end{array}$ & $R_{2}^{1}$$\begin{array}{cc}
\\
\\
\end{array}$ & $\cdot$$\begin{array}{cc}
\\
\\
\end{array}$ & $R_{12}^{1}$$\begin{array}{cc}
\\
\\
\end{array}$ & $c_{1}$\tabularnewline
\cline{1-4} 
$R_{1}^{2}$$\begin{array}{cc}
\\
\\
\end{array}$ & $R_{2}^{2}$$\begin{array}{cc}
\\
\\
\end{array}$ & $R_{3}^{2}$$\begin{array}{cc}
\\
\\
\end{array}$ & $R_{12}^{2}$$\begin{array}{cc}
\\
\\
\end{array}$ & $c_{2}$\tabularnewline
\cline{1-4} 
$R_{1}^{3}$$\begin{array}{cc}
\\
\\
\end{array}$ & $\cdot$$\begin{array}{cc}
\\
\\
\end{array}$ & $R_{3}^{3}$$\begin{array}{cc}
\\
\\
\end{array}$ & $\cdot$$\begin{array}{cc}
\\
\\
\end{array}$ & $c_{3}$\tabularnewline
\cline{1-4} 
\multicolumn{1}{c}{$q_{1}$} & \multicolumn{1}{c}{$q_{2}$} & \multicolumn{1}{c}{$q_{3}$} & \multicolumn{1}{c}{$q_{12}$} & $\boxed{\boxed{\mathcal{B}'}}$\tabularnewline
\end{tabular}
\par\end{centering}

\caption{\label{fig: alternative1-1}If one wishes to establish connection
between the pair $\left(R_{1}^{1},R_{2}^{1}\right)$ and the pair
$\left(R_{1}^{2},R_{2}^{2}\right)$, the system $\mathcal{B}$ of
Fig.~\ref{fig: system general} has to be redefined: a new conteNt,
$q_{12}$, should be created that adds $R_{12}^{1}=\left(R_{1}^{1},R_{2}^{1}\right)$
and $R_{12}^{2}=\left(R_{1}^{2},R_{2}^{2}\right)$ to the two respective
bunches (for $c_{1}$ and $c_{2}$) as separate variables, each of
which is jointly distributed with the other variables within the same
bunch. }
\end{figure}

\begin{figure}
\begin{centering}
\begin{tabular}{|c|c|c|c|c|c}
\cline{1-5} 
$R_{1}^{1}$$\begin{array}{cc}
\\
\\
\end{array}$ & $R_{2}^{1}$$\begin{array}{cc}
\\
\\
\end{array}$ & $\cdot$$\begin{array}{cc}
\\
\\
\end{array}$ & $R_{12}^{1}$$\begin{array}{cc}
\\
\\
\end{array}$ & $\cdot$$\begin{array}{cc}
\\
\\
\end{array}$ & $c_{1}$\tabularnewline
\cline{1-5} 
$R_{1}^{2}$$\begin{array}{cc}
\\
\\
\end{array}$ & $R_{2}^{2}$$\begin{array}{cc}
\\
\\
\end{array}$ & $R_{3}^{2}$$\begin{array}{cc}
\\
\\
\end{array}$ & $R_{12}^{2}$$\begin{array}{cc}
\\
\\
\end{array}$ & $R_{13}^{2}$$\begin{array}{cc}
\\
\\
\end{array}$ & $c_{2}$\tabularnewline
\cline{1-5} 
$R_{1}^{3}$$\begin{array}{cc}
\\
\\
\end{array}$ & $\cdot$$\begin{array}{cc}
\\
\\
\end{array}$ & $R_{3}^{3}$$\begin{array}{cc}
\\
\\
\end{array}$ & $\cdot$$\begin{array}{cc}
\\
\\
\end{array}$ & $R_{13}^{3}$$\begin{array}{cc}
\\
\\
\end{array}$ & $c_{3}$\tabularnewline
\cline{1-5} 
\multicolumn{1}{c}{$q_{1}$} & \multicolumn{1}{c}{$q_{2}$} & \multicolumn{1}{c}{$q_{3}$} & \multicolumn{1}{c}{$q_{12}$} & \multicolumn{1}{c}{$q_{13}$} & $\boxed{\boxed{\mathcal{B}''}}$\tabularnewline
\end{tabular}
\par\end{centering}

\caption{\label{fig: alternative1-2}Further redefinition of system $\mathcal{B}$:
in addition to the new conteNt and connection shown in Fig.~\ref{fig: alternative1-1},
a conteNt $q_{13}$ is introduced adding $R_{13}^{2}=\left(R_{1}^{2},R_{3}^{2}\right)$
and $R_{13}^{3}=\left(R_{1}^{3},R_{3}^{3}\right)$ to the two respective
bunches (for $c_{2}$ and $c_{3}$). }
\end{figure}

One can analogously introduce a new connection formed by $R_{13}^{2}=\left(R_{1}^{2},R_{3}^{2}\right)$
and $R_{13}^{3}=\left(R_{1}^{3},R_{3}^{3}\right)$ to bridge the bunches
for the conteXts $c_{2}$ and $c_{3}$. One can combine this connection
with the one for $q_{12}$. (One could even add a new variable $R_{123}=\left(R_{1}^{2},R_{2}^{2},R_{3}^{2}\right)$
to the second bunch, for $c_{2}$, but a connection consisting of
a single bunch never affects contextuality analysis and can be dropped.)
The approach to contextuality presented in this paper allows for any
such modifications. However, we do not consider them obligatory, and
in some cases, as in the system $\mathcal{A}$ of Fig. \ref{fig:system 2by2},
they may be considered too restrictive (see the discussion in Section
\ref{sub: (Non)contextuality-of-consistent}).

\subsection{\label{sub: The-intuition-of}The intuition for (non)contextuality}

The main idea can be intuitively presented as follows. We have defined
a \mbox{c-c} system as a set of conteXt-representing bunches together
with connections between these bunches that reflect commonality of
conteNts. It is equally possible, however, to view a \mbox{c-c} system
as a set of conteNt-representing connections related to each other
by bunches that reflect commonality of conteXts. Thus, the system
$\mathcal{B}$ depicted in Fig.~\ref{fig: system general} consists
of the three connections, the elements of each of which are pairwise
stochastically unrelated random variables. However, the element $R_{1}^{1}$
of the first connection is stochastically related to the element $R_{2}^{1}$
of the second connection because they share a conteXt; and analogously
for the other two conteXts. 

We distinguish two forms of the dependence of random variables on
their conteXts. One of them is ``\emph{contextuality proper},''
the other one we call ``\emph{direct influenc}e\emph{s}.'' Let us
begin with the latter. Direct influences are reflected in the differences,
if any, between the distributions of the elements of the same connection.
For instance, if $R_{1}^{1}$ and $R_{1}^{2}$ have different distributions,
then the change of the conteXt from $c_{1}$ to $c_{2}$ directly
influences the random variable representing the conteNt $q_{1}$.
Direct influences are important, but they are of the same nature as
the dependence of random variables on their conteNt. Thus, if the
distribution of responses to the question about bees changes depending
on what movie has been previously watched, the influence of the movie
on the response is not any more puzzling than the influence of the
question itself. We prefer not to use the term ``contextuality''
for such forms of conteXt-dependence. Instead we describe them by
saying that the system of random variables is \emph{inconsistently
connected} if some of the elements of some of the connections in it
have different distributions. If in a system all elements of any connection
have the same distribution, then we call the system \emph{consistently
connected}. 

There is a simple and universal (applicable to all systems) way to
measure the degree of direct influences within each connection. Although
the connections in a \mbox{c-c} system consist of pairwise stochastically
unrelated random variables, nothing prevents one from thinking \emph{counterfactually}
of what their joint distribution could be \emph{if they were} jointly
distributed. (The reader should for now suspend criticism, as later
on we will define this counterfactual reasoning rigorously.)

Let us take, e.g., the connection $\left(R_{1}^{1},R_{1}^{2},R_{1}^{3}\right)$
in the system $\mathcal{B}$ depicted in Fig.~\ref{fig: system general}.
If these random variables were jointly distributed, then we would
be able to compute the probability with which they assume identical
values, $\Pr\left[R_{1}^{1}=R_{1}^{2}=R_{1}^{3}\right]$. One may
ask: among all possible ``imaginary'' joint distributions of $\left(R_{1}^{1},R_{1}^{2},R_{1}^{3}\right)$
(in which all three of them retain their individual distributions),
what is the maximal possible value for $\Pr\left[R_{1}^{1}=R_{1}^{2}=R_{1}^{3}\right]$?
As it turns out, this maximal probability is well-defined and uniquely
determinable: let us denote it $\max{}_{1}$. We define $\max{}_{2}$
for the connection $\left(R_{2}^{1},R_{2}^{2}\right)$ as the maximal
``imaginary'' value for $\Pr\left[R_{2}^{1}=R_{2}^{2}\right]$ given
the individual distributions of $R_{2}^{1}$ and $R_{2}^{2}$; and
we define $\max_{3}$ for the connection $\left(R_{3}^{2},R_{3}^{3}\right)$
analogously. These maximal probabilities of coincidence can be viewed
as reflecting the degree of direct influences of conteXts upon the
random variables: e.g., $\max_{1}=1$ if and only if the distributions
of all three random variables in $\left(R_{1}^{1},R_{1}^{2},R_{1}^{3}\right)$
are the same; and $\max{}_{1}=0$ if and only if the direct influences
are so prominent that the equality $R_{1}^{1}=R_{1}^{2}=R_{1}^{3}$
becomes ``unimaginable'': they cannot occur in any of the imagined
joint distributions because the supports of the three variable (the
subsets of possible values that have nonzero probability masses) do
not have elements in common. 

Now we come to ``contextuality proper.'' The maximal probabilities
just discussed are computed for each connection taken separately,
without taking into account the bunches that reflect the commonality
of conteXts across the connections. The question arises: are these
maximal probability values, 
\begin{equation}
\begin{array}{c}
\Pr\left[R_{1}^{1}=R_{1}^{2}=R_{1}^{3}\right]=\max{}_{1},\\
\\
\Pr\left[R_{2}^{1}=R_{2}^{2}\right]=\max{}_{2},\\
\\
\Pr\left[R_{3}^{2}=R_{3}^{3}\right]=\max{}_{3},
\end{array}
\end{equation}
compatible with the observed bunches of the system? In other words,
can one achieve these maximal (imaginary) probabilities in all three
connections simultaneously if one takes into account all the known
(not imagined) joint distributions in the bunches of the system? If
the answer is affirmative, then we can say that the knowledge of the
bunches representing different conteXts adds nothing to what we already
know of the direct influences by having considered the connections
separately --- we call such a system \emph{noncontextual}. If the
answer is negative, however, then the conteXts do influence the random
variables beyond any direct influences they exert on them --- the
system is \emph{contextual}. \footnote{In reference to footnote \ref{fn: As-the-reviewing} below, in the
newer version of CbD (Dzhafarov \& Kujala, 2016a,b), in the case of
more than two random variables, as in $\left(R_{1}^{1},R_{1}^{2},R_{1}^{3}\right)$,
the maximum probability should be considered not only for $R_{1}^{1}=R_{1}^{2}=R_{1}^{3}$
but also for each of $R_{1}^{1}=R_{1}^{2}$, $R_{1}^{2}=R_{1}^{3}$,
and $R_{1}^{1}=R_{1}^{3}$.}

\subsection{\label{sub: (Non)contextuality-of-consistent}(Non)contextuality
of consistently connected systems}

A system that exhibits no direct influences at all (i.e., is consistently
connected) may very well be contextual. In a consistently connected
version of our system $\mathcal{B}$ the three maximal probability
values will all be 1, and a system will be contextual if the ``imaginary''
equalities 
\begin{equation}
\begin{array}{c}
\Pr\left[R_{1}^{1}=R_{1}^{2}=R_{1}^{3}\right]=1,\\
\\
\Pr\left[R_{2}^{1}=R_{2}^{2}\right]=1,\\
\\
\Pr\left[R_{3}^{2}=R_{3}^{3}\right]=1,
\end{array}\label{eq: B consistent}
\end{equation}
are incompatible with the observed bunches of the system. This will
be the case when one can say, by abuse of language, that the system
is contextual because the elements in each of its connections cannot
be viewed as being essentially one and the same random variable. 

Consistent connectedness, in special forms, is known under a variety
of other names. In psychology, within the framework of so-called \emph{selective
influences}, the term describing consistent connectedness is ``marginal
selectivity'' (Townsend \& Schweickert, 1989). In quantum physics
it is often called ``no-signaling'' property, especially when dealing
with the EPR-type paradigms discussed in Section \ref{sec: Cyclic-c-c-systems}
(Popescu \& Rohrlich, 1994; Masanes, Acin, \& Gisin, 2006), or somewhat
more generally, ``no-disturbance'' property (Kurzynski, Cabello,
\& Kaszlikowski, 2014). Cereceda (2000) lists several other terms. 

Many scholars, especially in quantum mechanics, have considered contextuality
for consistently connected systems only (e.g., Dzhafarov \& Kujala,
2014c; Fine, 1982; Kurzynski, Ramanathan, Kaszlikowski, 2012; Kurzynski
et al., 2014). The same is true for the contextuality theory of Abramsky
and colleagues when it is applied to systems of random variables (Abramsky
\& Brandenburger, 2011; Abramsky et al., 2015). As a rule, however,
consistent connectedness is considered in a \emph{strong version},
wherein a consistently connected system should satisfy the following
property: in any two bunches $R^{1},R^{2}$ that share a set of conteNts
$q_{1},\ldots,q_{k}$, the corresponding sets of random variables
$\left(R_{1}^{1},\ldots,R_{k}^{1}\right)$ and $\left(R_{1}^{2},\ldots,R_{k}^{2}\right)$
have one and the same joint distribution. In the theory of selective
influences this property, or requirement, is called ``\emph{complete
marginal selectivity}'' (Dzhafarov, 2003). 

When applied to the system $\mathcal{B}$ in Fig.~\ref{fig: system general},
the strong form of consistent connectedness means that, in addition
to the same distribution of the random variables in each of the three
connections of $\mathcal{B}$, we also posit the same distribution
for $\left(R_{1}^{1},R_{2}^{1}\right)$ and $\left(R_{1}^{2},R_{2}^{2}\right)$
and for $\left(R_{1}^{2},R_{3}^{2}\right)$ and $\left(R_{1}^{3},R_{3}^{3}\right)$.
It is easy to see that this amounts to replacing the system $\mathcal{B}$
with the redefined system $\mathcal{B}''$ shown in Fig.~\ref{fig: alternative1-2},
and assuming that it is consistently connected in the ``ordinary''
sense. The CbD theory allows for both $\mathcal{B}$ and $\mathcal{B}''$
to represent one and the same empirical situation, the choice between
them being outside the scope of the theory. Therefore the notion of
consistent connectedness in this paper includes the strong version
thereof as a special case. 

\begin{figure}
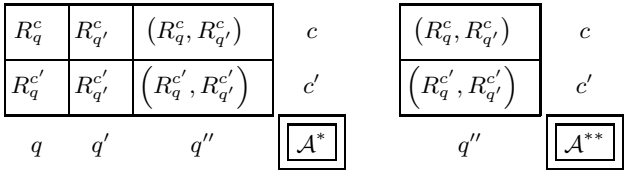

\begin{centering}
\begin{tabular}{|c|c|c|c}
\cline{1-3} 
$R_{q}^{c}$$\begin{array}{cccc}
\\
\\
\end{array}$ & $R_{q'}^{c}$$\begin{array}{cccc}
\\
\\
\end{array}$ & $\left(R_{q}^{c},R_{q'}^{c}\right)$$\begin{array}{cc}
\\
\\
\end{array}$ & $c$\tabularnewline
\cline{1-3} 
$R_{q}^{c'}$$\begin{array}{cccc}
\\
\\
\end{array}$ & $R_{q'}^{c'}$$\begin{array}{cccc}
\\
\\
\end{array}$ & $\left(R_{q}^{c'},R_{q'}^{c'}\right)$$\begin{array}{cc}
\\
\\
\end{array}$ & $c'$\tabularnewline
\cline{1-3} 
\multicolumn{1}{c}{$q$} & \multicolumn{1}{c}{$q'$} & \multicolumn{1}{c}{$q''$} & $\boxed{\boxed{\mathcal{A}^{*}}}$\tabularnewline
\end{tabular}$\qquad$%
\begin{tabular}{|c|c}
\cline{1-1} 
$\left(R_{q}^{c},R_{q'}^{c}\right)$$\begin{array}{cc}
\\
\\
\end{array}$ & $c$\tabularnewline
\cline{1-1} 
$\left(R_{q}^{c'},R_{q'}^{c'}\right)$$\begin{array}{cc}
\\
\\
\end{array}$ & $c'$\tabularnewline
\cline{1-1} 
\multicolumn{1}{c}{$q''$} & $\boxed{\boxed{\mathcal{A}^{**}}}$\tabularnewline
\end{tabular}
\par\end{centering}

\caption{\label{fig:system 2by2 strong}Two modifications of system $\mathcal{A}$
in our opening example that are required to consider this original
system consistently connected in the strong sense (``non-signaling,''
or ``complete marginal selectivity''). The two modifications are
equivalent if they are consistently connected in the ``ordinary''
sense.}
 
\end{figure}

The difference between the strong and weaker forms of consistent connectedness
is especially transparent if we consider the system $\mathcal{A}$
of our opening example (Fig.~\ref{fig:system 2by2}). Its consistent
connectedness means that the distribution of responses to a given
question, $q$ or $q'$, is the same irrespective of the conteXt,
$c$ or $c'$. The correlations between the two responses, however,
may very well be different in the two conteXts. If this is the case,
the consistently connected system $\mathcal{A}$ can be shown to be
contextual (see Section \ref{sec: Cyclic-c-c-systems}). By contrast,
if one assumes the strong form of consistent connectedness, the system
$\mathcal{A}$ is replaced with the system $\mathcal{A}^{*}$ shown
in Fig.~\ref{fig:system 2by2 strong}, consistently connected in
the ``ordinary'' sense. This system is trivially noncontextual,
as its two bunches have the same distribution. In Fig.~\ref{fig:system 2by2 strong}
this system is shown together with the system $\mathcal{A}^{**}$
in which the first two columns of the \mbox{c-c} matrix representing
$\mathcal{A}^{*}$ are dropped as redundant. Note, however, that the
systems $\mathcal{A}^{*}$ and $\mathcal{A}^{**}$ are not equivalent
if they are not consistently connected: the single-connection system
$\mathcal{A}^{**}$, as should be clear from Section \ref{sub: The-intuition-of},
is always noncontextual, whereas the system $\mathcal{A}^{*}$ may
very well be contextual.

\section{ConteXts and conteNts: a formal treatment}

Here, we present the basic conceptual set-up of our theory: a random
variable (confined to categorical random variables), jointly distributed
random variables (confined to finite sets thereof), functions of random
variables, and systems of random variables, with bunches and connections.
The reader who is not interested in a systematic introduction may
just skim through Sections \ref{sub:Base-sets-of} and \ref{sub: Systems-of-random}
and proceed to Section \ref{sec: Contextuality-analysis}.

Our view of random variables and relations among them is ``discourse-relative,''
in the sense that the existence of these variables and relations depends
on what other random variables are ``in play.''

\subsection{\label{sub:(Categorical)-random-variables}Categorical random variables}

We begin with a class $\mathsf{E}$ of (categorical) random variables
that we consider ``existing'' (or ``defined,'' or ``introducible,''
etc.). We need not be concerned with the cardinality of $\E$ as in
this paper we will always deal with finite subsets thereof.\footnote{The cardinality need not even be defined, as we consider $\E$ a class
rather than a set.} A random variable $X$ is a pair
\begin{equation}
X=\left(\id X,\di X\right),
\end{equation}
where $\id X$ is its unique \emph{identity label} (within the class
$\E$), whereas $\di X$ (to be read as a single symbol) is its \emph{distribution}.
The latter in turn is defined as a function
\begin{equation}
\di X:V_{X}\rightarrow[0,1],
\end{equation}
where $V_{X}$ is a \emph{finite set} (called the \emph{set of possible
values} of the random variable $X$), and
\begin{equation}
\sum_{v\in V_{X}}\di X\left(v\right)=1.
\end{equation}
The value $\di X\left(v\right)$ for any $v\in V_{X}$ is referred
to as the \emph{probability mass} of $X$ at its value $v$. For any
subset $W$ of $V_{X}$ we define the probability of $X\in W$ as
\begin{equation}
\Pr\left[X\in W\right]=\sum_{v\in W}\di X\left(v\right).
\end{equation}
In particular, for $v\in V_{X}$,
\begin{equation}
\Pr\left[X\in\left\{ v\right\} \right]=\di X\left(v\right),
\end{equation}
and we may also write $\Pr\left[X=v\right]$ instead of $\Pr\left[X\in\left\{ v\right\} \right]$.

Note that we impose no restrictions on the nature of the values $v$,
only that their set $V_{X}$ is finite. In particular, if $V_{1},\ldots,V_{n}$
are finite sets, then a random variable $Z\in\E$ with a distribution
\begin{equation}
\di Z:V_{1}\times\ldots\times V_{n}\rightarrow\left[0,1\right]
\end{equation}
is a categorical random variable. It can be denoted $Z=\left(X_{1},\ldots,X_{n}\right)$,
where $X_{i}$ is called the $i$th component (or the $i$th 1-marginal)
of $Z$, with the distribution defined by
\begin{equation}
\sum_{\substack{\left(v_{1},\ldots,v_{i},\ldots,v_{n}\right)\\
\in V_{1}\times\ldots\times V_{i-1}\\
\times\{v_{i}\}\\
\times V_{i+1}\times\ldots\times V_{n}
}
}\di Z\left(v_{1},\ldots,v_{n}\right)=\di X_{i}\left(v_{i}\right),\label{eq: 1-marginals}
\end{equation}
for any $v_{i}\in V_{i}$. The summation in this formula is across
all possible $n$-tuples $\left(v_{1},\ldots,v_{n}\right)$ with the
value of $v_{i}$ being fixed.
\begin{defn}
\label{def: Jointly distributed}We will say that $X_{1},\ldots,X_{n}$
in $\E$ are \emph{jointly distributed} if they are 1-marginals of
some $Z=\left(X_{1},\ldots,X_{n}\right)$ in $\E$. The random variable
$Z$ then can be called a vector (sequence, $n$-tuple), of \emph{jointly
distributed} $X_{1},\ldots,X_{n}$. If $X_{1},\ldots,X_{n}$ are not
jointly distributed, they are \emph{stochastically unrelated} (in
$\E$). 
\end{defn}
Note that according to this definition, $X_{1},\ldots,X_{n}$ in $\E$
are not jointly distributed if $\E$ does not contain $\left(X_{1},\ldots,X_{n}\right)$,
even though one can always conceive of a joint distribution for them.
This reflects our interpretation of $\E$ as the class of the variables
that ``exist'' (rather than just ``imagined,'' as discussed in
Section \ref{sub: The-intuition-of}). 

For any subsequence $\left(i_{1},\ldots,i_{k}\right)$ of $\left(1,\ldots,n\right)$
one can compute the corresponding \emph{$k$-marginal} of $Z$. Without
loss of generality, let $\left(i_{1},\ldots,i_{k}\right)=\left(1,\ldots,k\right)$.
Then the $k$-marginal $Y=\left(X_{1},\ldots,X_{k}\right)$ has the
distribution defined by
\begin{equation}
\sum_{\substack{\left(v_{1},\ldots,v_{k},v_{k+1}\ldots,v_{n}\right)\\
\in\left\{ v_{1}\right\} \times\ldots\times\left\{ v_{k}\right\} \\
\times V_{k+1}\ldots\times V_{X_{n}}
}
}\di Z\left(v_{1},\ldots,v_{n}\right)=\di Y\left(v_{1},\ldots,v_{k}\right),\label{eq: k-marginals}
\end{equation}
for any $\left(v_{1},\ldots,v_{k}\right)\in V_{1}\times\ldots\times V_{k}$.
The summation in this formula is across all possible $n$-tuples $\left(v_{1},\ldots,v_{n}\right)$
with the values of $v_{1},\ldots,v_{k}$ being fixed. This distribution
of the $k$-marginal $Y$ is referred to as a \emph{$k$-marginal
distribution}.

\subsection{\label{sub: Functions-of-random}Functions of random variables}

Let $X\in\E$ be a random variable with the distribution $\di X:V_{X}\rightarrow[0,1]$,
and let $f:V_{X}\rightarrow f\left(V_{X}\right)$ be some function.
The function $f\left(X\right)$ of a random variable $X$ is a random
variable $Y$ such that $X$ and $Y$ are 1-marginals of some random
variable $Z=\left(X,Y\right)$ with the distribution $\di Z:V_{X}\times f\left(V_{X}\right)\rightarrow\left[0,1\right]$
defined by
\begin{equation}
\di Z\left(v,w\right)=\left\{ \begin{array}{cc}
\di X\left(v\right) & \textnormal{if }w=f\left(v\right)\\
\\
0 & \textnormal{if otherwise}
\end{array}.\right.\label{eq: function f pair}
\end{equation}
It follows that the distribution of $Y$ as a 1-marginal of $Z$ is
defined by
\begin{equation}
\di Y\left(w\right)=\sum_{v\in f^{-1}\left(\left\{ w\right\} \right)}\di X\left(v\right),\label{eq: function f marginal}
\end{equation}
for any $w\in f\left(V_{X}\right)$. 

We stipulate as the main property of the class $\E$ that \emph{any
function of $X$ in $\E$ belongs to $\E$}. This property together
with the definition of 1-marginals implies that $Z=\left(X,f\left(X\right)\right)$
belongs to $\E$, i.e., $X$ and $f\left(X\right)$ are jointly distributed.

If $Y_{1}=f_{1}\left(X\right)$ and $Y_{2}=f_{2}\left(X\right)$,
we can consider $\left(f_{1},f_{2}\right)$ as a function $f$ mapping
$V_{X}$ into $f_{1}\left(V_{X}\right)\times f_{2}\left(V_{X}\right)$.
Then $Y=\left(Y_{1},Y_{2}\right)$ being a function of $X$ is merely
a special case of the situation considered above. Its meaning is that
$X$ and $\left(Y_{1},Y_{2}\right)$ are 1-marginals of some random
variable $Z=\left(X,\left(Y_{1},Y_{2}\right)\right)$ (that belongs
to $\E$) whose distribution is defined by
\begin{equation}
\di Z\left(v,\left(w_{1},w_{2}\right)\right)=\left\{ \begin{array}{cl}
\di X\left(v\right) & \textnormal{if }\left(w_{1},w_{2}\right)\\
 & =f\left(v\right)=\left(f_{1}\left(v\right),f_{2}\left(v\right)\right)\\
\\
0 & \textnormal{if otherwise}
\end{array}.\right.\label{eq: function f1,f2 pair}
\end{equation}
The (1-marginal) distribution of $Y=\left(Y_{1},Y_{2}\right)$ is
defined by
\begin{equation}
\begin{array}{l}
\di Y\left(w_{1},w_{2}\right)=\sum_{v\in f^{-1}\left(\left\{ \left(w_{1},w_{2}\right)\right\} \right)}\di X\left(v\right)\\
\\
=\sum_{v\in f_{1}^{-1}\left(\left\{ w_{1}\right\} \right)\cap f_{2}^{-1}\left(\left\{ w_{2}\right\} \right)}\di X\left(v\right).
\end{array}\label{eq: function f1,f2 marginal}
\end{equation}
The random variables $Y_{1}$ and $Y_{2}$ themselves are 1-marginals
of $Y=\left(Y_{1},Y_{2}\right)$ just defined. Indeed, the separate
distribution of $Y_{1}$ computed in accordance with (\ref{eq: function f marginal})
is
\begin{equation}
\di Y_{1}\left(w_{1}\right)=\sum_{v\in f_{1}^{-1}\left(\left\{ w_{1}\right\} \right)}\di X\left(v\right),
\end{equation}
for any $w_{1}\in f_{1}\left(V_{X}\right)$. We get the same formula
from (\ref{eq: function f1,f2 marginal}) by applying to it the formula
for computing 1-marginals, (\ref{eq: 1-marginals}): 
\begin{equation}
\begin{array}{l}
\di Y_{1}\left(w_{1}\right)=\sum_{w_{2}}\sum_{v\in f_{1}^{-1}\left(\left\{ w_{1}\right\} \right)\cap f_{2}^{-1}\left(\left\{ w_{2}\right\} \right)}\di X\left(v\right)\\
\\
=\sum_{v\in f_{1}^{-1}\left(\left\{ w_{1}\right\} \right)\cap\left(\bigcup_{w_{2}}f_{2}^{-1}\left(\left\{ w_{2}\right\} \right)\right)}\di X\left(v\right)\\
\\
=\sum_{v\in f_{1}^{-1}\left(\left\{ w_{1}\right\} \right)}\di X\left(v\right),
\end{array},
\end{equation}
because the union of the sets $f_{2}^{-1}\left(\left\{ w_{2}\right\} \right)$
across all values of $w_{2}\in f_{2}\left(V_{X}\right)$ is the entire
set $V_{X}$. Analogous reasoning applies to $Y_{2}$.

This shows that $Y_{1}$ and $Y_{2}$ defined as functions of some
$X\in\E$ are jointly distributed in the sense of Definition \ref{def: Jointly distributed}:
they are 1-marginals of some $Y=\left(Y_{1},Y_{2}\right)\in\E$. It
is easy to show that the converse holds true as well: if $Y_{1}$
and $Y_{2}$ in $\E$ are jointly distributed, then they are functions
of one and the same random variable that belongs to $\E$. Indeed,
in accordance with Definition \ref{def: Jointly distributed}, they
are 1-marginals of some $Y=\left(Y_{1},Y_{2}\right)\in\E$. But then
$Y_{1}$ and $Y_{2}$ are functions of this $Y$. Specifically, denoting
the sets of possible values for $Y_{1},Y_{2}$ by $W_{1},W_{2}$,
respectively, we have $Y_{1}=f_{1}\left(Y\right)$, where 
\begin{equation}
f_{1}:W_{1}\times W_{2}\rightarrow W_{1}
\end{equation}
is defined by $f_{1}\left(w_{1},w_{2}\right)=w_{1}$ (a projection
function). The computations of the distribution of $Y_{1}$ in accordance
with (\ref{eq: 1-marginals}) coincides with that in accordance with
(\ref{eq: function f marginal}),
\begin{equation}
\di Y_{1}\left(w_{1}\right)=\sum_{w_{2}\in W_{2}}\di Y\left(w_{1},w_{2}\right)=\sum_{v\in f_{1}^{-1}\left(\left\{ w_{1}\right\} \right)}\di Y\left(v\right),
\end{equation}
for any $w_{1}\in W_{1}$. Analogous reasoning applies to $Y_{2}$.

This result is trivially generalized to an arbitrary finite set of
random variables.\footnote{In fact it holds for any set of any random entities (Dzhafarov \&
Kujala, 2010), but our focus in this paper is on finite sets of categorical
random variables.}
\begin{thm}
\label{thm: JDC}Random variables $X_{1},\ldots,X_{n}\in\E$ are jointly
distributed if and only if they are representable as functions of
one and the same random variable $X\in\E$.
\end{thm}
Note that this $X$ may very well equal one of the $X_{1},\ldots,X_{n}$.
More generally, one can add $X$ to its functions $X_{1},\ldots,X_{n}$
to create a jointly distributed set $X,X_{1},\ldots,X_{n}$ in which
all elements are, obviously, functions of one of its elements. Note
also that if $X_{i}=f_{i}\left(X\right)$, $i=1,\ldots,n$, then $f=\left(f_{1},\ldots,f_{n}\right)$
is a function, and we can equivalently reformulate Theorem \ref{thm: JDC}
as saying that a vector of random variables $\left(X_{1},\ldots,X_{n}\right)$
is a random variable (by definition, with jointly distributed components)
if and only if it is a function of some random variable $X$.

In spite of its simplicity, Theorem \ref{thm: JDC} was discovered,
in various special forms, only in the 1980s (Suppes \& Zanotti, 1981;
Fine, 1982). It has a direct bearing on the problem of ``hidden variables''
in quantum mechanics: given a set of random variables, is there a
random entity of which these random variables are functions? To formulate
this problem rigorously and to enable the use of Theorem \ref{thm: JDC}
for solving it we will need the notion of a coupling, introduced below
(Section \ref{sub: (Probabilistic)-couplings}).

\subsection{\label{sub: Two-meanings-of}Two meanings of equality of random variables}

The following remark may prevent possible confusions. Given a random
variable $Z$ and a measurable set $E$, the expression $Z\in E$
clearly does not mean that $Z$ as a random variables (with its identity
$\id Z$ and distribution $\di Z$) is an element of $E$. Rather
this expression is a way of saying that we are considering an event
$E$ in the measure space $\di Z=\left(S_{Z},\Sigma_{Z},\mu_{Z}\right)$
associated with a random variable $Z$. Thus, $\Pr\left[Z\in E\right]$
is $\mu_{Z}\left(E\right)$. As a special case, given jointly distributed
$X,Y$, the expression $X=Y$ is merely a shortcut for $\left(X,Y\right)\in W$,
where $W=\left\{ \left(v_{1},v_{2}\right)\in S_{X}\times S_{Y}:v_{1}=v_{2}\right\} $.
This meaning of equality should not be confused with another meaning:
$X=Y$ can also mean that these two symbols refer to one and the same
random variable, so that $\id X=\id Y$ and (consequently) $\di X=\di Y$.
We think that the meaning of $X=Y$ in this paper is always clear
from the context (now using this word without capital X).

\subsection{\label{sub:Base-sets-of}Base sets of random variables\label{sub: Base set}}

How does one construct the class $\E$? For instance, with $X_{1},\ldots,X_{n}$
all in $\E$, how do we know whether they are jointly distributed,
i.e., whether $\E$ contains a $Z=\left(X_{1},\ldots,X_{n}\right)$?
Can we simply declare that any random variables $X_{1},\ldots,X_{n}$
are jointly distributed? The answer to the last question is negative:
to be able to model empirical phenomena one needs to keep the meaning
of joint distribution tied to the empirical meaning of ``co-occurrence''
--- which means that joint distribution cannot be imposed arbitrarily. 

To make all of this clear, let us construct the class $\E$ of ``existing''
random variables systematically. The construction is simple: we introduce
a nonempty \emph{base set} $\mathsf{\R}$ of (categorical) random
variables (in this paper we assume this set to be finite, but this
need not be so generally), and we posit that
\begin{description}
\item [{(P1)}] a random variable belongs to $\E$ if and only if it is
a function of any one of the elements of $\mathsf{\R}$;
\item [{(P2)}] no random variable in $\E$ is a function of two distinct
elements of $\mathsf{\R}$.
\end{description}
The constraints (P1-P2) ensure that no two random variables existing
in the sense of belonging to $\E$ may have a joint distribution unless
they are functions of one and the same element of $\mathsf{\R}$.
Indeed, let some transformations $\alpha\left(A\right)$ and $\beta\left(B\right)$
have a joint distribution, for $A,B\in\mathsf{R}$. Then a random
variable $\left(\alpha\left(A\right),\beta\left(B\right)\right)$
exists, which means that this pair is a function of some $C\in\mathsf{R}$.
But then $\alpha\left(A\right)$ is a function of $A$ and $C$, whence
$A=C$, and $\beta\left(B\right)$ is a function of both $B$ and
$C$, whence $B=C=A$. (In reference to Section \ref{sub: Two-meanings-of},
the equalities here are used in the sense of ``one and the same random
variable.'') 

Instead of ``$X$ belongs to $\E$'' we can also say ``$X$ exists
with respect to $\R$.'' This is preferable if one deals with different
base sets $\R$ inducing different classes $\E$, as we do in the
subsequent sections.

Consider an example: let $\mathsf{R}$ consist of the four random
variables 
\[
X=\left(X_{1},X_{2},X_{3}\right),Y=\left(Y_{1},Y_{2}\right),Z,U=\left(U_{1},U_{2}\right).
\]
These random variables are declared to exist, and then so are functions
of these random variables. Thus, $X_{2}$ exists because $X$ exists
and $X_{2}$ is its function (second projection). Analogously, if
the values of $Y_{1},Y_{2}$ are numerical, the variable $Y_{1}+Y_{2}$
exists. However, no component of one of the four random variables,
say, $X_{2}$, is jointly distributed with any component of another,
say, $U_{1}$, and no function $f\left(U_{1},X_{2}\right)$ is a random
variable (its distribution is undefined). By the same logic, no two
different vectors in $\R$ can share a component: if they did, this
component would be a function of both of them, contravening (P2).

\subsection{\label{sub: Systems-of-random}Systems of random variables}

The example of $\R$ at the end of the previous section is in fact
how we introduce our main object: \emph{conteXt-conteNt systems of
random variables}.
\begin{defn}
\label{def: c-c}Let $\R$ be a base set of (categorical) random variables
each element of which, called a \emph{bunch}, is a vector of random
variables. Let $U_{\R}$ be the union of all components of these bunches.
A \emph{conteXt-conteNt (c-c) system }$\mathcal{R}$ of random variables
based on $\mathsf{\R}$ is created by endowing $\mathsf{\R}$ with
a partition of $U_{\R}$ into subsets called \emph{connections} and
satisfying the following two properties: 

\emph{(intersection property}) a bunch and a connection do not have
more than one component of $U_{\R}$ in common; and

\emph{(comparability property) }elements of a connection have the
same set of possible values.\footnote{In a more general treatment this translates into the same set and
the same sigma algebra of events.}
\end{defn}
Recall that partitioning of a set means creating a set of pairwise
disjoint subsets whose union is the entire set. Note that due to the
intersection property in Definition \ref{def: c-c}, any two elements
of a connection are stochastically unrelated (they are 1-marginals
of different bunches). Note also that due to the comparability property
the elements of a connection may (but generally do not) have the same
distribution. 

Let the bunches of the \mbox{c-c} system be enumerated $1,\ldots,n$,
and the connections be enumerated $1,\ldots,m$. Due to the intersection
property in Definition \ref{def: c-c}, any random variable in the
set $U_{\R}$ of a \mbox{c-c} system $\mathcal{R}$ can be uniquely
identified by the labels of the bunch and of the connection it belongs
to. These labels (or some symbols in a one-to-one correspondence with
them) are referred to as \emph{conteXts} (labels for bunches) and
\emph{conteNts} (labels for connections). As we see, in the formal
theory bunches and connection define rather than are defined by the
conteXts and conteNts, respectively. It is the other way around in
empirical applications (see the introductory section), where our understanding
of what constitutes a given conteNt under different conteXts guides
the creation of the bunches and connections.

The unique labeling of the random variables by the conteXts and conteNts
means that any \mbox{c-c} system can be presented in the form already
familiar to us from the introduction: a conteXt-conteNt (\mbox{c-c})
matrix. An example of a \mbox{c-c} system presented in the form of
a \mbox{c-c} matrix is given in Fig.~\ref{fig: system general}.
The initial base set of random variables is
\begin{equation}
\mathsf{\R}=\left\{ \begin{array}{l}
R^{1}=\left(R_{1}^{1},R_{2}^{1}\right)\\
\\
R^{2}=\left(R_{1}^{2},R_{2}^{2},R_{3}^{2}\right)\\
\\
R^{3}=\left(R_{1}^{3},R_{3}^{3}\right)
\end{array}\right\} ,\label{eq: three bunches in R}
\end{equation}
the union set is 
\begin{equation}
U_{\R}=\left\{ R_{1}^{1},R_{2}^{1},R_{1}^{2},R_{2}^{2},R_{3}^{2},R_{1}^{3},R_{3}^{3}\right\} \label{eq: union of all in R}
\end{equation}
with the lower indexes already chosen in view of the partitioning
into connections,
\begin{equation}
\left\{ \begin{array}{l}
\left(R_{1}^{1},R_{1}^{2},R_{1}^{3}\right)\\
\\
\left(R_{2}^{1},R_{2}^{2}\right)\\
\\
\left(R_{3}^{2},R_{3}^{3}\right)
\end{array}\right\} .\label{eq: three connections in R}
\end{equation}

The intersection property in Definition \ref{def: c-c} is critical:
if a conteXt and a connection could have more than one random variable
in common, both the double-indexing of the random variables by conteXts
and conteNts and the subsequent contextuality analysis of the system
would be impossible.

\subsection{Kolmogorovian Probability Theory and Contextuality-by-Default}

In this section we briefly discuss the relationship between KPT and
CbD. This discussion is not needed for understanding the subsequent
sections. We will assume the reader's familiarity with the basics
of measure theory.

The definition of (categorical) random variables in KPT is as follows.
Let $\left(S,\Sigma,\mu\right)$ be a domain probability space, and
let $\left(V_{X},\Sigma_{X}\right)$ be a codomain measurable space,
with $V_{X}$ a finite set and $\Sigma_{X}$ usually (and here) defined
as its power set. A random variable $X$ is a function $S\rightarrow V_{X}$
such that $X^{-1}\left(\left\{ v\right\} \right)\in\Sigma$, for any
$v\in V_{X}$. The probability mass $p_{X}\left(v\right)$ is defined
as $\mu\left(X^{-1}\left(\left\{ v\right\} \right)\right)$, and for
any subset $V\subset V_{X},$ the probability of $X$ falling in $V$
is computed as 
\begin{equation}
\Pr\left[X\in V\right]=\mu\left(X^{-1}\left(V\right)\right)=\sum_{v\in V}p_{X}\left(v\right).
\end{equation}
We call $\left(S,\Sigma,\mu\right)$ the \emph{sample space}\footnote{It seems common to use this term for the set $S$ alone; but the term
``space'' in mathematics means a set with some structure imposed
on it, and the structure here is the sigma algebra and the measure.
We prefer therefore to use the term ``sample space'' for the entire
domain probability space. $S$ alone can be referred to as the \emph{sample
set} for $X$. } for $X$. 

The great conceptual convenience of KPT is that the joint distribution
of two random variables taken as two functions defined on the same
sample space is uniquely determined by these two functions: if $X$
is as above and $Y$ is another random variable, then its joint distribution
with $X$ above is defined by
\begin{equation}
p_{XY}\left(v,w\right)=\mu\left(X^{-1}\left(\left\{ v\right\} \right)\cap Y^{-1}\left(\left\{ w\right\} \right)\right),
\end{equation}
for any $\left(v,w\right)\in V_{X}\times V_{Y}$. 

In CbD, random variables are considered only with respect to a specified
base set. A random variable exists if it is a function of one and
only one of the elements of this base set; and functions of different
base random variables are considered stochastically unrelated. Is
this picture compatible with KPT? We think it is, provided KPT is
not naively thought of as positing the existence of a single sample
space for all imaginable random variables. Such a view can be shown
to be mathematically flawed (Dzhafarov \& Kujala, 2014a-b).

Every sample space $\left(S,\Sigma,\mu\right)$ corresponds to a random
variable $Z$ defined as the identity mapping $S\rightarrow S$ from
$\left(S,\Sigma,\mu\right)$ to $\left(S,\Sigma\right)$, and then
every random variable defined on $\left(S,\Sigma,\mu\right)$ is representable
as a transformation of $Z$. If we consider a set of sample spaces
unrelated to each other, then the corresponding identity functions
form a base set of random variables, and what we get is essentially
the same picture as in CbD. 

We need one qualification though: even if all the functions considered
are categorial random variables, the base set itself need not be a
finite set of categorical random variables, as it is in Section \ref{sub: Systems-of-random}.
This is not, however, a restriction inherent in CbD but the choice
we have made in this paper. A finite number of categorical base variables
are sufficient if one only considers a finite set of functions thereof,
which is the case we deal with.

\section{\label{sec: Contextuality-analysis}Contextuality analysis}

In this section we give the definitions and introduce the conceptual
apparatus involved in determining whether a \mbox{c-c} system is
contextual or noncontextual.

\subsection{\label{sub: (Probabilistic)-couplings}Probabilistic couplings}

Imagining joint distributions for things that are not jointly distributed,
as it was presented in Section \ref{sub: The-intuition-of}, is not
rigorous mathematics. The latter requires that we use the mathematical
tool of \emph{(probabilistic) couplings}.
\begin{defn}
A coupling of a set of random variables $X_{1},\ldots,X_{n}$ is a
random variable $\left(Y_{1},\ldots,Y_{n}\right)$ (with jointly distributed
components) such that $Y_{i}$ has the same distribution as $X_{i}$,
for all $i=1,\ldots,n$.
\end{defn}
As an illustration, let $X_{1}$ and $X_{2}$ be distributed as

\begin{center}
\begin{tabular}{c|c|c|c|c}
\multicolumn{1}{c}{} & \multicolumn{1}{c}{$X_{1}=1$} & \multicolumn{1}{c}{$X_{1}=2$} & \multicolumn{1}{c}{$X_{1}=3$} & \tabularnewline
\cline{2-4} 
$\textnormal{pr. mass}$ & 0.3 & 0.3 & 0.4 & \tabularnewline
\cline{2-4} 
\end{tabular}
\par\end{center}

\noindent and

\begin{center}
\begin{tabular}{c|c|c|c}
\multicolumn{1}{c}{} & \multicolumn{1}{c}{$X_{2}=1$} & \multicolumn{1}{c}{$X_{2}=2$} & \multirow{2}{*}{.}\tabularnewline
\cline{2-3} 
$\textnormal{pr. mass}$ & 0.7 & 0.3 & \tabularnewline
\cline{2-3} 
\end{tabular}
\par\end{center}

\noindent Then $\left(Y_{1},Y_{2}\right)$ with the distribution

\begin{center}
\begin{tabular}{c|c|c|c|c}
\multicolumn{1}{c}{} & \multicolumn{1}{c}{$Y_{1}=1$} & \multicolumn{1}{c}{$Y_{1}=2$} & \multicolumn{1}{c}{$Y_{1}=3$} & \tabularnewline
\cline{2-4} 
$Y_{2}=1$ & 0.3 & 0.2 & 0.2 & 0.7\tabularnewline
\cline{2-4} 
$Y_{2}=2$ & 0 & 0.1 & 0.2 & 0.3\tabularnewline
\cline{2-4} 
\multicolumn{1}{c}{} & \multicolumn{1}{c}{0.3} & \multicolumn{1}{c}{0.3} & \multicolumn{1}{c}{0.4} & \tabularnewline
\end{tabular}
\par\end{center}

\noindent is a coupling for $X_{1}$ and $X_{2}$. And so is $\left(Y'_{1},Y'_{2}\right)$
with the distribution

\begin{center}
\begin{tabular}{c|c|c|c|c}
\multicolumn{1}{c}{} & \multicolumn{1}{c}{$Y'_{1}=1$} & \multicolumn{1}{c}{$Y'_{1}=2$} & \multicolumn{1}{c}{$Y'_{1}=3$} & \tabularnewline
\cline{2-4} 
$Y'_{2}=1$ & 0.3 & 0 & 0.4 & 0.7\tabularnewline
\cline{2-4} 
$Y'_{2}=2$ & 0 & 0.3 & 0 & 0.3\tabularnewline
\cline{2-4} 
\multicolumn{1}{c}{} & \multicolumn{1}{c}{0.3} & \multicolumn{1}{c}{0.3} & \multicolumn{1}{c}{0.4} & \tabularnewline
\end{tabular}.
\par\end{center}

\noindent Generally, the number of couplings of a given set of random
variables is infinite. 

In our paper couplings are constructed in two ways only: either for
connections in a \mbox{c-c} system, taken separately, or for the
entire set of bunches in the system. In relation to the system $\mathcal{B}$
in Fig. \ref{fig: system general} and (\ref{eq: three bunches in R}),
a coupling for the connection $\left(R_{1}^{1},R_{1}^{2},R_{1}^{3}\right)$
is a triple $\left(T_{1}^{1},T_{1}^{2},T_{1}^{3}\right)$ such that
$T_{1}^{j}$ and $R_{1}^{j}$ have the same distribution, for $j=1,2,3$;
and analogously for the other two connections. 

The set of the three bunches $\left(R^{1},R^{2},R^{3}\right)$ in
(\ref{eq: three bunches in R}) is coupled by $S=\left(S^{1},S^{2},S^{3}\right)$
where 
\begin{equation}
\left\{ \begin{array}{l}
S^{1}=\left(S_{1}^{1},S_{2}^{1}\right)\\
\\
S^{2}=\left(S_{1}^{2},S_{2}^{2},S_{3}^{2}\right)\\
\\
S^{3}=\left(S_{1}^{3},S_{3}^{3}\right)
\end{array}\right\} 
\end{equation}
such that $S^{j}$ and $R^{j}$ have the same distribution, for $j=1,2,3$.

In the following we will freely use phrases indicating that a coupling
for some random variables ``exists,'' or ``can be constructed,''
or that these random variables ``can be coupled.'' Note, however,
that the couplings do not ``exist'' with respect to the base set
of random variables formed by the bunches of a \mbox{c-c} system,
as no coupling of the bunches can be presented as a function of just
one of these bunches. If the bunches are assumed to have links to
empirical observations, then the couplings can be said to have no
empirical meaning. A coupling forms a base set of its own, consisting
of itself. Its marginals (or \emph{subcouplings}) corresponding to
the bunches of the \mbox{c-c} system do ``exist'' with respect
to this new base set, as they are functions of its only element. However,
the bunches of the \mbox{c-c} system themselves do not ``exist''
with respect to the base set formed by this coupling. One can add
the coupling $S=\left(S^{1},S^{2},S^{3}\right)$ to the set $\left(R^{1},R^{2},R^{3}\right)$
of the three bunches of our system $\mathcal{B}$ as a fourth element
of a new base set, stochastically unrelated to the bunches.

\subsection{\label{sub: Flattening-convention}``Flattening'' convention}

Let us adopt the following simplifying convention in regard to couplings
(and more generally, vectors of jointly distributed variables): a
vector of jointly distributed random variables $\left(A^{1},\ldots,A^{n}\right)$
in which $A^{i}=\left(A_{1}^{i},\ldots,A_{k_{i}}^{i}\right)$, for
each $i=1,\ldots,n$, is considered equivalent (replaceable by) the
vector 
\begin{equation}
\left(A_{1}^{1},\ldots,A_{k_{1}}^{1},\ldots,A_{1}^{n},\ldots,A_{k_{n}}^{n}\right).
\end{equation}
As each of the random variables is assumed to be uniquely indexed
(in our analysis, double-indexed), the order in which they are shown
in any given vector is usually arbitrary. As an example, a coupling
$S=\left(S^{1},S^{2},S^{3}\right)$ of the three bunches $\mathcal{R}=\left(R^{1},R^{2},R^{3}\right)$
in (\ref{eq: three bunches in R}), written in extenso, is
\begin{equation}
S=\left(\left(S_{1}^{1},S_{2}^{1}\right),\left(S_{1}^{2},S_{2}^{2},S_{3}^{2}\right),\left(S_{1}^{3},S_{3}^{3}\right)\right)
\end{equation}
with $\left(S_{1}^{1},S_{2}^{1}\right)$ distributed as the bunch
$\left(R_{1}^{1},R_{2}^{1}\right)$, etc. In accordance with our agreement,
this coupling can be equivalently written as
\begin{equation}
S=\left(S_{1}^{1},S_{2}^{1},S_{1}^{2},S_{2}^{2},S_{3}^{2},S_{1}^{3},S_{3}^{3}\right),\label{eq: coupling for R}
\end{equation}
in which $\left(S_{1}^{1},S_{2}^{1}\right)$ distributed as the bunch
$\left(R_{1}^{1},R_{2}^{1}\right)$, etc. The ``flattening'' convention
makes it easier to compare couplings of a connection taken in isolation
with the \emph{subcoupling} of the coupling (\ref{eq: coupling for R})
corresponding to the same connection. Thus, for the connection $\left(R_{1}^{1},R_{1}^{2},R_{1}^{3}\right)$
taken separately we can consider all possible couplings $\left(T_{1}^{1},T_{1}^{2},T_{1}^{3}\right)$
and then compare them with the subcouplings $\left(S_{1}^{1},S_{1}^{2},S_{1}^{3}\right)$
extracted as 3-marginals from all possible couplings (\ref{eq: coupling for R}).
We will need such comparisons for determining whether the system in
question is contextual.

\subsection{Maximal couplings for connections}
\begin{defn}
Let $R_{j}^{1},\ldots,R_{j}^{k}$ be a connection (for a conteNt $q_{j}$)
in a \mbox{c-c} system. A coupling $\left(T_{j}^{1},\ldots,T_{j}^{k}\right)$
of $R_{j}^{1},\ldots,R_{j}^{k}$ is a \emph{maximal coupling} if the
value of 
\[
\Pr\left[T_{j}^{1}=\ldots=T_{j}^{k}\right]
\]
is the largest possible among all couplings of $R_{j}^{1},\ldots,R_{j}^{k}$.
\end{defn}
(In relation to Section \ref{sub: Two-meanings-of}, the equality
here clearly is not the identity of the random variables but a description
of an event associated with jointly distributed variables.) Theorem
\ref{thm: maximal-coupling} below ensures that the maximum mentioned
in the definition always exist. The notion of a maximal coupling is
well-defined for arbitrary sets of arbitrary random variables (see
Thorisson, 2000), but we will only need it for connections formed
by categorical variables.\footnote{It is a common mistake to think that this notion may be useless outside
the class of categorical variables: thus, one might erroneously assume
that if the distributions of $T_{j}^{1},\ldots,T_{j}^{k}$ are absolutely
continuous with respect to the Lebesgue measure on the set of reals,
then $\Pr\left[T_{j}^{1}=\ldots=T_{j}^{k}\right]$ must be zero. In
fact, this probability can be any number between 0 and 1, because
the joint distribution of $T_{j}^{1},\ldots,T_{j}^{k}$ within the
intersection of their supports may very well be concentrated on the
diagonal $T_{j}^{1}=\ldots=T_{j}^{k}$. In particular, if $R_{j}^{1},\ldots,R_{j}^{k}$
are identically distributed, then, irrespective of this distribution,
they have a maximal coupling $\left(T_{j}^{1},\ldots,T_{j}^{k}\right)$
with $\Pr\left[T_{j}^{1}=\ldots=T_{j}^{k}\right]=1$. } 

For any coupling $\left(T_{j}^{1},\ldots,T_{j}^{k}\right)$ of $\left(R_{j}^{1},\ldots,R_{j}^{k}\right)$
,

\begin{equation}
\Pr\left[T_{j}^{1}=\ldots=T_{j}^{k}\right]=\sum_{v\in V}\Pr\left[T_{j}^{1}=v,\ldots,T_{j}^{k}=v\right],
\end{equation}
where $V$ is the set of possible values shared by the elements of
the connection. This sum is maximized across all possible couplings
if each of the summands on the right-hand side is maximized separately.
The maximal possible value for $\Pr\left[T_{j}^{1}=v,\ldots T_{j}^{k}=v\right]$
(with the individual distributions of $Y_{i}$ being fixed) is
\begin{equation}
\max p_{v}=\min\left(\Pr\left[T_{j}^{1}=v\right],\ldots,\Pr\left[T_{j}^{k}=v\right]\right).\label{eq: max coupling mins}
\end{equation}
Indeed, the probability of a joint event can never exceed any of the
probabilities of the component events. To prove that a maximal coupling
exists for any connection, we need to show that every value of $\left(T_{j}^{1},\ldots,T_{j}^{k}\right)$
can be assigned a probability so that, for all $v\in V$, 
\begin{equation}
\Pr\left[T_{j}^{1}=v,\ldots T_{j}^{k}=v\right]=\max p_{v},
\end{equation}
and
\begin{equation}
\sum_{\substack{v_{1},\ldots,v_{i-1},\\
v_{i}=v,\\
v_{i+1},\ldots,v_{k}
}
}\Pr\left[T_{j}^{1}=v_{1},\ldots,T_{j}^{k}=v_{k}\right]=\Pr\left[Y_{i}=v\right],
\end{equation}
for any $i=1,\ldots,k$. A simple proof that this is always possible
can be found in Thorisson (2000, pp. 7-8 and 104-107).
\begin{thm}
\label{thm: maximal-coupling}A maximal coupling $\left(T_{j}^{1},\ldots,T_{j}^{k}\right)$
can be constructed for any connection $\left(R_{j}^{1},\ldots,R_{j}^{k}\right)$
in a \mbox{c-c} system, with
\begin{equation}
\begin{array}{l}
\Pr\left[T_{j}^{1}=\ldots=T_{j}^{k}=v\right]\\
=\min\left(\Pr\left[T_{j}^{1}=v\right],\ldots,\Pr\left[T_{j}^{k}=v\right]\right),
\end{array}\label{eq: theorem}
\end{equation}
for any $v$ in the set $V$ of possible values of (each of) $R_{j}^{1},\ldots,R_{j}^{k}$.
\end{thm}
As an example, let the variables $R_{1}^{1},R_{1}^{2},R_{1}^{3}$
in the first connection of the system $\mathcal{B}$ of Fig.~\ref{fig: system general}
be binary, with the possible values 1 and 2. Let

\begin{center}
\begin{tabular}{cccc}
 & $R_{1}^{1}=1$ & $R_{1}^{1}=2$ & \multirow{2}{*}{,}\tabularnewline
\cline{2-3} 
\multicolumn{1}{c|}{$\textnormal{pr. mass}$} & \multicolumn{1}{c|}{0.3} & \multicolumn{1}{c|}{0.7} & \tabularnewline
\cline{2-3} 
 &  &  & \tabularnewline
\end{tabular}
\par\end{center}

\begin{center}
\begin{tabular}{cccc}
 & $R_{1}^{2}=1$ & $R_{1}^{2}=2$ & \multirow{2}{*}{,}\tabularnewline
\cline{2-3} 
\multicolumn{1}{c|}{$\textnormal{pr. mass}$} & \multicolumn{1}{c|}{0.4} & \multicolumn{1}{c|}{0.6} & \tabularnewline
\cline{2-3} 
 &  &  & \tabularnewline
\end{tabular}
\par\end{center}

\begin{center}
\begin{tabular}{c|c|c|c}
\multicolumn{1}{c}{} & \multicolumn{1}{c}{$R_{1}^{3}=1$} & \multicolumn{1}{c}{$R_{1}^{3}=2$} & \multirow{2}{*}{.}\tabularnewline
\cline{2-3} 
$\textnormal{pr. mass}$ & 0.7 & 0.3 & \tabularnewline
\cline{2-3} 
\end{tabular}
\par\end{center}

\noindent Then, in the maximal coupling $\left(T_{1}^{1},T_{1}^{2},T_{1}^{3}\right)$,
\[
\begin{array}{l}
\max p_{1}=\Pr\left[T_{1}^{1}=T_{1}^{2}=T_{1}^{3}=1\right]\\
\\
=\min\left(0.3,0.4,0.7\right)=0.3,\\
\\
\max p_{2}=\Pr\left[T_{1}^{1}=T_{1}^{2}=T_{1}^{3}=2\right]\\
\\
=\min\left(0.7,0.6,0.3\right)=0.3.
\end{array},
\]
We can now assign probabilities to the rest of the values of $\left(T_{1}^{1},T_{1}^{2},T_{1}^{3}\right)$
in an infinity of possible ways, e.g., as shown below,

\begin{center}
\begin{tabular}{c|c|c|c|c|c|c|c|c|c}
\multicolumn{1}{c}{value:} & \multicolumn{1}{c}{111} & \multicolumn{1}{c}{112} & \multicolumn{1}{c}{121} & \multicolumn{1}{c}{211} & \multicolumn{1}{c}{221} & \multicolumn{1}{c}{212} & \multicolumn{1}{c}{122} & \multicolumn{1}{c}{222} & \multirow{2}{*}{,}\tabularnewline
\cline{2-9} 
$\textnormal{pr. mass}$ & 0.3 & 0 & 0 & 0.1 & 0.3 & 0 & 0 & 0.3 & \tabularnewline
\cline{2-9} 
\end{tabular}
\par\end{center}

\noindent so that these joint probabilities yield the right values
of $\Pr\left[T_{1}^{j}=1\right]=\Pr\left[R_{1}^{j}=1\right]$, for
$j=1,2,3$.

Consider another example, using the second connection in the system
$\mathcal{B}$, $\left(R_{2}^{1},R_{2}^{2}\right)$. Suppose that
both these variables have the same distribution:

\begin{center}
\begin{tabular}{c|c|c|c|c}
\multicolumn{1}{c}{value} & \multicolumn{1}{c}{$1$} & \multicolumn{1}{c}{$2$} & \multicolumn{1}{c}{3} & \multirow{2}{*}{,}\tabularnewline
\cline{2-4} 
$\textnormal{pr. mass}$ & 0.3 & 0.2 & 0.5 & \tabularnewline
\cline{2-4} 
\end{tabular}
\par\end{center}

\noindent in which case we can say that this connection is \emph{consistent}
(and, to remind, if this is the case for all connections, then the
system is consistently connected). The maximal coupling $\left(T_{2}^{1},T_{2}^{2}\right)$
here has a uniquely determined distribution

\begin{center}
\begin{tabular}{c|c|c|c|c}
\multicolumn{1}{c}{} & \multicolumn{1}{c}{$T_{2}^{1}=1$} & \multicolumn{1}{c}{$T_{2}^{1}=2$} & \multicolumn{1}{c}{$T_{2}^{1}=3$} & \tabularnewline
\cline{2-4} 
$T_{2}^{2}=1$ & 0.3 & 0 & 0 & 0.3\tabularnewline
\cline{2-4} 
$T_{2}^{2}=2$ & 0 & 0.2 & 0 & 0.2\tabularnewline
\cline{2-4} 
$T_{2}^{2}=3$ & 0 & 0 & 0.5 & 0.5\tabularnewline
\cline{2-4} 
\multicolumn{1}{c}{} & \multicolumn{1}{c}{0.3} & \multicolumn{1}{c}{0.2} & \multicolumn{1}{c}{0.5} & \tabularnewline
\end{tabular},
\par\end{center}

\noindent with $\Pr\left[T_{2}^{1}=T_{2}^{2}\right]=1$.

\subsection{Contextuality}

\begin{figure}
\begin{centering}
\begin{tabular}{|c|c|c|c}
\cline{1-3} 
$R_{1}^{1}$$\begin{array}{cc}
\\
\\
\end{array}$ & $R_{2}^{1}$$\begin{array}{cc}
\\
\\
\end{array}$ & $\cdot$$\begin{array}{cc}
\\
\\
\end{array}$ & $c_{1}$\tabularnewline
\cline{1-3} 
$\cdot$$\begin{array}{cc}
\\
\\
\end{array}$ & $R_{2}^{2}$$\begin{array}{cc}
\\
\\
\end{array}$ & $R_{3}^{2}$$\begin{array}{cc}
\\
\\
\end{array}$ & $c_{2}$\tabularnewline
\cline{1-3} 
$R_{1}^{3}$$\begin{array}{cc}
\\
\\
\end{array}$ & $\cdot$$\begin{array}{cc}
\\
\\
\end{array}$ & $R_{3}^{3}$$\begin{array}{cc}
\\
\\
\end{array}$ & $c_{3}$\tabularnewline
\cline{1-3} 
\multicolumn{1}{c}{$q_{1}$} & \multicolumn{1}{c}{$q_{2}$} & \multicolumn{1}{c}{$q_{3}$} & $\boxed{\boxed{\mathsf{SZLG}}}$\tabularnewline
\end{tabular}
\par\end{centering}

\caption{\label{fig: SZLG}The \mbox{c-c} matrix for a Suppes-Zanotti-Leggett-Garg-type
system. This is a cyclic system of rank 3, in the terminology of Section
\ref{sec: Cyclic-c-c-systems}. Each conteXt includes two conteNts,
and each conteNt is included in two conteXts. All random variables
are binary.}
\end{figure}

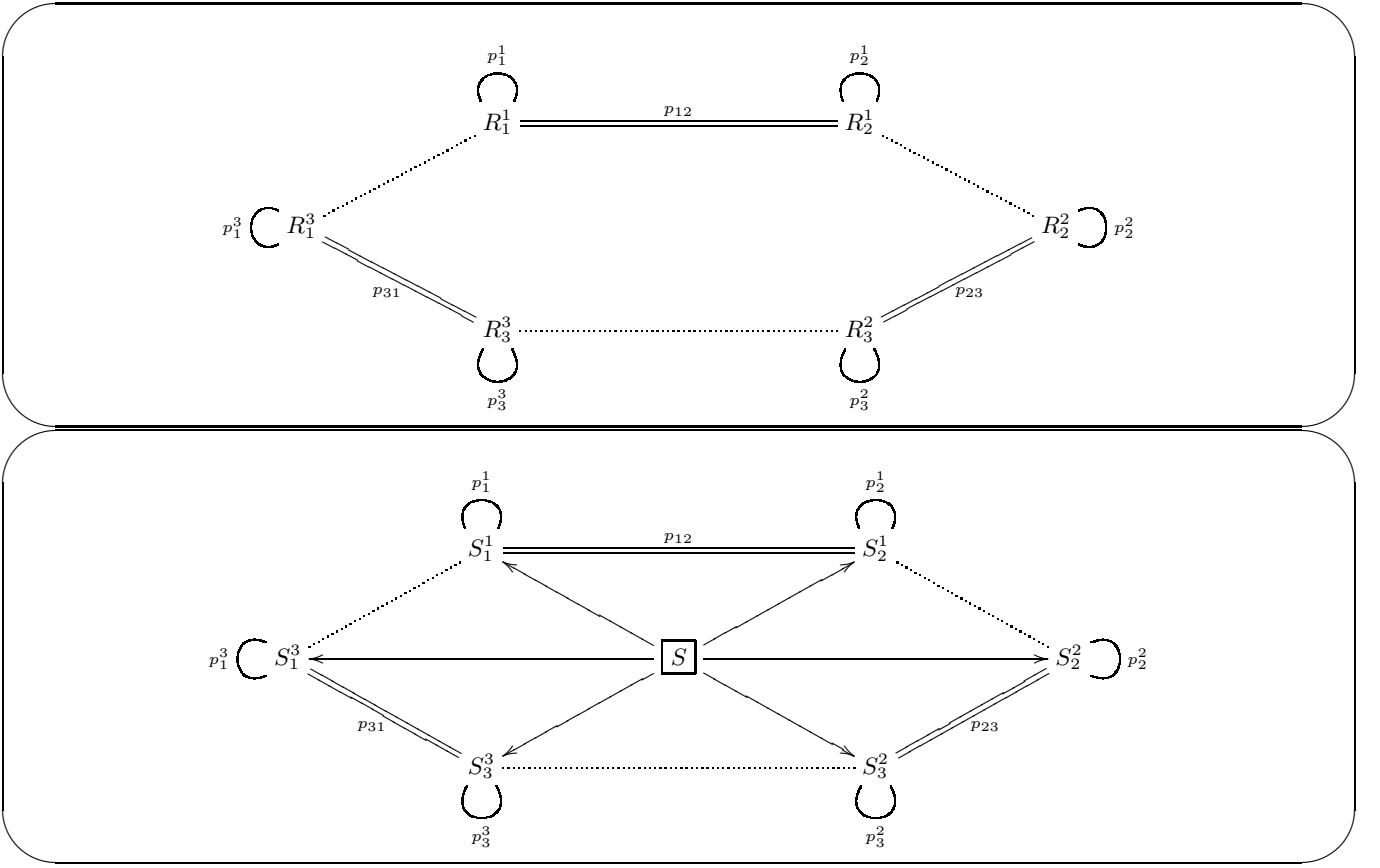
\begin{figure*}
\begin{centering}
\ovalbox{\begin{minipage}[t]{0.99\textwidth}%
\begin{center}
\[
\xymatrix@C=2cm{ & R_{1}^{1}\ar@{-}@(ul,ur)^{p_{1}^{1}}\ar@{=}[rr]^{p_{12}} &  & R_{2}^{1}\ar@{-}@(ur,ul)_{p_{2}^{1}}\ar@{.}[dr]\\
R_{1}^{3}\ar@{-}@(lu,ld)_{p_{1}^{3}}\ar@{.}[ur] &  &  &  & R_{2}^{2}\ar@{-}@(ru,rd)^{p_{2}^{2}}\ar@{=}[dl]^{p_{23}}\\
 & R_{3}^{3}\ar@{-}@(dr,dl)^{p_{3}^{3}}\ar@{=}[ul]^{p_{31}} &  & R_{3}^{2}\ar@{-}@(dl,dr)_{p_{3}^{2}}\ar@{.}[ll]
}
\]

\par\end{center}%
\end{minipage}}
\par\end{centering}

\begin{centering}
\ovalbox{\begin{minipage}[t]{0.99\textwidth}%
\begin{center}
\[
\xymatrix@C=2cm{ & S_{1}^{1}\ar@{-}@(ul,ur)^{p_{1}^{1}}\ar@{=}[rr]^{p_{12}} &  & S_{2}^{1}\ar@{-}@(ur,ul)_{p_{2}^{1}}\ar@{.}[dr]\\
S_{1}^{3}\ar@{-}@(lu,ld)_{p_{1}^{3}}\ar@{.}[ur] &  & \boxed{S}\ar[ru]\ar[ul]\ar[ll]\ar[rr]\ar[dl]\ar[dr] &  & S_{2}^{2}\ar@{-}@(ru,rd)^{p_{2}^{2}}\ar@{=}[dl]^{p_{23}}\\
 & S_{3}^{3}\ar@{-}@(dr,dl)^{p_{3}^{3}}\ar@{=}[ul]^{p_{31}} &  & S_{3}^{2}\ar@{-}@(dl,dr)_{p_{3}^{2}}\ar@{.}[ll]
}
\]

\par\end{center}%
\end{minipage}}
\par\end{centering}

\caption{\label{fig: SZLG and coupling}A schematic representation of a SZLG
\mbox{c-c} system (top) and of its coupling (bottom). The coupling
is shown as a set of random variables that are functions of some random
variable $S$. The symbols $p_{i}^{j}$ attached to $R_{i}^{j}$ show
the probability with which this variable equals 1; the symbols $p_{kl}$
attached to double-lines show the joint probability with which the
flanking variables equal 1. All these probabilities are preserved
in the coupling because, by definition, the bunches of the system
(the pairs of random variables connected by double-lines) have the
same distribution as the corresponding subcouplings of the coupling.
The elements of the connections of the system (the pairs of random
variables connected by dotted lines) are stochastically unrelated,
but the corresponding subcouplings of the coupling are jointly distributed,
as 2-marginals of the coupling. }
\end{figure*}

\begin{defn}
\label{def: cntx}A coupling for (the bunches of) a \mbox{c-c} system
is \emph{maximally connected} if its subcouplings corresponding to
the connections of the system are maximal couplings of these connections.
If a system has a maximally connected coupling, it is \emph{noncontextual}.
Otherwise it is \emph{contextual}. \footnote{In reference to footnote \ref{fn: As-the-reviewing} below, in a newer
version of CbD (Dzhafarov \& Kujala, 2016a,b) (non)contextuality is
defined in terms of multimaximal couplings, rather than merely maximal
ones. This does not, however, affects the logic of CbD in any nontrivial
way.}
\end{defn}
For motivation of this definition, see Section \ref{sub: The-intuition-of}.
Let us illustrate this definition using the system of binary random
variables first considered, in abstract, by Suppes and Zanotti (1981)
and then, as a paradigm in quantum mechanics, by Leggett and Garg
(1985). The \mbox{c-c} matrix for this system is presented in Fig.~\ref{fig: SZLG}.
There are three conteNts here, any two of which are represented (measured,
responded to) in one of three possible conteXts. Figure \ref{fig: SZLG and coupling}
(top panel) shows this system schematically: a set of three bunches
stochastically unrelated to each other, and three connections ``bridging''
them. Since any of the six random variables in the system has two
possible values, any coupling
\begin{equation}
S=\left(S_{1}^{1},S_{2}^{1},S_{2}^{2},S_{3}^{2},S_{3}^{3},S_{1}^{3}\right)
\end{equation}
of this system has $2^{6}$ possible values, and its distribution
is defined by assigning $2^{6}$ probability masses to them. 

It is difficult to see how one could show graphically that six random
variables are jointly distributed. In the bottom panel of Fig.~\ref{fig: SZLG and coupling}
this problem is solved by invoking Theorem \ref{thm: JDC}, according
to which the random variables in a coupling are all functions of one
and the same, ``hidden'' random variable. We do not need to specify
this random variable and the functions producing the components of
a coupling explicitly. It is always possible, however, to choose this
random variable to be the coupling $S$ itself, and treat the random
variables in the coupling as projection functions: $S_{1}^{1}$ is
the first projection of $S$, $S_{2}^{1}$ is the second projection
of $S$, etc.

The distribution of the six random variables in the coupling should,
by definition, agree with the bunches of the system $\mathcal{B}$,
each of which is uniquely characterized by three probabilities:
\begin{equation}
\begin{array}{c}
\Pr\left[R_{1}^{1}=1\right]=p_{1}^{1}=\Pr\left[S_{1}^{1}=1\right],\\
\\
\Pr\left[R_{2}^{1}=1\right]=p_{2}^{1}=\Pr\left[S_{2}^{1}=1\right],\\
\\
\Pr\left[R_{1}^{1}=1,R_{2}^{1}=1\right]=p_{12}=\Pr\left[S_{1}^{1}=1,S_{2}^{1}=1\right]
\end{array}
\end{equation}
for the bunch $\left(R_{1}^{1},R_{2}^{1}\right)$, and analogously
for the other two. There are generally an infinity of couplings that
satisfy these equations. 

\begin{figure*}
\begin{centering}
\ovalbox{\begin{minipage}[t]{0.99\textwidth}%
\begin{center}
\[
\xymatrix@C=2cm{ & R_{1}^{1}\ar@{-}@(ul,ur)^{p_{1}^{1}} &  & R_{2}^{1}\ar@{-}@(ur,ul)_{p_{2}^{1}}\ar@{.}[dr]\\
R_{1}^{3}\ar@{-}@(lu,ld)_{p_{1}^{3}}\ar@{.}[ur] &  &  &  & R_{2}^{2}\ar@{-}@(ru,rd)^{p_{2}^{2}}\\
 & R_{3}^{3}\ar@{-}@(dr,dl)^{p_{3}^{3}} &  & R_{3}^{2}\ar@{-}@(dl,dr)_{p_{3}^{2}}\ar@{.}[ll]
}
\]

\par\end{center}%
\end{minipage}}
\par\end{centering}

\begin{centering}
\ovalbox{\begin{minipage}[t]{0.99\textwidth}%
\begin{center}
\[
\xymatrix@C=2cm{ & T_{1}^{1}\ar@{-}@(ul,ur)^{p_{1}^{1}}\ar@{}[rr] &  & T_{2}^{1}\ar@{-}@(ur,ul)_{p_{2}^{1}}\ar@{.}[dr]^{p_{2}=\min\left(p_{2}^{1},p_{2}^{2}\right)}\\
T_{1}^{3}\ar@{-}@(lu,ld)_{p_{1}^{3}}\ar@{.}[ur]^{p_{1}=\min\left(p_{1}^{1},p_{1}^{3}\right)} & \boxed{T_{1}}\ar[u]\ar[l] & \boxed{T_{3}}\ar[dr]\ar[dl] & \boxed{T_{2}}\ar[r]\ar[u] & T_{2}^{2}\ar@{-}@(ru,rd)^{p_{2}^{2}}\ar@{}[dl]\\
 & T_{3}^{3}\ar@{-}@(dr,dl)^{p_{3}^{3}}\ar@{}[ul] &  & T_{3}^{2}\ar@{-}@(dl,dr)_{p_{3}^{2}}\ar@{.}[ll]^{p_{3}=\min\left(p_{3}^{2},p_{3}^{3}\right)}
}
\]

\par\end{center}%
\end{minipage}}
\par\end{centering}

\begin{centering}
\ovalbox{\begin{minipage}[t]{0.99\textwidth}%
\begin{center}
\[
\xymatrix@C=2cm{ & S_{1}^{1}\ar@{-}@(ul,ur)^{p_{1}^{1}}\ar@{=}[rr]^{p_{12}} &  & S_{2}^{1}\ar@{-}@(ur,ul)_{p_{2}^{1}}\ar@{.}[dr]^{p_{2}=\min\left(p_{2}^{1},p_{2}^{2}\right)}\\
S_{1}^{3}\ar@{-}@(lu,ld)_{p_{1}^{3}}\ar@{.}[ur]^{p_{1}=\min\left(p_{1}^{1},p_{1}^{3}\right)} &  & \boxed{S}\ar[ru]\ar[ul]\ar[ll]\ar[rr]\ar[dl]\ar[dr] &  & S_{2}^{2}\ar@{-}@(ru,rd)^{p_{2}^{2}}\ar@{=}[dl]^{p_{23}}\\
 & S_{3}^{3}\ar@{-}@(dr,dl)^{p_{3}^{3}}\ar@{=}[ul]^{p_{31}} &  & S_{3}^{2}\ar@{-}@(dl,dr)_{p_{3}^{2}}\ar@{.}[ll]^{p_{3}=\min\left(p_{3}^{2},p_{3}^{3}\right)}
}
\]

\par\end{center}%
\end{minipage}}
\par\end{centering}

\caption{\label{fig: SZLG connections}Each of the three connections of the
SZLG system of Fig.~\ref{fig: SZLG and coupling}, taken separately
(top), can be coupled by a maximal coupling (middle). A hypothetical
maximally connected coupling of the SZLG system (bottom) is one in
which the subcouplings corresponding to the connections of the system
are their maximal couplings. If such a coupling can be constructed
(equivalently, if the maximal couplings in the middle panel are compatible
with the system's bunches), then the system is noncontextual. It is
possible that such a coupling does not exist, in which case the system
is contextual.}
\end{figure*}
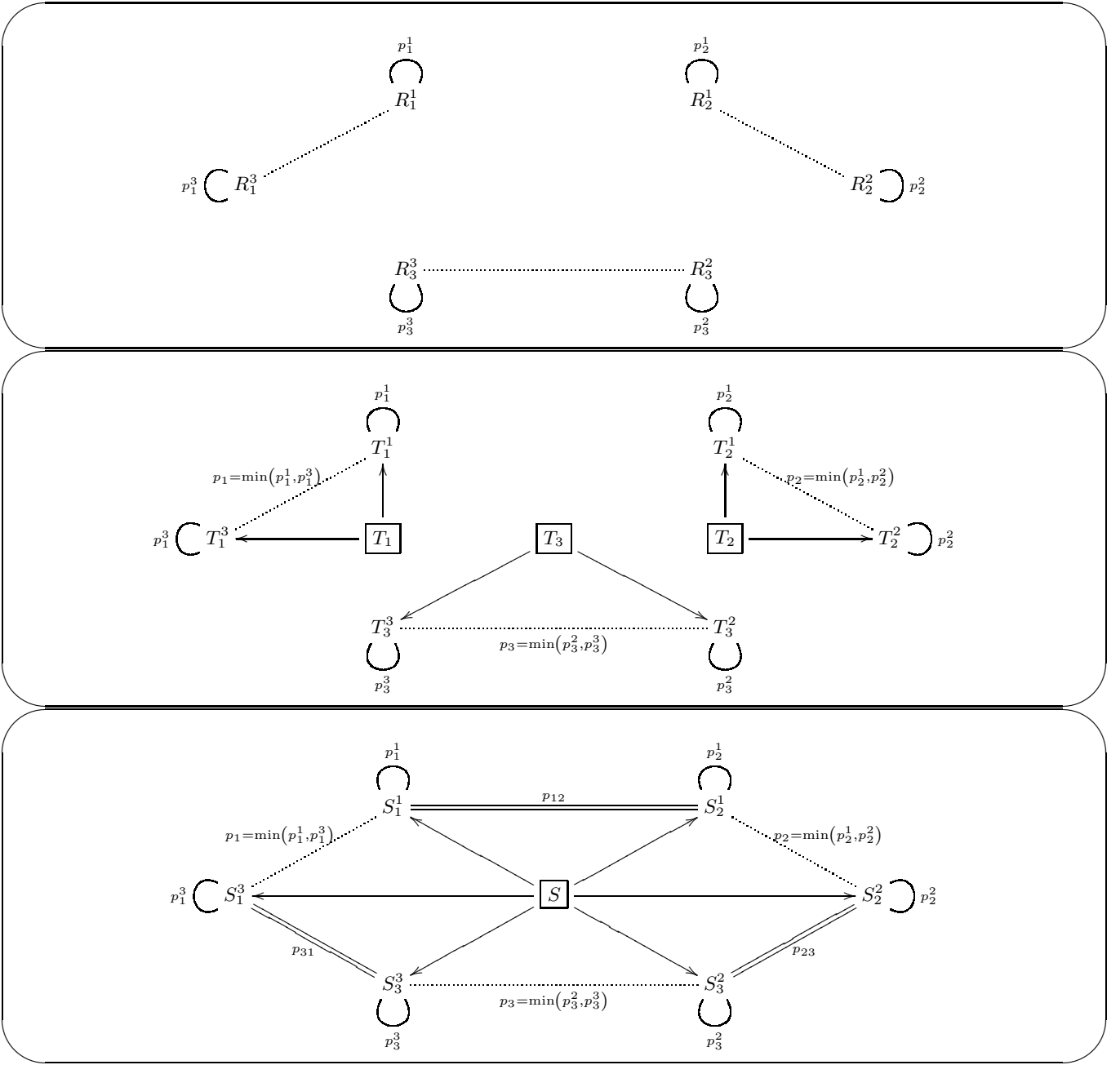

Consider now the three connections of the system as three separate
pairs of random variables (Fig.~\ref{fig: SZLG connections}, top
panel), and for each of them consider its coupling (Fig.~\ref{fig: SZLG connections},
middle panel). The distributions of the elements of a connection are
fixed, and its coupling should, by definition, preserve them. Thus,
for the connection $\left(R_{1}^{1},R_{1}^{3}\right)$,
\begin{equation}
\begin{array}{c}
\Pr\left[R_{1}^{1}=1\right]=p_{1}^{1}=\Pr\left[T_{1}^{1}=1\right],\\
\\
\Pr\left[R_{1}^{3}=1\right]=p_{1}^{3}=\Pr\left[T_{1}^{3}=1\right].
\end{array}
\end{equation}
With $p_{1}^{1}$ and $p_{1}^{3}$ given, the distribution of the
coupling $\left(T_{1}^{1},T_{1}^{3}\right)$ of $\left(R_{1}^{1},R_{1}^{3}\right)$
is determined by
\begin{equation}
\Pr\left[T_{1}^{1}=1,T_{1}^{3}=1\right]=p_{1}.
\end{equation}
By Theorem \ref{thm: maximal-coupling}, $p_{1}$ can be chosen so
that $\Pr\left[T_{1}^{1}=T_{1}^{3}\right]$ attains its maximal possible
value, and this choice is
\begin{equation}
p_{1}=\min\left(p_{1}^{1},p_{1}^{3}\right).
\end{equation}
We know that such a coupling is a maximal coupling of $\left(R_{1}^{1},R_{1}^{3}\right)$,
in this simple case, uniquely determined. We choose values $p_{2}$
and $p_{3}$ for the maximal couplings of the remaining two connections
analogously.

In accordance with Definition \ref{def: cntx}, the question now is
whether these three values of the joint probabilities, $p_{1},p_{2},p_{3}$,
are compatible with the bunches of the system. Put differently, can
one construct a maximally connected coupling shown in Fig.~\ref{fig: SZLG connections},
bottom panel, a coupling in which all the probabilities shown are
achieved together? The system is noncontextual if and only if the
answer to this question is affirmative. 

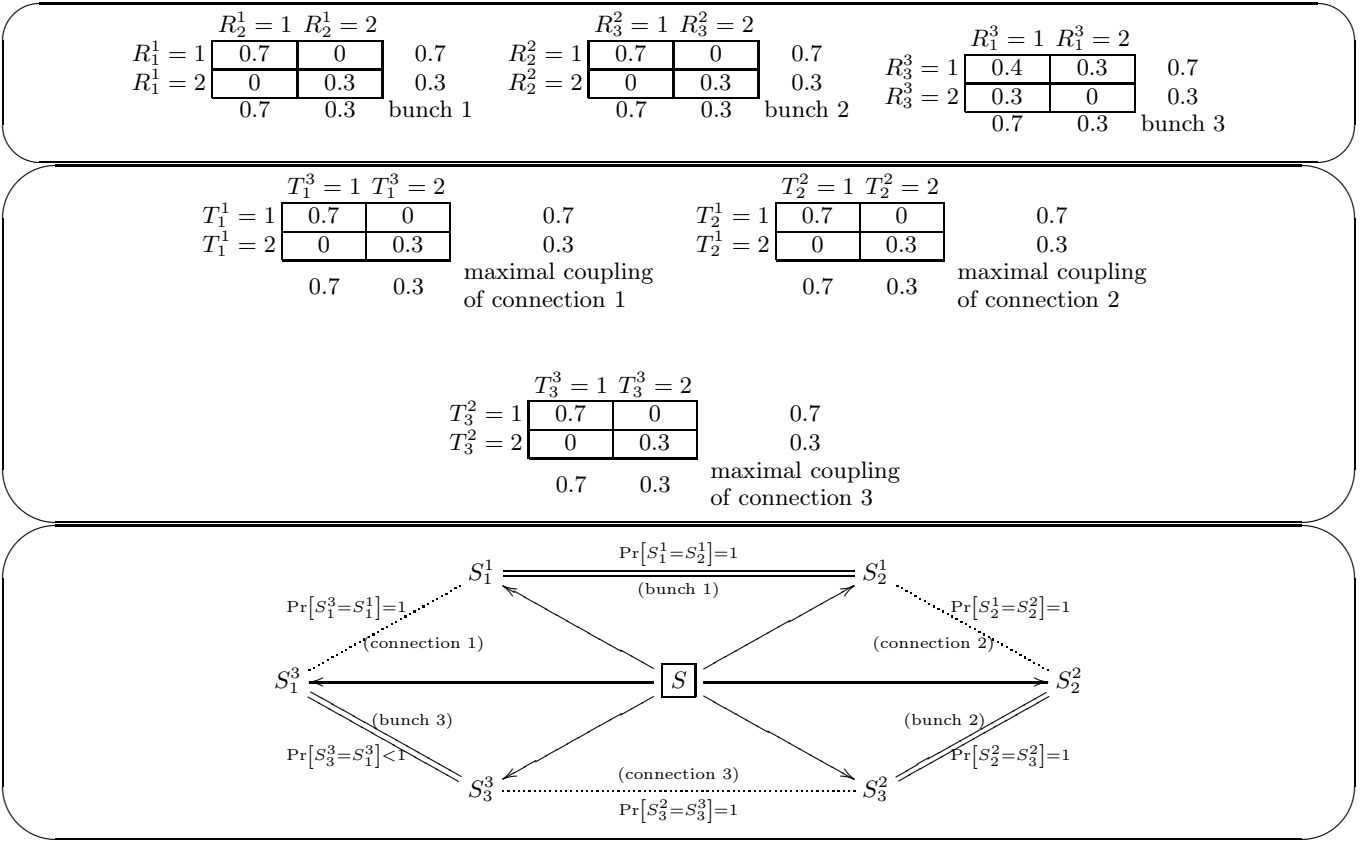
\begin{figure*}
\begin{centering}
\ovalbox{\begin{minipage}[t]{0.99\textwidth}%
\begin{center}
\begin{tabular}{cccc}
 & $R_{2}^{1}=1$ & $R_{2}^{1}=2$ & \tabularnewline
\cline{2-3} 
\multicolumn{1}{c|}{$R_{1}^{1}=1$} & \multicolumn{1}{c|}{0.7} & \multicolumn{1}{c|}{0} & 0.7\tabularnewline
\cline{2-3} 
\multicolumn{1}{c|}{$R_{1}^{1}=2$} & \multicolumn{1}{c|}{0} & \multicolumn{1}{c|}{0.3} & 0.3\tabularnewline
\cline{2-3} 
 & 0.7 & 0.3 & bunch 1\tabularnewline
 &  &  & \tabularnewline
\end{tabular}$\quad$%
\begin{tabular}{cccc}
 & $R_{3}^{2}=1$ & $R_{3}^{2}=2$ & \tabularnewline
\cline{2-3} 
\multicolumn{1}{c|}{$R_{2}^{2}=1$} & \multicolumn{1}{c|}{0.7} & \multicolumn{1}{c|}{0} & 0.7\tabularnewline
\cline{2-3} 
\multicolumn{1}{c|}{$R_{2}^{2}=2$} & \multicolumn{1}{c|}{0} & \multicolumn{1}{c|}{0.3} & 0.3\tabularnewline
\cline{2-3} 
 & 0.7 & 0.3 & bunch 2\tabularnewline
 &  &  & \tabularnewline
\end{tabular}$\quad$%
\begin{tabular}{c|c|c|c}
\multicolumn{1}{c}{} & \multicolumn{1}{c}{$R_{1}^{3}=1$} & \multicolumn{1}{c}{$R_{1}^{3}=2$} & \tabularnewline
\cline{2-3} 
$R_{3}^{3}=1$ & 0.4 & 0.3 & 0.7\tabularnewline
\cline{2-3} 
$R_{3}^{3}=2$ & 0.3 & 0 & 0.3\tabularnewline
\cline{2-3} 
\multicolumn{1}{c}{} & \multicolumn{1}{c}{0.7} & \multicolumn{1}{c}{0.3} & bunch 3\tabularnewline
\end{tabular}
\par\end{center}%
\end{minipage}}
\par\end{centering}

\begin{centering}
\ovalbox{\begin{minipage}[t]{0.99\textwidth}%
\begin{center}
\begin{tabular}{cccc}
 & $T_{1}^{3}=1$ & $T_{1}^{3}=2$ & \tabularnewline
\cline{2-3} 
\multicolumn{1}{c|}{$T_{1}^{1}=1$} & \multicolumn{1}{c|}{0.7} & \multicolumn{1}{c|}{0} & 0.7\tabularnewline
\cline{2-3} 
\multicolumn{1}{c|}{$T_{1}^{1}=2$} & \multicolumn{1}{c|}{0} & \multicolumn{1}{c|}{0.3} & 0.3\tabularnewline
\cline{2-3} 
 & 0.7 & 0.3 & $\begin{array}{l}
\textnormal{maximal coupling}\\
\textnormal{of connection 1}
\end{array}$\tabularnewline
 &  &  & \tabularnewline
\end{tabular}$\quad$%
\begin{tabular}{cccc}
 & $T_{2}^{2}=1$ & $T_{2}^{2}=2$ & \tabularnewline
\cline{2-3} 
\multicolumn{1}{c|}{$T_{2}^{1}=1$} & \multicolumn{1}{c|}{0.7} & \multicolumn{1}{c|}{0} & 0.7\tabularnewline
\cline{2-3} 
\multicolumn{1}{c|}{$T_{2}^{1}=2$} & \multicolumn{1}{c|}{0} & \multicolumn{1}{c|}{0.3} & 0.3\tabularnewline
\cline{2-3} 
 & 0.7 & 0.3 & $\begin{array}{l}
\textnormal{maximal coupling}\\
\textnormal{of connection 2}
\end{array}$\tabularnewline
 &  &  & \tabularnewline
\end{tabular}
\par\end{center}

\begin{center}
\begin{tabular}{c|c|c|c}
\multicolumn{1}{c}{} & \multicolumn{1}{c}{$T_{3}^{3}=1$} & \multicolumn{1}{c}{$T_{3}^{3}=2$} & \tabularnewline
\cline{2-3} 
$T_{3}^{2}=1$ & 0.7 & 0 & 0.7\tabularnewline
\cline{2-3} 
$T_{3}^{2}=2$ & 0 & 0.3 & 0.3\tabularnewline
\cline{2-3} 
\multicolumn{1}{c}{} & \multicolumn{1}{c}{0.7} & \multicolumn{1}{c}{0.3} & $\begin{array}{l}
\textnormal{maximal coupling}\\
\textnormal{of connection 3}
\end{array}$\tabularnewline
\end{tabular}
\par\end{center}%
\end{minipage}}
\par\end{centering}

\begin{centering}
\ovalbox{\begin{minipage}[t]{0.99\textwidth}%
\begin{center}
$\xymatrix@C=2cm{ & S_{1}^{1}\ar@{=}[rr]_{\textnormal{(bunch 1)}}^{\Pr\left[S_{1}^{1}=S_{2}^{1}\right]=1} &  & S_{2}^{1}\ar@{.}[dr]_{\textnormal{(connection 2)}}^{\Pr\left[S_{2}^{1}=S_{2}^{2}\right]=1}\\
S_{1}^{3}\ar@{.}[ur]_{\textnormal{(connection 1)}}^{\Pr\left[S_{1}^{3}=S_{1}^{1}\right]=1} &  & \boxed{S}\ar[ru]\ar[ul]\ar[ll]\ar[rr]\ar[dl]\ar[dr] &  & S_{2}^{2}\ar@{=}[dl]_{\textnormal{(bunch 2)}}^{\Pr\left[S_{2}^{2}=S_{3}^{2}\right]=1}\\
 & S_{3}^{3}\ar@{=}[ul]_{\textnormal{(bunch 3)}}^{\Pr\left[S_{3}^{3}=S_{1}^{3}\right]<1} &  & S_{3}^{2}\ar@{.}[ll]_{\textnormal{(connection 3)}}^{\Pr\left[S_{3}^{2}=S_{3}^{3}\right]=1}
}
$
\par\end{center}%
\end{minipage}}
\par\end{centering}

\caption{\label{fig: SZLG exmaple}An example of a contextual SZLG system.
The bunches are shown in the top panel. The middle panel shows the
computed maximal couplings for the connections. The bottom panel shows
that a maximally connected coupling does not exist because it is internally
contradictory. To each edge of the diagram in this panel we append
the probability with which the flanking random variables are equal
to each other, and we show from where this probability is derived.
Thus, $\Pr\left[S_{1}^{1}=S_{2}^{1}\right]=1$ is derived from the
distribution of the first bunch, while $\Pr\left[S_{1}^{3}=S_{1}^{1}\right]=1$
is derived from the distribution of the maximal coupling of the first
connection.}
\end{figure*}

To see that it does not have to be affirmative, consider the example
presented in Fig.~\ref{fig: SZLG exmaple}. The maximally connected
coupling does not exist because, in the bottom-panel diagram, 
\begin{lyxlist}{00.00.0000}
\item [{(i)}] going clockwise from $S_{1}^{1}$ and using the transitivity
of the relation ``always equals,'' we conclude that $\Pr\left[S_{1}^{1}=S_{3}^{3}\right]=1$;
\item [{(ii)}] going counterclockwise from $S_{1}^{1}$ we see that $\Pr\left[S_{1}^{1}=S_{1}^{3}\right]=1$;
\item [{(iii)}] but then $\Pr\left[S_{1}^{3}=S_{3}^{3}\right]$ must also
be 1, which it is not. 
\end{lyxlist}

\section{\label{sec: Contextuality-as-LP}Contextuality as a linear programming
problem}

Is there a general method for establishing contextuality or lack thereof
in a given \mbox{c-c} system? It turns out that such a method exists,
and insofar as finite sets of categorical random variables are involved,
it is a simple linear programming method. A maximally connected coupling
of a \mbox{c-c} system is uniquely associated with a certain underdetermined
system of linear equations, and the \mbox{c-c} system is contextual
if and only if this system of linear equations has no nonnegative
solutions. The theory of these equations generalizes the Linear Feasibility
Test described in Dzhafarov and Kujala (2012).

\subsection{Notation and conventions}

We need to introduce or recall some notation and conventions. Let
a system $\mathcal{R}$ involve conteXts $c_{1},\ldots,c_{n}$ ($n>1$)
and conteNts $q_{1},\ldots,q_{m}$ ($m>1$). 
\begin{enumerate}
\item Notation related to a \mbox{c-c} system $\mathcal{R}$:
\begin{equation}
\mathcal{R}=\left(\begin{array}{cccc}
\ddots & \vdots & \iddots\\
\cdots & R_{j}^{i} & \cdots & \left(\textnormal{bunch }R^{i}\right)\\
\iddots & \vdots & \ddots\\
 & \left(\begin{array}{c}
\textnormal{connection}\\
\mathcal{R}_{j}
\end{array}\right)
\end{array}\right).
\end{equation}

\item Corresponding notation for a (maximally connected) coupling $S$ of
$\mathcal{R}$: 
\begin{equation}
S=\left(\begin{array}{cccc}
\ddots & \vdots & \iddots\\
\cdots & S_{j}^{i} & \cdots & \left(\begin{array}{c}
\textnormal{subcoupling }\\
S^{i}
\end{array}\right)\\
\iddots & \vdots & \ddots\\
 & \left(\begin{array}{c}
\textnormal{subcoupling}\\
S_{j}
\end{array}\right)
\end{array}\right).
\end{equation}

\item Notation for any (maximal) coupling $T_{j}$ of a connection $\mathcal{R}_{j}$
taken separately:
\begin{equation}
T_{j}=\left(\begin{array}{c}
\vdots\\
T_{j}^{i}\\
\vdots
\end{array}\right).
\end{equation}

\item A value of a random variable $R_{j}^{i}$ (hence also of $S_{j}^{i}$
or $T_{j}^{i}$) is denoted $v_{j}^{i}$ or $w_{j}^{i}$. A value
of a bunch $R^{i}$ (hence also of the subcoupling $S^{i}$ of $S$)
is denoted $v^{i}$ or $w^{i}$. We use $v_{j}$ or $w_{j}$ to denote
values of couplings $T_{j}$ and the corresponding subcouplings $S_{j}$
of $S$ (assumed to be maximally connected). The value $v$ of $S$
has the structure 
\begin{equation}
v=\left(\begin{array}{cccc}
\ddots & \vdots & \iddots\\
\cdots & v_{j}^{i} & \cdots & \left(\begin{array}{c}
\textnormal{bunch}\\
\textnormal{ value }v^{i}
\end{array}\right)\\
\iddots & \vdots & \ddots\\
 & \left(\begin{array}{c}
\textnormal{connection}\\
\textnormal{ value }v_{j}
\end{array}\right)
\end{array}\right).\label{eq: value v}
\end{equation}
As is customary, we use $v,v^{i},v_{j},v_{j}^{i}$ sometimes as variables
and sometimes as specific values of these variables.
\item Recall that in these matrices and vectors some entries are not defined:
not every conteNt is paired with every conteXt. If $q_{j}$ is measured
(responded to) in conteXt $c_{i}$, the random variable $R_{j}^{i}$
exists, and the elements of the set $V_{j}$ of its possible values
can be enumerated $1,\ldots,k_{j}$. Denoting 
\begin{equation}
k=\max_{j=1,\ldots,m}k_{j},
\end{equation}
without loss of generality, we can assume that 
\begin{equation}
V_{j}=\left\{ 1,\ldots,k\right\} ,
\end{equation}
for every $j=1,\ldots,m$. Indeed, one can always add values to $V_{j}$
that occur with probability zero. The set of all values $v$ of a
coupling $S$ is therefore $\left\{ 1,\ldots,k\right\} ^{N}$, where
$N$ is the number of all random variables in $S$.
\item We will refer to the values $v$ of $S$ as \emph{hidden outcomes}.
The term derives from quantum mechanics, where the problem of contextuality
was initially presented as that of hidden variables (see the last
paragraph of Section \ref{sub: Functions-of-random}).
\end{enumerate}

\begin{figure*}
\begin{footnotesize}%
\begin{tabular}{|c|c|c|c|c|c|c|c|c|c|c|c|c|c|c|c|c|c|}
\cline{1-16} 
\multicolumn{16}{|c|}{hidden outcome, $\left(\begin{array}{cc}
v_{1}^{1} & v_{2}^{1}\\
v_{1}^{2} & v_{2}^{2}
\end{array}\right)$} & \multicolumn{1}{c}{\multirow{3}{*}{}} & \multicolumn{1}{c}{}\tabularnewline
$\!\begin{array}{cc}
+ & +\\
+ & +
\end{array}\!$ & $\!\begin{array}{cc}
+ & +\\
+ & -
\end{array}\!$ & $\!\begin{array}{cc}
+ & +\\
- & +
\end{array}\!$ & $\!\begin{array}{cc}
+ & +\\
- & -
\end{array}\!$ & $\!\begin{array}{cc}
+ & -\\
+ & +
\end{array}\!$ & $\!\begin{array}{cc}
+ & -\\
+ & -
\end{array}\!$ & $\!\begin{array}{cc}
+ & -\\
- & +
\end{array}\!$ & $\!\begin{array}{cc}
+ & -\\
- & -
\end{array}\!$ & $\!\begin{array}{cc}
- & +\\
+ & +
\end{array}\!$ & $\!\begin{array}{cc}
- & +\\
+ & -
\end{array}\!$ & $\!\begin{array}{cc}
- & +\\
- & +
\end{array}\!$ & $\!\begin{array}{cc}
- & +\\
- & -
\end{array}\!$ & $\!\begin{array}{cc}
- & -\\
+ & +
\end{array}\!$ & $\!\begin{array}{cc}
- & -\\
+ & -
\end{array}\!$ & $\!\begin{array}{cc}
- & -\\
- & +
\end{array}\!$ & $\!\begin{array}{cc}
- & -\\
- & -
\end{array}\!$ &  & \multicolumn{1}{c}{}\tabularnewline
\cline{1-16} \cline{18-18} 
1 & 1 & 1 & 1 &  &  &  &  &  &  &  &  &  &  &  &  &  & $\begin{array}{l}
\Pr\left[S^{1}=\left(+,+\right)\right]\\
=\Pr\left[R^{1}=\left(+,+\right)\right]
\end{array}$\tabularnewline
\cline{1-16} \cline{18-18} 
 &  &  &  & 1 & 1 & 1 & 1 &  &  &  &  &  &  &  &  &  & $\begin{array}{l}
\Pr\left[S^{1}=\left(+,-\right)\right]\\
=\Pr\left[R^{1}=\left(+,-\right)\right]
\end{array}$\tabularnewline
\cline{1-16} \cline{18-18} 
 &  &  &  &  &  &  &  & 1 & 1 & 1 & 1 &  &  &  &  &  & $\begin{array}{l}
\Pr\left[S^{1}=\left(-,+\right)\right]\\
=\Pr\left[R^{1}=\left(-,+\right)\right]
\end{array}$\tabularnewline
\cline{1-16} \cline{18-18} 
 &  &  &  &  &  &  &  &  &  &  &  & 1 & 1 & 1 & 1 &  & $\begin{array}{l}
\Pr\left[S^{1}=\left(-,-\right)\right]\\
=\Pr\left[R^{1}=\left(-,-\right)\right]
\end{array}$\tabularnewline
\cline{1-16} \cline{18-18} 
1 &  &  &  & 1 &  &  &  & 1 &  &  &  & 1 &  &  &  &  & $\begin{array}{l}
\Pr\left[S^{2}=\left(+,+\right)\right]\\
=\Pr\left[R^{2}=\left(+,+\right)\right]
\end{array}$\tabularnewline
\cline{1-16} \cline{18-18} 
 & 1 &  &  &  & 1 &  &  &  & 1 &  &  &  & 1 &  &  &  & $\begin{array}{l}
\Pr\left[S^{2}=\left(+,-\right)\right]\\
=\Pr\left[R^{2}=\left(+,-\right)\right]
\end{array}$\tabularnewline
\cline{1-16} \cline{18-18} 
 &  & 1 &  &  &  & 1 &  &  &  & 1 &  &  &  & 1 &  &  & $\begin{array}{l}
\Pr\left[S^{2}=\left(-,+\right)\right]\\
=\Pr\left[R^{2}=\left(-,+\right)\right]
\end{array}$\tabularnewline
\cline{1-16} \cline{18-18} 
 &  &  & 1 &  &  &  & 1 &  &  &  & 1 &  &  &  & 1 &  & $\begin{array}{l}
\Pr\left[S^{2}=\left(-,-\right)\right]\\
=\Pr\left[R^{2}=\left(-,-\right)\right]
\end{array}$\tabularnewline
\cline{1-16} \cline{18-18} 
1 & 1 &  &  & 1 & 1 &  &  &  &  &  &  &  &  &  &  &  & $\begin{array}{l}
\Pr\left[S_{1}=\left(+,+\right)\right]\\
=\Pr\left[T_{1}=\left(+,+\right)\right]
\end{array}$\tabularnewline
\cline{1-16} \cline{18-18} 
 &  &  &  &  &  &  &  &  &  & 1 & 1 &  &  & 1 & 1 &  & $\begin{array}{l}
\Pr\left[S_{1}=\left(-,-\right)\right]\\
=\Pr\left[T_{1}=\left(-,-\right)\right]
\end{array}$\tabularnewline
\cline{1-16} \cline{18-18} 
1 &  & 1 &  &  &  &  &  & 1 &  & 1 &  &  &  &  &  &  & $\begin{array}{l}
\Pr\left[S_{2}=\left(+,+\right)\right]\\
=\Pr\left[T_{2}=\left(+,+\right)\right]
\end{array}$\tabularnewline
\cline{1-16} \cline{18-18} 
 &  &  &  &  & 1 &  & 1 &  &  &  &  &  & 1 &  & 1 &  & $\begin{array}{l}
\Pr\left[S_{2}=\left(-,-\right)\right]\\
=\Pr\left[T_{2}=\left(-,-\right)\right]
\end{array}$\tabularnewline
\cline{1-16} \cline{18-18} 
\end{tabular}\end{footnotesize}

\caption{\label{tab:The-Boolean-matrix}The Boolean matrix $\mathbf{M}$ (left)
and vector $\mathbf{P}$ (right) for the \mbox{c-c} system $\mathcal{A}$
in Fig.~\ref{fig:system 2by2}. The values of the variables are encoded
by $\pm1$, with $+1$ shown by plus sign and $-1$ shown by minus
sign. The first 8 elements of $\mathbf{P}$ are bunch probabilities,
the last 4 elements are connection probabilities.}
\end{figure*}

\begin{figure*}
\begin{centering}
\ovalbox{\begin{minipage}[t]{0.7\textwidth}%
\begin{center}
\begin{tabular}{cccc}
 & $R_{2}^{1}=+1$ & $R_{2}^{1}=-1$ & \tabularnewline
\cline{2-3} 
\multicolumn{1}{c|}{$R_{1}^{1}=+1$} & \multicolumn{1}{c|}{$\frac{1}{2}$$\begin{array}{c}
\\
\\
\end{array}$} & \multicolumn{1}{c|}{$0$$\begin{array}{c}
\\
\\
\end{array}$} & $\frac{1}{2}$\tabularnewline
\cline{2-3} 
\multicolumn{1}{c|}{$R_{1}^{1}=-1$} & \multicolumn{1}{c|}{$0$$\begin{array}{c}
\\
\\
\end{array}$} & \multicolumn{1}{c|}{$\frac{1}{2}$$\begin{array}{c}
\\
\\
\end{array}$} & $\frac{1}{2}$\tabularnewline
\cline{2-3} 
 & $\frac{1}{2}$ & $\frac{1}{2}$ & $\boxed{b1}$\tabularnewline
 &  &  & \tabularnewline
\end{tabular}$\quad$%
\begin{tabular}{cccc}
 & $R_{2}^{2}=+1$ & $R_{2}^{2}=-1$ & \tabularnewline
\cline{2-3} 
\multicolumn{1}{c|}{$R_{1}^{2}=+1$} & \multicolumn{1}{c|}{0$\begin{array}{c}
\\
\\
\end{array}$} & \multicolumn{1}{c|}{$\frac{1}{2}$$\begin{array}{c}
\\
\\
\end{array}$} & $\frac{1}{2}$\tabularnewline
\cline{2-3} 
\multicolumn{1}{c|}{$R_{1}^{2}=-1$} & \multicolumn{1}{c|}{$\frac{1}{2}$$\begin{array}{c}
\\
\\
\end{array}$} & \multicolumn{1}{c|}{0$\begin{array}{c}
\\
\\
\end{array}$} & $\frac{1}{2}$\tabularnewline
\cline{2-3} 
 & $\frac{1}{2}$ & $\frac{1}{2}$ & $\boxed{b2}$\tabularnewline
 &  &  & \tabularnewline
\end{tabular}
\par\end{center}%
\end{minipage}}
\par\end{centering}

\begin{centering}
\ovalbox{\begin{minipage}[t]{0.7\textwidth}%
\begin{center}
\begin{tabular}{cccc}
 & $T_{1}^{2}=+1$ & $T_{1}^{2}=-1$ & \tabularnewline
\cline{2-3} 
\multicolumn{1}{c|}{$T_{1}^{1}=+1$} & \multicolumn{1}{c|}{$\frac{1}{2}$$\begin{array}{c}
\\
\\
\end{array}$} & \multicolumn{1}{c|}{$0$$\begin{array}{c}
\\
\\
\end{array}$} & $\frac{1}{2}$\tabularnewline
\cline{2-3} 
\multicolumn{1}{c|}{$T_{1}^{1}=-1$} & \multicolumn{1}{c|}{$0$$\begin{array}{c}
\\
\\
\end{array}$} & \multicolumn{1}{c|}{$\frac{1}{2}$$\begin{array}{c}
\\
\\
\end{array}$} & $\frac{1}{2}$\tabularnewline
\cline{2-3} 
 & $\frac{1}{2}$ & $\frac{1}{2}$ & $\boxed{c1}$\tabularnewline
 &  &  & \tabularnewline
\end{tabular}$\quad$%
\begin{tabular}{cccc}
 & $T_{2}^{2}=+1$ & $T_{2}^{2}=-1$ & \tabularnewline
\cline{2-3} 
\multicolumn{1}{c|}{$T_{2}^{1}=+1$} & \multicolumn{1}{c|}{$\frac{1}{2}$$\begin{array}{c}
\\
\\
\end{array}$} & \multicolumn{1}{c|}{$0$$\begin{array}{c}
\\
\\
\end{array}$} & $\frac{1}{2}$\tabularnewline
\cline{2-3} 
\multicolumn{1}{c|}{$T_{2}^{1}=-1$} & \multicolumn{1}{c|}{$0$$\begin{array}{c}
\\
\\
\end{array}$} & \multicolumn{1}{c|}{$\frac{1}{2}$$\begin{array}{c}
\\
\\
\end{array}$} & $\frac{1}{2}$\tabularnewline
\cline{2-3} 
 & $\frac{1}{2}$ & $\frac{1}{2}$ & $\boxed{c2}$\tabularnewline
 &  &  & \tabularnewline
\end{tabular}
\par\end{center}%
\end{minipage}}
\par\end{centering}

\begin{centering}
\ovalbox{\begin{minipage}[t]{0.7\textwidth}%
\begin{center}
$\xymatrix@C=2cm{S_{1}^{1}\ar@{=}[rr]_{\Pr\left[S_{1}^{1}=S_{2}^{1}\right]=1}^{\boxed{b1}} &  & S_{2}^{1}\ar@{.}[dd]_{\Pr\left[S_{2}^{1}=S_{2}^{2}\right]=1}^{\boxed{c2}}\\
 & \boxed{S}\ar[ru]\ar[ul]\ar[dl]\ar[dr]\\
S_{1}^{2}\ar@{.}[uu]_{\Pr\left[S_{1}^{1}=S_{1}^{2}\right]=1}^{\boxed{c1}} &  & S_{2}^{2}\ar@{=}[ll]_{\Pr\left[S_{2}^{2}=S_{1}^{2}\right]=0}^{\boxed{b2}}
}
$
\par\end{center}%
\end{minipage}}
\par\end{centering}

\caption{\label{fig: QQ system}An example of a contextual \mbox{c-c} system
of the $\mathcal{A}$-type (Fig. \ref{fig:system 2by2}). The bunch
probabilities are shown in the top panel. The middle panel shows the
computed maximal connection probabilities. The bottom panel shows
(using the same format and logic as in Fig.~\ref{fig: SZLG exmaple})
that a maximally connected coupling $S$ would be internally contradictory. }
\end{figure*}
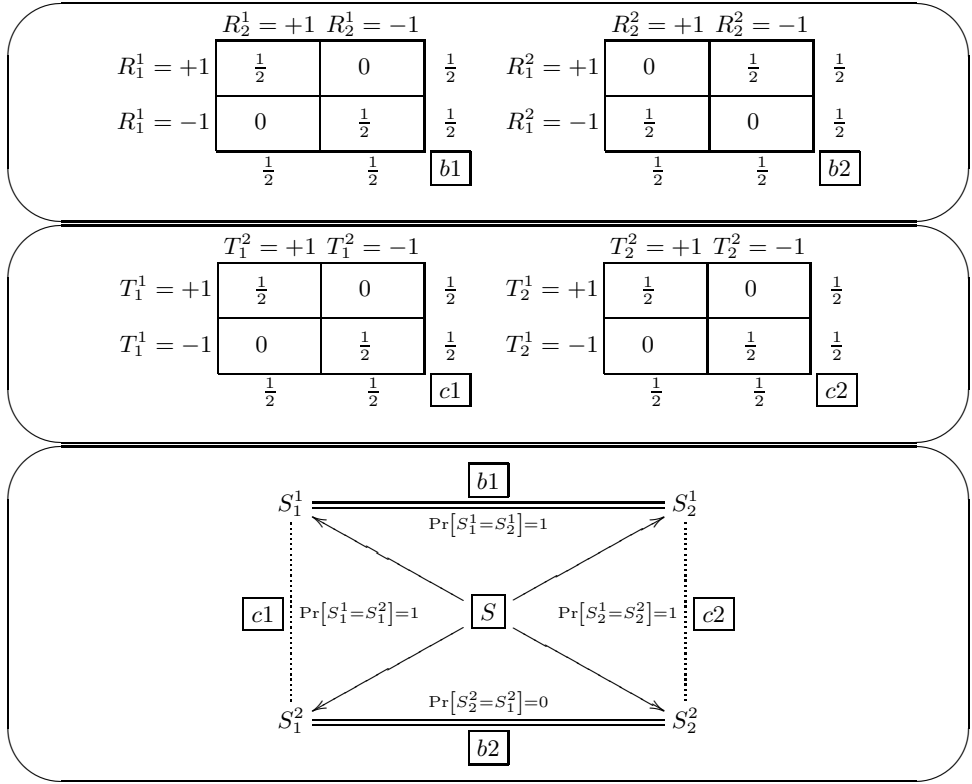

\subsection{Linear equations associated with a c-c system}

To specify a distribution of $S$, each of the hidden outcomes $v$
should be assigned a probability mass $\gamma\left(v\right)$. Let
us form a column vector $\mathbf{Q}$ by arranging these $\gamma\left(v\right)$-values
in some, say, lexicographic order of $v$. Let us also form a column
vector $\mathbf{P}$ with the following structure:
\begin{equation}
\mathbf{P}=\left(\begin{array}{c}
\textnormal{\emph{bunch}}\\
\textnormal{\emph{probabilities}}\\
\textnormal{for }c_{1},\ldots,c_{n}
\end{array},\ldots,\begin{array}{c}
\textnormal{\emph{connection}}\\
\textnormal{\emph{probabilities}}\\
\textnormal{for }q_{1},\ldots,q_{m}
\end{array}\right).
\end{equation}
Here,
\begin{equation}
\begin{array}{c}
\textnormal{bunch}\\
\textnormal{probabilities}\\
\textnormal{for }c_{i}
\end{array}=\left(\Pr\left[R^{i}=v^{i}\right]:v^{i}\in\left\{ 1,\ldots,k\right\} ^{n_{i}}\right),
\end{equation}
where $n_{i}$ is the number of elements in $v^{i}$. That is, the
bunch probabilities for $c{}_{i}$ are the joint probabilities that
determine the distribution of the bunch $R^{i}$. The connection probabilities
for $q_{j}$ are the probabilities imposed by the maximal coupling
$T_{j}$ of the connection $\mathcal{R}_{j}$ taken separately: 
\begin{equation}
\begin{array}{c}
\textnormal{connection}\\
\textnormal{probabilities}\\
\textnormal{for }q_{j}
\end{array}=\left(\Pr\left[T_{j}=\left(l,\ldots,l\right)\right]:l\in\left\{ 1,\ldots,k\right\} \right).
\end{equation}

Since $S$ is a coupling of $\mathcal{R}$, we should have, for every
value $w^{i}$ of every bunch $R^{i}$, 
\begin{equation}
\sum_{v}\lambda^{i}\left(v,w^{i}\right)\gamma\left(v\right)=\Pr\left[R^{i}=w^{i}\right],\label{eq: bunch equations}
\end{equation}
where $\lambda^{i}\left(v,w^{i}\right)=1$ if $v^{i}=w^{i}$ (i.e.,
if the $i$th row of $v$, in reference to (\ref{eq: value v}), equals
$w^{i}$), and $\lambda^{i}\left(v,w^{i}\right)=0$ otherwise. Since
$S$ is a maximally connected coupling of $\mathcal{R}$, we should
have, for every value $w_{j}=\left(l,\ldots,l\right)$ of every maximal
coupling $T_{j}$,

\begin{equation}
\sum_{v}\lambda_{j}\left(v,w_{j}\right)\gamma\left(v\right)=\Pr\left[T_{j}=w_{j}=\left(l,\ldots,l\right)\right]\label{eq: connection equations}
\end{equation}
where $\lambda_{j}\left(v,w_{j}\right)=1$ if the $j$th column $v_{j}$
of $v$ in (\ref{eq: value v}) equals $w_{j}$, and $\lambda_{j}\left(v,w_{j}\right)=0$
otherwise. 

In we list the hidden outcomes $v$ in the same order as in the vector
$\mathbf{Q}$, the $1/0$ values of $\lambda^{i}\left(v,w^{i}\right)$
and $1/0$ values of $\lambda_{j}\left(v,w_{j}\right)$ in (\ref{eq: bunch equations})
and (\ref{eq: connection equations}) form rows of a Boolean matrix
$\mathbf{M}$, one row per each $\left(i,w^{i}\right)$ and each $\left(j,w_{j}\right)$,
such that (\ref{eq: bunch equations}) and (\ref{eq: connection equations})
can be written together as
\begin{equation}
\mathbf{MQ}=\mathbf{P}.\label{eq: associated eqs}
\end{equation}
We will refer to this matrix equation as the \emph{system of equations
associated with the \mbox{c-c} system} $\mathcal{R}$. In Section
\ref{sec: How-to-measure} below we will show that a vector of real
numbers $\mathbf{Q}$ satisfying this equation always exist. To form
a distribution for a maximally connected coupling $S$, however, $\mathbf{Q}$
also has to satisfy the following two constrains:
\begin{description}
\item [{(a)}] all components of $\mathbf{Q}$ are nonnegative, and
\item [{(b)}] they sum to 1.
\end{description}
The latter requirement is satisfied ``automatically.'' Indeed, by
construction, the rows of $\mathbf{M}$ corresponding to all possible
values of any given bunch have pairwise disjoint cells containing
1's: a hidden outcome $v$ in (\ref{eq: value v}) contains in its
$i$th row one and only one value of the $i$th bunch. This means
that if one adds all the rows of $\mathbf{M}$ corresponding to the
$i$th bunch one will get a row with 1's in all cells. The scalar
product of this row and $\mathbf{Q}$ equals both the sum of the elements
in $\mathbf{Q}$ and the sum of all bunch probabilities in the $i$th
bunch, which is 1.

The nonnegativity constraint, however, does not have to be satisfied:
it is possible that every one of the infinite set of solutions for
$\mathbf{Q}$ contains some negative components. This is the case
when the \mbox{c-c} system for which we have constructed the equations
is contextual.

We can formulate now the main statement of this section.
\begin{thm}
\label{thm: LP}A \mbox{c-c} system is noncontextual (i.e., it has
a maximally connected coupling) if and only if the associated system
of equations $\mathbf{MQ=P}$ has a solution for $\mathbf{Q}$ with
nonnegative components. Any such a solution defines a distribution
of the hidden outcomes of the coupling.
\end{thm}
The task of finding solutions for (\ref{eq: associated eqs}) subject
to the nonnegativity constraint is a linear programming task. It is
always well-defined and leads to an answer (an example of a solution
or the determination that it does not exist) in polynomial time with
respect to the number of the elements in $\mathbf{Q}$ (Karmarkar,
1984). This is all that matters to us theoretically. In practice,
some algorithms are more efficient than Karmarkar's in most cases
(e.g., the simplex algorithm).

The linear programming problem in Theorem \ref{thm: LP} is especially
transparent when all variables in a \mbox{c-c} system are binary
with the same possible values, say, 1 and -1. The reader may find
it useful to check, using Fig.~\ref{tab:The-Boolean-matrix}, all
the steps of the derivation of the linear equations (\ref{eq: associated eqs})
using the \mbox{c-c} system $\mathcal{A}$ of our opening example
(Fig.~\ref{fig:system 2by2}). Whether this system is contextual
depends on $\mathbf{P}$, specifically, on the bunch probabilities
in $\mathbf{P}$. Recall that the connection probabilities, the last
four elements of $\mathbf{P}$ in Fig. \ref{tab:The-Boolean-matrix},
are computed from the bunch probabilities using Theorem \ref{thm: maximal-coupling}.
Thus, if the bunch probabilities in $\mathbf{P}$ are as shown in
the upper panel of Fig.~\ref{fig: QQ system}, then the connection
probabilities should be as in the middle panel, and it can be shown
by applying a linear programming algorithm that the matrix equation
$\mathbf{MQ}=\mathbf{P}$ does not have a solution with nonnegative
elements. In this simple case we can confirm this result by a direct
observation of the internal contradiction in the maximally connected
coupling shown in the bottom panel.

\section{\label{sec: Cyclic-c-c-systems}Cyclic c-c systems}

The question we pose now is: is there a shortcut to find out if a
\mbox{c-c} system is contextual, without resorting to linear programming?
As it turns out, for an important class of so-called cyclic systems
with binary variables (Dzhafarov, Zhang, \& Kujala, 2015; Dzhafarov,
Kujala, \& Cervantes, 2016; Kujala, Dzhafarov, \& Larsson, 2015; Kujala
\& Dzhafarov, 2016) the answer to this question is affirmative.

\subsection{Contextuality criterion for cyclic c-c systems}

\begin{figure}
\begin{centering}
\includegraphics[scale=0.25]{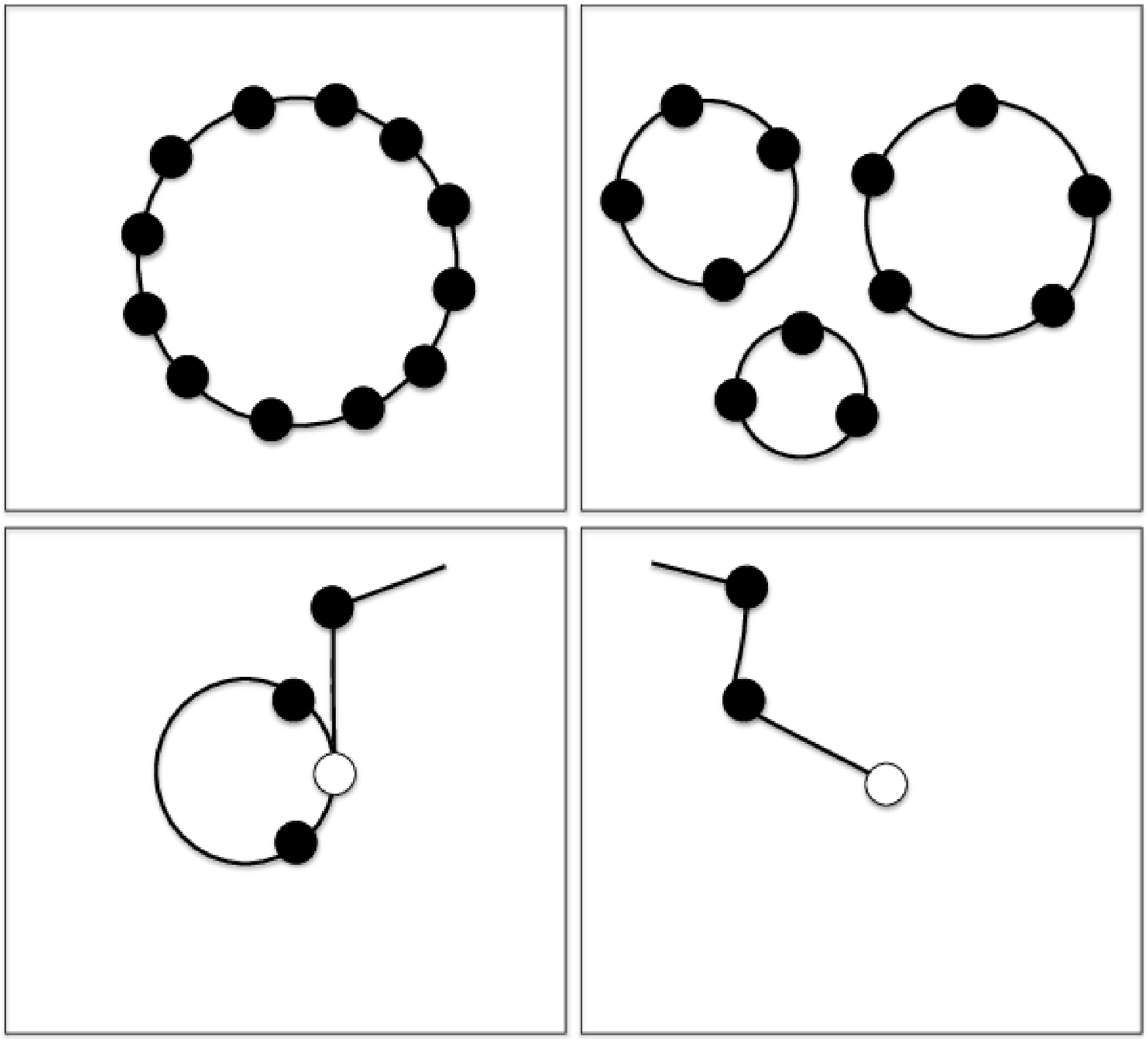}
\par\end{centering}

\caption{\label{fig: cyclicity proof}The conteNts (small circles) in a system
satisfying the conditions CYC1-CYC2 can be arranged in a cycle (upper
left panel) or several disjoint cycles (upper right panel), such that
any two adjacent conteNts define one conteXt. The proof is in the
lower panels: for the conteNts not to form cycles, some conteNt should
be placed in a position like the ones of the open circles in the two
lower panels. This is, however, impossible: in the left one it belongs
to more than two conteXts (identified with a pair of conteNts), in
the right one it belongs to a single conteXt.}
\end{figure}

A \emph{cyclic} system is defined as a \mbox{c-c} system in which
\begin{description}
\item [{(CYC1)}] each conteXt includes precisely two conteNts,
\item [{(CYC2)}] each conteNt is included in precisely two conteXts.
\end{description}
We will also assume that
\begin{description}
\item [{(CYC3)}] all random variables are binary with the same two possible
values (traditionally, $1$ and $-1$).
\end{description}
Fig.~\ref{fig: cyclicity proof} makes it clear why such a system
is called cyclic: to satisfy the properties above, the conteNts should
be arrangeable in one or more cycles in which a conteXt corresponds
to any two adjacent conteNts. If the conteNts are arranged into several
cycles, from the point of view of contextuality analysis each cycle
forms a separate system, with no information regarding one of them
being relevant for another's analysis. We will therefore, with no
loss of generality, assume that a cyclic system involves a single
cycle. 

The number of conteNts (or connections) in a cyclic system equals
the number of conteXts (or bunches) in it, and it is referred to as
the \emph{rank} of the system. The \mbox{c-c} matrix for the cyclic
system has the form shown in Fig.~\ref{fig: cyclic n}, generalizing
the matrices in Fig.~\ref{fig:system 2by2} (cyclic system of rank
2) and Fig.~\ref{fig: SZLG} (cyclic system of rank 3).

\begin{figure}
\begin{centering}
\begin{tabular}{|c|c|c|c|c|c|c|c}
\cline{1-7} 
$R_{1}^{1}$$\begin{array}{c}
\\
\\
\end{array}$ & $R_{2}^{1}$$\begin{array}{c}
\\
\\
\end{array}$ & $\cdot$$\begin{array}{c}
\\
\\
\end{array}$ & $\cdot$$\begin{array}{c}
\\
\\
\end{array}$ & $\cdots$$\begin{array}{c}
\\
\\
\end{array}$ & $\cdot$$\begin{array}{c}
\\
\\
\end{array}$ & $\cdot$$\begin{array}{c}
\\
\\
\end{array}$ & $c_{1}$\tabularnewline
\cline{1-7} 
$\cdot$$\begin{array}{c}
\\
\\
\end{array}$ & $R_{2}^{2}$$\begin{array}{c}
\\
\\
\end{array}$ & $R_{3}^{2}$$\begin{array}{c}
\\
\\
\end{array}$ & $\cdot$$\begin{array}{c}
\\
\\
\end{array}$ & $\cdots$$\begin{array}{c}
\\
\\
\end{array}$ & $\cdot$$\begin{array}{c}
\\
\\
\end{array}$ & $\cdot$$\begin{array}{c}
\\
\\
\end{array}$ & $c_{2}$\tabularnewline
\cline{1-7} 
$\cdot$$\begin{array}{c}
\\
\\
\end{array}$ & $\cdot$$\begin{array}{c}
\\
\\
\end{array}$ & $R_{3}^{3}$$\begin{array}{c}
\\
\\
\end{array}$ & $R_{4}^{3}$$\begin{array}{c}
\\
\\
\end{array}$ & $\cdots$$\begin{array}{c}
\\
\\
\end{array}$ & $\cdot$$\begin{array}{c}
\\
\\
\end{array}$ & $\cdot$$\begin{array}{c}
\\
\\
\end{array}$ & $c_{3}$\tabularnewline
\cline{1-7} 
$\vdots$$\begin{array}{c}
\\
\\
\end{array}$ & $\vdots$$\begin{array}{c}
\\
\\
\end{array}$ & $\vdots$$\begin{array}{c}
\\
\\
\end{array}$ & $\vdots$$\begin{array}{c}
\\
\\
\end{array}$ & $\iddots$$\begin{array}{c}
\\
\\
\end{array}$ & $\vdots$$\begin{array}{c}
\\
\\
\end{array}$ & $\vdots$$\begin{array}{c}
\\
\\
\end{array}$ & $\vdots$\tabularnewline
\cline{1-7} 
$\cdot$$\begin{array}{c}
\\
\\
\end{array}$ & $\cdot$$\begin{array}{c}
\\
\\
\end{array}$ & $\cdot$$\begin{array}{c}
\\
\\
\end{array}$ & $\cdot$$\begin{array}{c}
\\
\\
\end{array}$ & $\cdots$$\begin{array}{c}
\\
\\
\end{array}$ & $R_{n-1}^{n-1}$$\begin{array}{c}
\\
\\
\end{array}$ & $R_{n}^{n-1}$$\begin{array}{c}
\\
\\
\end{array}$ & $c_{n-1}$\tabularnewline
\cline{1-7} 
$R_{1}^{n}$$\begin{array}{c}
\\
\\
\end{array}$ & $\cdot$$\begin{array}{c}
\\
\\
\end{array}$ & $\cdot$$\begin{array}{c}
\\
\\
\end{array}$ & $\cdot$$\begin{array}{c}
\\
\\
\end{array}$ & $\cdots$$\begin{array}{c}
\\
\\
\end{array}$ & $\cdot$$\begin{array}{c}
\\
\\
\end{array}$ & $R_{n}^{n}$$\begin{array}{c}
\\
\\
\end{array}$ & $c_{n}$\tabularnewline
\cline{1-7} 
\multicolumn{1}{c}{$q_{1}$} & \multicolumn{1}{c}{$q_{2}$} & \multicolumn{1}{c}{$q_{3}$} & \multicolumn{1}{c}{$q_{4}$} & \multicolumn{1}{c}{$\cdots$} & \multicolumn{1}{c}{$q_{n-1}$} & \multicolumn{1}{c}{$q_{n}$} & $\boxed{\boxed{\mathsf{CYC}}}$\tabularnewline
\end{tabular}
\par\end{centering}

\caption{\label{fig: cyclic n}The \mbox{c-c} matrix for a cyclic system of
an arbitrary rank $n$ (shown here for a sufficiently large $n$,
but $n$ can be as small as 2 or 3, shown in Figs. \ref{fig:system 2by2}
and \ref{fig: SZLG}, respectively). Each conteXt includes two conteNts,
and each conteNt is included in two conteXts. Bunches of the system
are formed by the pairs of random variables $\left(R_{1}^{1},R_{2}^{1}\right),\left(R_{2}^{2},R_{3}^{2}\right),\ldots,\left(R_{n-1}^{n-1},R_{n}^{n-1}\right),\left(R_{n}^{n},R_{1}^{n}\right)$,
whereas the connections of the system are formed by the pairs $\left(R_{1}^{1},R_{1}^{n}\right),\left(R_{2}^{2},R_{2}^{1}\right),\ldots,\left(R_{n-1}^{n-1},R_{n-1}^{n-2}\right),\left(R_{n}^{n},R_{n}^{n-1}\right)$.
All random variables are binary, with possible values denoted $+1$
and $-1$. }
\end{figure}

In the presentation below we use $\left\langle X\right\rangle $ to
denote the expected value of a random variable $X$ with possible
values $+1$ and $-1$:
\begin{equation}
\left\langle X\right\rangle =\Pr\left[X=1\right]-\Pr\left[X=-1\right].
\end{equation}

Given $k>0$ real numbers $x_{1},\ldots,x_{k}$, we define the function
\begin{equation}
\sodd\left(x_{1},\ldots,x_{k}\right)=\max_{\left(\iota_{1},\ldots,\iota_{k}\right)\in\left\{ -1,1\right\} ^{k}:\prod_{i=1}^{k}\iota_{i}=-1}\left(\sum_{\begin{array}{c}
i=1\end{array}}^{k}\iota_{i}x_{i}\right).
\end{equation}
This means that one takes each argument $x_{i}$ in the sum with either
$+$ sign or $-$ sign, tries all combinations in which the number
of minuses is odd, and chooses the largest sum. For example, 
\[
\begin{array}{c}
\sodd\left(5,6\right)=-5+6,\\
\\
\sodd\left(5,-6\right)=5-\left(-6\right),\\
\\
\sodd\left(1,2,-3,-10,100\right)=-1+2-\left(-3\right)-\left(-10\right)+100.
\end{array}
\]
 Finally, we have to introduce the cyclic addition and subtraction
operations, $\oplus1$ and $\ominus1$: if the numbers $1,\ldots,n$
are arranged circularly like on a clock dial, then $\oplus1$ and
$\ominus1$ mean, respectively, clockwise and counterclockwise shift
to the next position. The only difference of these operations from
the usual $+1$ and $-1$ is that $n\oplus1=1$ and $1\ominus1=n$. 

Now we can formulate a criterion of (i.e., a necessary and sufficient
condition for) contextuality of a cyclic system. 
\begin{thm}[Kujala and Dzhafarov, 2016]
A cyclic system of rank $n$ is noncontextual if and only if
\begin{equation}
\sodd\left(\left\langle R_{i}^{i}R_{i\oplus1}^{i}\right\rangle :i=1,\ldots,n\right)\leq n-2+\sum_{i=1}^{n}\left|\left\langle R_{i}^{i}\right\rangle -\left\langle R_{i}^{i\ominus1}\right\rangle \right|.\label{eq: criterion gen}
\end{equation}

\end{thm}
In the left-hand side expression, the arguments of the function $\sodd$
are the expected products for the $n$ bunches of the system: $\left\langle R_{1}^{1}R_{2}^{1}\right\rangle $,
$\left\langle R_{2}^{2}R_{3}^{2}\right\rangle $, etc., the last one,
due to the cyclicality, being $\left\langle R_{n}^{n}R_{1}^{n}\right\rangle $.
In the right-hand side of the inequality, the summation sign operates
over the $n$ connections of the system: for each connection, $\left(R_{1}^{1},R_{1}^{n}\right),\left(R_{2}^{2},R_{2}^{1}\right),\ldots,\left(R_{n}^{n},R_{n}^{n-1}\right)$,
we take the distance between the expectations of its elements. If
the system is consistently connected, all these distances are zero,
and the criterion acquires the form
\begin{equation}
\sodd\left(\left\langle R_{i}^{i}R_{i\oplus1}^{i}\right\rangle :i=1,\ldots,n\right)\leq n-2.\label{eq: criterion cons}
\end{equation}

\subsection{Examples of cyclic systems}

It has been mentioned in the introduction that cyclic systems of rank
2 have been prominently studied in a behavioral setting, in the paradigm
where the conteNts are two Yes/No questions and conteXts are defined
by two orders in which these questions are asked. The noncontextuality
criterion (\ref{eq: criterion gen}) for $n=2$ acquires the form
\begin{equation}
\left|\left\langle R_{1}^{1}R_{2}^{1}\right\rangle -\left\langle R_{2}^{2}R_{1}^{2}\right\rangle \right|\leq\left|\left\langle R_{1}^{1}\right\rangle -\left\langle R_{1}^{2}\right\rangle \right|+\left|\left\langle R_{2}^{1}\right\rangle -\left\langle R_{2}^{2}\right\rangle \right|.\label{eq: QQ}
\end{equation}
It is known (Moore, 2002) that the distributions of responses to the
same question depend on whether the question is asked first or second.
In our terminology, this means that the system is inconsistently connected,
and the right-hand side of the inequality above is greater than zero.
At the same time, as Wang and Busemeyer (2013) have discovered in
their analysis of a large body of question pairs, the probability
with which the answer to the two questions is one and the same does
not depend on the order in which they are asked. To the extent this
generalization holds, it means that the left-hand side of the inequality
(\ref{eq: QQ}) is zero. In turn, this means that the system describing
responses to two questions asked in two orders cannot be contextual
(see Dzhafarov, Zhang, Kujala, 2015, for a detailed discussion).

\begin{figure}
\begin{centering}
\begin{tabular}{|c|c|c|c|c}
\cline{1-4} 
$R_{1}^{1}$$\begin{array}{c}
\\
\\
\end{array}$ & $R_{2}^{1}$$\begin{array}{c}
\\
\\
\end{array}$ & $\cdot$$\begin{array}{c}
\\
\\
\end{array}$ & $\cdot$$\begin{array}{c}
\\
\\
\end{array}$ & $c_{1}$\tabularnewline
\cline{1-4} 
$\cdot$$\begin{array}{c}
\\
\\
\end{array}$ & $R_{2}^{2}$$\begin{array}{c}
\\
\\
\end{array}$ & $R_{3}^{2}$$\begin{array}{c}
\\
\\
\end{array}$ & $\cdot$$\begin{array}{c}
\\
\\
\end{array}$ & $c_{2}$\tabularnewline
\cline{1-4} 
$\cdot$$\begin{array}{c}
\\
\\
\end{array}$ & $\cdot$$\begin{array}{c}
\\
\\
\end{array}$ & $R_{3}^{3}$$\begin{array}{c}
\\
\\
\end{array}$ & $R_{4}^{3}$$\begin{array}{c}
\\
\\
\end{array}$ & $c_{3}$\tabularnewline
\cline{1-4} 
$R_{1}^{4}$$\begin{array}{c}
\\
\\
\end{array}$ & $\cdot$$\begin{array}{c}
\\
\\
\end{array}$ & $\cdot$$\begin{array}{c}
\\
\\
\end{array}$ & $R_{4}^{4}$$\begin{array}{c}
\\
\\
\end{array}$ & $c_{4}$\tabularnewline
\cline{1-4} 
\multicolumn{1}{c}{$q_{1}$} & \multicolumn{1}{c}{$q_{2}$} & \multicolumn{1}{c}{$q_{3}$} & \multicolumn{1}{c}{$q_{4}$} & $\boxed{\boxed{\mathsf{EPR-B}}}$\tabularnewline
\end{tabular}
\par\end{centering}

\begin{centering}
\includegraphics[scale=0.4]{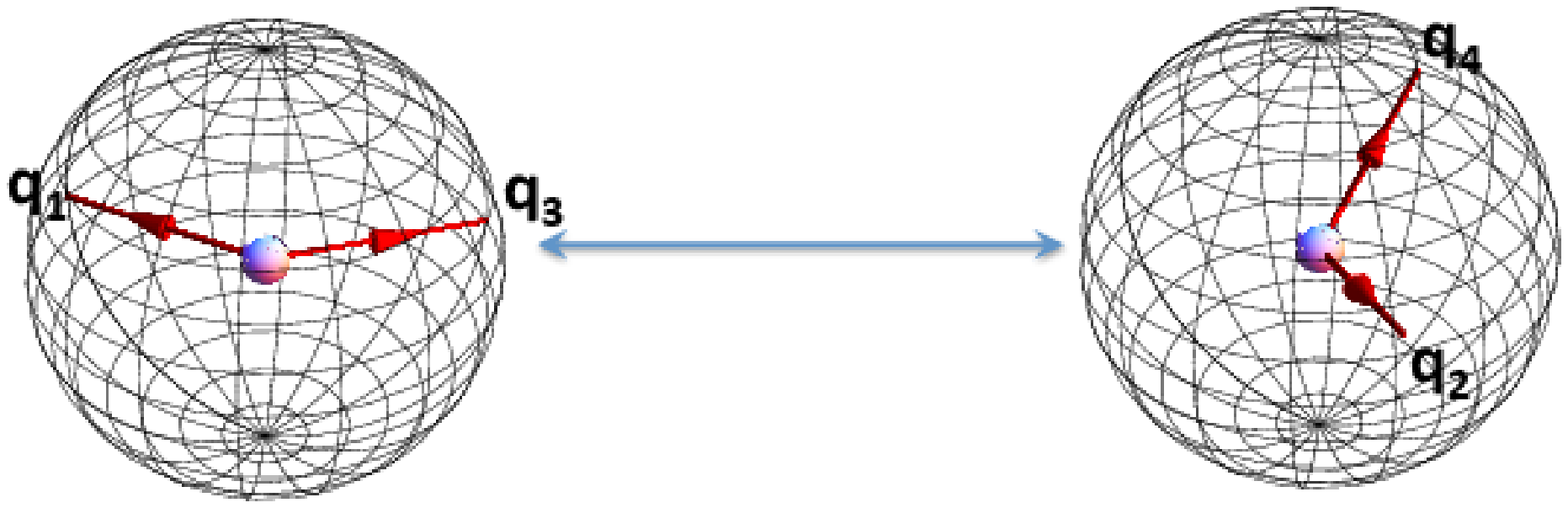}
\par\end{centering}

\caption{\label{fig:Einstein-Podolsky-Rosen-Bohm-sys}Einstein-Podolsky-Rosen-Bohm-system
(Bohm \& Aharonov, 1957), for which the celebrated Bell-CHSH inequalities
were derived (Bell, 1964, 1966; Clauser et al., 1969; Clauser \& Horne,
1974). This is a cyclic system of rank 4. Two spin-$\nicefrac{1}{2}$
particles (e.g., electrons) are created in what is called a ``singlet
state'' and move away from each other while remaining entangled.
Alice chooses one of the two axes denoted $1$ and $3$ and measures
the spin of the left particle, ``up'' ($+1$) or ``down'' ($-1$).
Bob chooses one of the two axes denoted $2$ and $4$ and measures
the spin of the right particle, $+1$ or $-1$. The conteNts here
are the (choices of the) four axes, $q_{1},q_{2},q_{3},q_{4}$, the
conteXts are defined by the pairs of the axes chosen, one by Alice
and another by Bob. The system is contextual for some combinations
of the four axes.}
\end{figure}

Perhaps the best known cyclic system is one of rank 4, whose quantum-mechanical
version is shown in Fig.~\ref{fig:Einstein-Podolsky-Rosen-Bohm-sys}.
According to the laws of quantum mechanics, the product expectation
$\left\langle R_{i}^{i}R_{i\oplus1}^{i}\right\rangle $ for Alice's
choice of axis $q_{i}$ and Bob's choice of axis $q_{i\oplus1}$ equals
$-\cos\theta_{i,i\oplus1}$, where $\theta_{i,i\oplus1}$ denotes
the angle between the two axes. Assume, e.g., that the four axes are
coplanar, and form the following angles with respect to some fixed
direction
\begin{equation}
\begin{array}{cccc}
q_{1} & q_{2} & q_{3} & q_{4}\\
0 & \nicefrac{\pi}{4} & \nicefrac{\pi}{2} & -\nicefrac{\pi}{4}
\end{array}.
\end{equation}
The calculation yields in this case 
\begin{equation}
\sodd\left(\left\langle R_{i}^{i}R_{i\oplus1}^{i}\right\rangle :i=1,2,3,4\right)=2\sqrt{2}.
\end{equation}
If any possibility of direct interaction between Alice and Bob is
excluded, i.e., Alice's measurements are not influenced by Bob's choices
of his axes and Bob's measurements are not influenced by Alice's choices
of her axes, and if we exclude any possibility of misrecording, then
\begin{equation}
\sum_{i=1}^{4}\left|\left\langle R_{i}^{i}\right\rangle -\left\langle R_{i}^{i\ominus1}\right\rangle \right|=0,
\end{equation}
and the inequality (\ref{eq: criterion gen}) acquires the form of
the inequality (\ref{eq: criterion cons}), for $n=4$. The value
$2\sqrt{2}$ for the left-hand side expression violates this inequality,
indicating that the system is contextual.

\begin{figure}
\begin{centering}
\includegraphics[bb=0bp 100bp 720bp 580bp,scale=0.25]{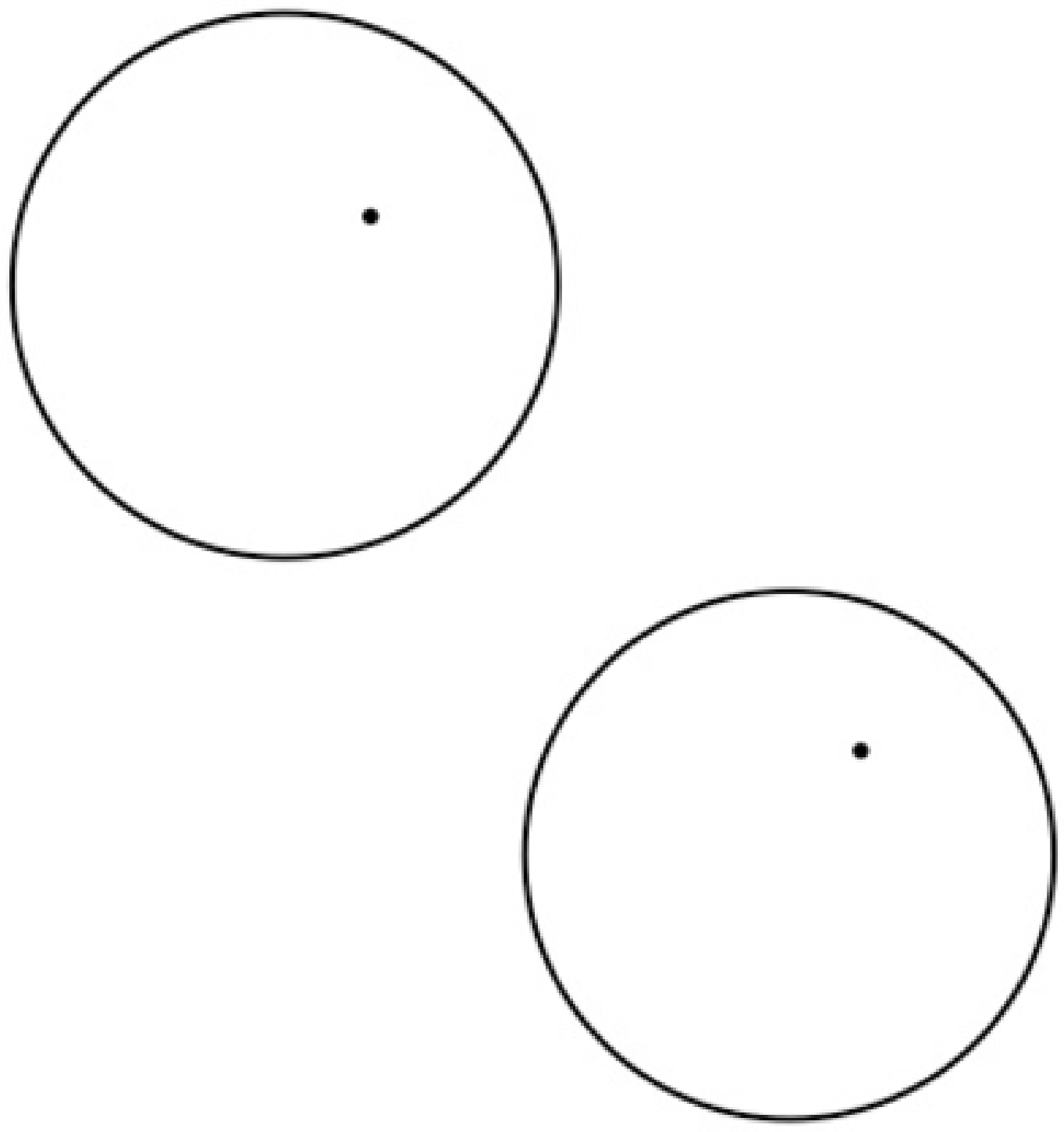}
\par\end{centering}

\caption{\label{fig:A-matching-experiment}A matching experiment: a participant
rotated a trackball that controlled the position of the dot within
a lower-right circle until she judged it to match the fixed position
of the target dot in the upper-left circle. The positions were described
in polar coordinates, and the target position could have one of two
radius values $q_{1}$ or $q_{3}$ combined with one of two polar
angles, $q_{2}$ or $q_{4}$. In each of the four conteXts, $\left(q_{i},q_{i\oplus1}\right)$,
$i=1,2,3,4$, the adjusted dot's position was described by two polar
coordinates (random variables) $Rad_{i,i\oplus1}$ and $Ang{}_{i,i\oplus1}$,
that were then dichotomized to create a cyclic system of rank 4.}
\end{figure}

There were several studies of systems having the cyclic rank 4 structure
in behavioral settings. Thus, Fig.~\ref{fig:A-matching-experiment}
describes one of the psychophysical matching experiments analyzed
in Dzhafarov, Ru, and Kujala (2015). The dichotomization of the response
variables was done as follows: we choose radial length values $rad_{1},rad_{3}$
(they may be the same) and polar angle values $ang_{2}$ and $ang_{4}$
(they also may be the same), and we define 
\begin{equation}
\begin{array}{l}
R_{i}^{i}=\left\{ \begin{array}{ccc}
+1 & \textnormal{if} & Rad_{i,i\oplus1}>rad_{i}\\
\\
-1 & \textnormal{if} & Rad_{i,i\oplus1}\leq rad_{i}
\end{array}\right.,\\
\\
R_{i\oplus1}^{i}=\left\{ \begin{array}{ccc}
+1 & if & Ang_{i,i\oplus1}>ang_{i\oplus1}\\
\\
-1 & if & Ang_{i,i\oplus1}\leq ang_{i\oplus1}
\end{array}\right.
\end{array}
\end{equation}
for $i=1,3$, and 
\begin{equation}
\begin{array}{l}
R_{i}^{i}=\left\{ \begin{array}{ccc}
+1 & if & Ang_{i,i\oplus1}>ang_{i}\\
\\
-1 & if & Ang_{i,i\oplus1}\leq ang_{i}
\end{array}\right.,\\
\\
R_{i\oplus1}^{i}=\left\{ \begin{array}{ccc}
+1 & \textnormal{if} & Rad_{i,i\oplus1}>rad_{i\oplus1}\\
\\
-1 & \textnormal{if} & Rad_{i,i\oplus1}\leq rad_{i\oplus1}
\end{array}\right.
\end{array}
\end{equation}
for $i=2,4.$ The parameters $rad_{1},rad_{3}$ and $ang_{2},ang_{4}$
can be chosen in multiple ways, paralleling various sets of four axes
in the Alice-Bob quantum-mechanical experiment. 

See Dzhafarov, Ru, and Kujala (2015) for other examples of behavioral
rank 4 cyclic systems. That paper also reviews a behavioral experiment
with a cyclic system of rank 3. Cyclic systems of rank 5 play an important
role in quantum theory (Klyachko et al., 2008). For the contextuality
analysis of an experiment designed to test (a special form of) the
inequality (\ref{eq: criterion gen}) for $n=5$ (Lapkiewitz et al.,
2011), see Kujala, Dzhafarov, and Larsson (2015).

\section{\label{sec: How-to-measure}How to measure degree of contextuality}

Intuitively, some contextual systems are more contextual than others.
For instance, a cyclic system of rank $n$ can violate the inequality
(\ref{eq: criterion gen}) ``grossly'' or ``slightly,'' and in
the latter case it may be considered less contextual. The question
we pose now is: if a \mbox{c-c} system is contextual, is there a
\emph{principled} way to measure the degree of contextuality in it?
The emphasis is on the qualifier ``principled,'' as one can easily
come up with various ad hoc measures for special types of \mbox{c-c}
systems.

In this section we describe one way of constructing such a measure.
It uses the notion of \emph{quasi-probability distributions} that
differ from the proper ones in that the probability masses in them
are replaced with arbitrary, possibly negative, real numbers summing
to unity. This conceptual tool has been previously used to deal with
contextuality in consistently connected systems (Abramsky \& Brandenburger,
2011; Al-Safi \& Short, 2013). A measure of contextuality based on
the notion of quasi-probability distributions, also for consistently
connected systems, was proposed by de Barros and Oas (2014) and investigated
in de Barros, Oas, and Suppes (2015) and de Barros et al. (2015).
Our measure is a generalization of the de Barros-Oas measure to arbitrary
\mbox{c-c} systems. Another generalization and a different way of
using quasi-probability distributions to measure contextuality in
arbitrary \mbox{c-c} systems was recently proposed by Kujala (2016).
This measure requires a modification of CbD and will not be discussed
here. 

\begin{figure*}
\begin{footnotesize}%
\begin{tabular}{|c|c|c|c|c|c|c|c|c|c|c|c|c|c|c|c|c|c|}
\cline{1-16} 
\multicolumn{16}{|c|}{$\mathbf{M}^{*}$} & \multicolumn{1}{c}{\multirow{11}{*}{}} & \multicolumn{1}{c}{}\tabularnewline
$\!\begin{array}{cc}
+ & +\\
+ & +
\end{array}\!$ & $\!\begin{array}{cc}
+ & +\\
+ & -
\end{array}\!$ & $\!\begin{array}{cc}
+ & +\\
- & +
\end{array}\!$ & $\!\begin{array}{cc}
+ & +\\
- & -
\end{array}\!$ & $\!\begin{array}{cc}
+ & -\\
+ & +
\end{array}\!$ & $\!\begin{array}{cc}
+ & -\\
+ & -
\end{array}\!$ & $\!\begin{array}{cc}
+ & -\\
- & +
\end{array}\!$ & $\!\begin{array}{cc}
+ & -\\
- & -
\end{array}\!$ & $\!\begin{array}{cc}
- & +\\
+ & +
\end{array}\!$ & $\!\begin{array}{cc}
- & +\\
+ & -
\end{array}\!$ & $\!\begin{array}{cc}
- & +\\
- & +
\end{array}\!$ & $\!\begin{array}{cc}
- & +\\
- & -
\end{array}\!$ & $\!\begin{array}{cc}
- & -\\
+ & +
\end{array}\!$ & $\!\begin{array}{cc}
- & -\\
+ & -
\end{array}\!$ & $\!\begin{array}{cc}
- & -\\
- & +
\end{array}\!$ & $\!\begin{array}{cc}
- & -\\
- & -
\end{array}\!$ &  & \multicolumn{1}{c}{$\mathbf{P}^{*}$}\tabularnewline
\cline{1-16} \cline{18-18} 
1 & 1 & 1 & 1 & 1 & 1 & 1 & 1 & 1 & 1 & 1 & 1 & 1 & 1 & 1 & 1 &  & $\begin{array}{l}
\\
\\
\end{array}\Pr\left[\right]=1$\tabularnewline
\cline{1-16} \cline{18-18} 
1 & 1 & 1 & 1 & 1 & 1 & 1 & 1 &  &  &  &  &  &  &  &  &  & $\begin{array}{l}
\Pr\left[S_{1}^{1}=+\right]\\
=\Pr\left[R_{1}^{1}=+\right]
\end{array}$\tabularnewline
\cline{1-16} \cline{18-18} 
1 & 1 &  &  & 1 & 1 &  &  & 1 & 1 &  &  & 1 & 1 &  &  &  & $\begin{array}{l}
\Pr\left[S_{2}^{1}=+\right]\\
=\Pr\left[R_{2}^{1}=+\right]
\end{array}$\tabularnewline
\cline{1-16} \cline{18-18} 
1 &  & 1 &  & 1 &  & 1 &  & 1 &  & 1 &  & 1 &  & 1 &  &  & $\begin{array}{l}
\Pr\left[S_{1}^{2}=+\right]\\
=\Pr\left[R_{1}^{2}=+\right]
\end{array}$\tabularnewline
\cline{1-16} \cline{18-18} 
1 &  & 1 &  & 1 &  & 1 &  & 1 &  & 1 &  & 1 &  & 1 &  &  & $\begin{array}{l}
\Pr\left[S_{2}^{2}=+\right]\\
=\Pr\left[R_{2}^{2}=+\right]
\end{array}$\tabularnewline
\cline{1-16} \cline{18-18} 
1 & 1 & 1 & 1 &  &  &  &  &  &  &  &  &  &  &  &  &  & $\begin{array}{l}
\Pr\left[S_{1}^{1}=+,S_{2}^{1}=+\right]\\
=\Pr\left[R_{1}^{1}=+,R_{2}^{1}=+\right]
\end{array}$\tabularnewline
\cline{1-16} \cline{18-18} 
1 &  &  &  & 1 &  &  &  & 1 &  &  &  & 1 &  &  &  &  & $\begin{array}{l}
\Pr\left[S_{1}^{2}=+,S_{2}^{2}=+\right]\\
=\Pr\left[R_{1}^{2}=+,R_{2}^{2}=+\right]
\end{array}$\tabularnewline
\cline{1-16} \cline{18-18} 
1 & 1 &  &  & 1 & 1 &  &  &  &  &  &  &  &  &  &  &  & $\begin{array}{l}
\Pr\left[S_{1}^{1}=+,S_{1}^{2}=+\right]\\
=\Pr\left[T_{1}^{1}=+,T_{1}^{2}=+\right]
\end{array}$\tabularnewline
\cline{1-16} \cline{18-18} 
1 &  & 1 &  &  &  &  &  & 1 &  & 1 &  &  &  &  &  &  & $\begin{array}{l}
\Pr\left[S_{2}^{1}=+,S_{2}^{2}=+\right]\\
=\Pr\left[T_{2}^{1}=+,T_{2}^{2}=+\right]
\end{array}$\tabularnewline
\cline{1-16} \cline{18-18} 
\end{tabular}\end{footnotesize}

\caption{\label{tab:The-Boolean-matrix-1}The Boolean matrix $\mathbf{M}^{*}$
and vector $\mathbf{P}^{*}$ for the \mbox{c-c} system $\mathcal{A}$
in Fig.~\ref{fig:system 2by2} (cyclic system of rank 2). The first
row of $\mathbf{M}^{*}$, formally, corresponds to the ``0-marginal''
probability in $\mathbf{P}$, and its role is to ensure that the elements
of $\mathbf{Q}$ in $\mathbf{M^{*}Q}=\mathbf{P^{*}}$ sum to 1. The
next 6 rows are modified bunch probabilities and the last 2 rows are
connection probabilities (``expanded,'' but in this case nominally
only). }
\end{figure*}

\begin{figure*}
\begin{footnotesize}

\begin{tabular}{|c|c|c|c|c|c|c|c|c|c|c|c|c|c|c|c|}
\hline 
\multicolumn{16}{|c|}{$\mathbf{M}^{*}$}\tabularnewline
$\!\begin{array}{cc}
+ & +\\
+ & +
\end{array}\!$ & $\!\begin{array}{cc}
+ & +\\
+ & -
\end{array}\!$ & $\!\begin{array}{cc}
+ & +\\
- & +
\end{array}\!$ & $\!\begin{array}{cc}
+ & +\\
- & -
\end{array}\!$ & $\!\begin{array}{cc}
+ & -\\
+ & +
\end{array}\!$ & $\!\begin{array}{cc}
+ & -\\
+ & -
\end{array}\!$ & $\!\begin{array}{cc}
+ & -\\
- & +
\end{array}\!$ & $\!\begin{array}{cc}
+ & -\\
- & -
\end{array}\!$ & $\!\begin{array}{cc}
- & +\\
+ & +
\end{array}\!$ & $\!\begin{array}{cc}
- & +\\
+ & -
\end{array}\!$ & $\!\begin{array}{cc}
- & +\\
- & +
\end{array}\!$ & $\!\begin{array}{cc}
- & +\\
- & -
\end{array}\!$ & $\!\begin{array}{cc}
- & -\\
+ & +
\end{array}\!$ & $\!\begin{array}{cc}
- & -\\
+ & -
\end{array}\!$ & $\!\begin{array}{cc}
- & -\\
- & +
\end{array}\!$ & $\!\begin{array}{cc}
- & -\\
- & -
\end{array}\!$\tabularnewline
\hline 
\hline 
1 & 1 & 1 & 1 & 1 & 1 & 1 & 1 & 1 & 1 & 1 & 1 & 1 & 1 & 1 & 1\tabularnewline
\hline 
1 & 1 & 1 & 1 & 1 & 1 & 1 & 1 &  &  &  &  &  &  &  & \tabularnewline
\hline 
1 & 1 &  &  & 1 & 1 &  &  & 1 & 1 &  &  & 1 & 1 &  & \tabularnewline
\hline 
1 &  & 1 &  & 1 &  & 1 &  & 1 &  & 1 &  & 1 &  & 1 & \tabularnewline
\hline 
1 &  & 1 &  & 1 &  & 1 &  & 1 &  & 1 &  & 1 &  & 1 & \tabularnewline
\hline 
1 & 1 & 1 & 1 &  &  &  &  &  &  &  &  &  &  &  & \tabularnewline
\hline 
1 &  &  &  & 1 &  &  &  & 1 &  &  &  & 1 &  &  & \tabularnewline
\hline 
1 & 1 &  &  & 1 & 1 &  &  &  &  &  &  &  &  &  & \tabularnewline
\hline 
1 &  & 1 &  &  &  &  &  & 1 &  & 1 &  &  &  &  & \tabularnewline
\hline 
\end{tabular}%
\begin{tabular}{|c|}
\hline 
$\mathbf{Q}$\tabularnewline
\hline 
\hline 
0\tabularnewline
\hline 
0\tabularnewline
\hline 
0\tabularnewline
\hline 
$\nicefrac{1}{2}$\tabularnewline
\hline 
0\tabularnewline
\hline 
$\nicefrac{1}{2}$\tabularnewline
\hline 
0\tabularnewline
\hline 
$-\nicefrac{1}{2}$\tabularnewline
\hline 
0\tabularnewline
\hline 
0\tabularnewline
\hline 
$\nicefrac{1}{2}$\tabularnewline
\hline 
$-\nicefrac{1}{2}$\tabularnewline
\hline 
0\tabularnewline
\hline 
0\tabularnewline
\hline 
0\tabularnewline
\hline 
$\nicefrac{1}{2}$\tabularnewline
\hline 
\end{tabular}%
\begin{tabular}{c}
\multirow{11}{*}{}\tabularnewline
\tabularnewline
\tabularnewline
\tabularnewline
\tabularnewline
\tabularnewline
\tabularnewline
\tabularnewline
\tabularnewline
\tabularnewline
\tabularnewline
\end{tabular}%
\begin{tabular}{|c|}
\multicolumn{1}{c}{}\tabularnewline
\multicolumn{1}{c}{$\mathbf{P}^{*}$}\tabularnewline
\hline 
$1$\tabularnewline
\hline 
$\nicefrac{1}{2}$\tabularnewline
\hline 
$\nicefrac{1}{2}$\tabularnewline
\hline 
$\nicefrac{1}{2}$\tabularnewline
\hline 
$\nicefrac{1}{2}$\tabularnewline
\hline 
$\nicefrac{1}{2}$\tabularnewline
\hline 
$0$\tabularnewline
\hline 
$\nicefrac{1}{2}$\tabularnewline
\hline 
$\nicefrac{1}{2}$\tabularnewline
\hline 
\end{tabular}\end{footnotesize}

\caption{\label{tab:The-Boolean-matrix-2}The same as Fig.~\ref{tab:The-Boolean-matrix-1},
but with the numerical values of the bunch and connection probabilities
of the contextual system shown in Fig.~\ref{fig: QQ system}. The
inserted column $\mathbf{Q}$ of quasi-probabilities is a solution
for $\mathbf{M}^{*}\mathbf{Q}=\mathbf{P^{*}}$ .}
\end{figure*}

\subsection{Dropping the nonnegativity constraint}

We have seen that a \mbox{c-c} system is contextual if and only if
the associated system of linear equations $\mathbf{MQ}=\mathbf{P}$
does not have a solution for $\mathbf{Q}$ with nonnegative components.
In this section we show that if we drop the nonnegativity constraint
the system of linear equations always has a solution (and generally
an infinity of them). Any such a solution assigns real numbers to
all hidden outcomes of the hypothetical maximally connected coupling.
Some of these numbers can be negative and some may exceed 1, but they
sum to 1. 

The existence of $\mathbf{Q}$ solving the linear equations $\mathbf{MQ}=\mathbf{P}$
follows from the existence of $\mathbf{Q}$ solving another system
of linear equations, 
\begin{equation}
\mathbf{M^{*}Q}=\mathbf{P}^{*}.
\end{equation}
We will refer to it as a \emph{modified-and-expanded} system of linear
equations associated with a \mbox{c-c} system. The term reflects
the fact that in $\mathbf{P}^{*}$ as compared to $\mathbf{P}$ the
bunch probabilities are presented in a modified form, and the connection
probabilities are expanded to specify entire distributions of the
(maximal) couplings of all connections. The rows of the Boolean matrix
$\mathbf{M}^{*}$ as compared to $\mathbf{M}$ change accordingly,
although its columns remain corresponding to the hidden outcomes $v$
ordered in the same way as $\gamma\left(v\right)$ are ordered in
$\mathbf{Q}$. 

The construction of $\mathbf{P}^{*}$ and $\mathbf{M}^{*}$ consists
of three parts.

\textbf{Part 1: first row.} The first element of $\mathbf{P}^{*}$
is 1, and the first row of $\mathbf{M}^{*}$ is filled with $1$'s.
This choice ensures
\begin{equation}
\sum_{v}\gamma\left(v\right)=1.
\end{equation}

\textbf{Part 2: bunch probabilities.} Next we include in $\mathbf{P}^{*}$
the $1$-marginal probabilities $\Pr\left[R_{j}^{i}=l\right]$ for
all random variables $R_{j}^{i}$ and all $l=1,\ldots,k-1$. The value
$l=k$ is excluded because $\Pr\left[R_{j}^{i}=k\right]$ is uniquely
determined as a linear combination of the probabilities already included.
The row of $\mathbf{M}^{*}$ corresponding to $\Pr\left[R_{j}^{i}=l\right]$
(i.e., the row whose scalar product by $\mathbf{Q}$ yields this probability)
has 1's in the cells for $v$ with $v_{j}^{i}=l$, and it has zeros
in other cells. 

The next set of elements of $\mathbf{P}^{*}$ are $2$-marginal probabilities
$\Pr\left[R_{j}^{i}=l,R_{j'}^{i}=l'\right]$ for all pairs of random
variables $R_{j}^{i},R_{j'}^{i}$ ($j<j'$) and $\left(l,l'\right)\in\left\{ 1,\ldots,k-1\right\} ^{2}$.
The $2$-marginal probabilities for $R_{j}^{i}=k$ or $R_{j'}^{i}=k$
are excluded because they are uniquely determinable as linear combinations
of the probabilities already included. The row of $\mathbf{M}^{*}$
corresponding to $\Pr\left[R_{j}^{i}=l,R_{j'}^{i}=l'\right]$ has
1's in the cells for $v$ with $v_{j}^{i}=l,v_{j'}^{i}=l'$, and it
has zeros in other cells.

Proceeding in this manner, we include in $\mathbf{P}^{*}$ the $r$-marginal
probabilities $\Pr\left[R_{j_{1}}^{i}=l_{1},\ldots,R_{j_{r}}^{i}=l_{r}\right]$
($j_{1}<\ldots<j_{r}$) for all bunches $R^{i}$ with at least $r$
distinct random variables, and for all $\left(l_{1},\ldots,l_{r}\right)\in\left\{ 1,\ldots,k-1\right\} ^{r}$.
We exclude all probabilities involving the value $k$ for at least
one of the random variables in the $r$-marginal. The row in $\mathbf{M}^{*}$
corresponding to $\Pr\left[R_{j_{1}}^{i}=l_{1},\ldots,R_{j_{r}}^{i}=l_{r}\right]$
has 1's in the cells for $v$ with $v_{j_{1}}^{i}=l_{1},\ldots,v_{j_{r}}^{i}=l_{r}$,
and it has zeros in other cells. The procedure stops at the smallest
$r$ such that no bunches in the system contain $r$ distinct random
variables.

\textbf{Part 3: connection probabilities.} This part of the construction
deals with the maximal couplings for connections. By Theorem \ref{thm: maximal-coupling},
a maximal coupling $T_{j}$ for a connection $\mathcal{R}_{j}$ always
exists, i.e., one can always find the joint probabilities $\Pr\left[T_{j}=v_{j}\right]$
for all $v_{j}$, so that 
\begin{equation}
\Pr\left[T_{j}^{i}=v_{j}^{i}\right]=\Pr\left[R_{j}^{i}=v_{j}^{i}\right]
\end{equation}
for all $R_{j}^{i}$ in $\mathcal{R}_{j}$, and 
\begin{equation}
\Pr\left[T_{j}=\left(l,\ldots,l\right)\right]=\min_{\textnormal{components }T_{j}^{i}\textnormal{ of }T_{j}}\Pr\left[T_{j}^{i}=l\right]\label{eq: maximality constarint}
\end{equation}
 for $l=1,\ldots,k$. We choose any maximal coupling $T_{j}$ for
each connection $\mathcal{R}_{j}$, and we treat it as if it were
an observed bunch. This allows us to repeat on $T_{j}$ the procedure
of Part 2, except that the $1$-marginal probabilities $\Pr\left[T_{j}^{i}=v_{j}^{i}\right]=\Pr\left[R_{j}^{i}=l\right]$,
which are the same for both bunches and connections, have already
been included in $\mathbf{P}^{*}$ (and the corresponding rows in
$M^{*}$ have been formed). We add, however, all higher-order marginals
$\Pr\left[T_{j}^{i_{1}}=l_{1},\ldots,T_{j}^{i_{r}}=l_{r}\right]$
($i_{1}<\ldots<i_{r}$) for all connections $R$$_{j}$ with at least
$r$ distinct random variables, for all $\left(l_{1},\ldots,l_{r}\right)\in\left\{ 1,\ldots,k-1\right\} ^{r}$.
The row in $\mathbf{M}^{*}$ corresponding to $\Pr\left[T_{j}^{i_{1}}=l_{1},\ldots,T_{j}^{i_{r}}=l_{r}\right]$
has 1's in the cells for $v$ with $v_{j}^{i_{1}}=l_{1},\ldots,v_{j}^{i_{r}}=l_{r}$,
and it has zeros in other cells. We apply this procedure to $r=2,3,\ldots$
until we reach the smallest $r$ such that no connection in the system
contains $r$ distinct random variables. 

This completes the construction of $\mathbf{P}^{*}$ and $\mathbf{M}^{*}$.
As an example, Fig.~\ref{tab:The-Boolean-matrix-1} shows $\mathbf{P}^{*}$
and $\mathbf{M}^{*}$ for the system $\mathcal{A}$ of Fig.~\ref{fig:system 2by2},
whose associated $\mathbf{P}$ and $\mathbf{M}$ are shown in Fig.~\ref{tab:The-Boolean-matrix}.
Since all connections in this system contain just two binary variables
(as in any cyclic system), the expanded connection probabilities in
this case are uniquely determined by the 1-marginal (bunch) probabilities
and the joint probabilities with which the coupled variables attain
the value 1. In more complex systems the expanded connection probabilities
for maximal couplings can be specified in an infinity of ways.

It can now be shown that the rows in \textbf{$\mathbf{M^{*}}$} are
linearly independent. Indeed, consider a linear combination of these
rows that equals the null vector.
\begin{equation}
\alpha_{1}\left(row_{1}\right)+\alpha_{2}\left(row_{2}\right)+\ldots+\alpha_{N_{rows}}\left(row_{N_{rows}}\right)=\mathbf{0}.\label{eq: ex-linear combination}
\end{equation}
The first row of \textbf{$\mathbf{M^{*}}$} consists of $1$'s only,
and this includes the entry $1$ in the column of \textbf{$\mathbf{M^{*}}$}
corresponding to the hidden outcome $v$ with all elements in it equal
to $k$. Any other row in \textbf{$\mathbf{M^{*}}$} contains zero
in this column. Indeed, entry 1 would have meant that the corresponding
probability in $\mathbf{P}^{*}$ was computed for at least one random
variable attaining a value belonging to $v$: $S_{j}^{i}=v_{j}^{i}=k$.
But this is impossible because the value $k$ is not used in any of
the probabilities in $\mathbf{P}^{*}$. It follows that $\alpha_{1}=0$.
Without loss of generality therefore, we can eliminate this row from
consideration.

The row corresponding to a 1-marginal $\Pr\left[S_{j}^{i}=l\right]$
($l<k$) contains $1$ in the column corresponding to the hidden outcome
$v$ with $v_{j}^{i}=l$ and other entries equal to $k$. All other
rows of \textbf{$\mathbf{M^{*}}$}, now that we have eliminated the
first row, contain zero in the column corresponding to this $v$.
Indeed, to have 1 for this $v$ in another row would have meant that
the corresponding probability in $\mathbf{P}^{*}$ was computed for
the conjunction of $S_{j}^{i}=l$ with at least one other random variable
attaining a value belonging to $v$: $S_{j'}^{i}=v_{j'}^{i}$ or $S_{j}^{i'}=v_{j}^{i'}$.
But all other elements of $v$ equal $k$, and the value $k$ is not
used in any of the probabilities in $\mathbf{P}^{*}$. It follows
that the $\alpha$-coefficients in (\ref{eq: ex-linear combination})
are zero for all rows corresponding to 1-marginal probabilities. Consequently
we can consider all these rows eliminated.

The row corresponding to a 2-marginal $\Pr\left[S_{j}^{i}=l,S_{j'}^{i'}=l'\right]$
(where $i=i'$ if this is a bunch probability, or $j=j'$ if it is
a connection probability, $l,l'<k$) contains $1$ in the column corresponding
to the hidden outcome $v$ with $v_{j}^{i}=l,v_{j'}^{i'}=l'$ and
all other values equal to $k$. All other rows in this column, now
that we have eliminated the first row and all 1-marginal rows, contain
zero in the column corresponding to this $v$. Indeed, to have 1 for
this $v$ in another row would have meant that the corresponding probability
in $\mathbf{P}^{*}$ was computed for the conjunction of $S_{j}^{i}=l,S_{j'}^{i'}=l'$
with at least one other random variable attaining a value belonging
to $v$: $S_{j''}^{i''}=v_{j''}^{i''}$ (where $i''=i'=i$ or $j''=j'=j$),
which is impossible as all other elements of $v$ equal $k$. It follows
that the $\alpha$-coefficients in (\ref{eq: ex-linear combination})
are zero for all rows corresponding to 2-marginal (bunch and connection)
probabilities. Consequently we can consider all these rows eliminated.

Proceeding in this manner to higher-order marginals until all of them
are exhausted, we prove that the rows in \textbf{$\mathbf{M^{*}}$}
are linearly independent. It follows that the system of equations
$\mathbf{M}^{*}\mathbf{Q}=\mathbf{P}^{*}$ always has solutions for
$\mathbf{Q}$ with real-number components (summing to 1). Since $\mbox{\textbf{M}}$
and $\mathbf{P}$ in the original system of linear equations $\mathbf{M}\mathbf{Q}=\mathbf{P}$
associated with a given \mbox{c-c} system are obtained as one and
the same linear combination of the rows of, respectively, \textbf{$\mathbf{M^{*}}$}
and $\mathbf{P}^{*}$, any solution of $\mathbf{M}^{*}\mathbf{Q}=\mathbf{P}^{*}$
is also a solution of $\mathbf{M}\mathbf{Q}=\mathbf{P}$.
\begin{thm}
\label{thm: quasi}Any modified-and-expanded system of equations $\mathbf{M^{*}Q=P}^{*}$
associated with a \mbox{c-c} system has a solution for $\mathbf{Q}$
whose components are real numbers summing to 1. Any such solution
is also a solution for the original system of equations $\mathbf{MQ=P}$
associated with the same \mbox{c-c} system.
\end{thm}
In relation to the possible generalization of CbD discussed in the
concluding section of this paper, note that nowhere in the proof of
this theorem did we use the fact that a quasi-coupling $S$ for the
system $\mathcal{R}$ is maximally connected. In other words, the
proof and the construction of \textbf{$\mathbf{M^{*}}$} and $\mathbf{P}^{*}$
make no use of the maximality constraint (\ref{eq: maximality constarint}).
The only fact that matters is that every connection taken separately
is coupled in some way, so that all connection probabilities in $\mathbf{P}^{*}$
are well-defined.

\subsection{Quasi-probabilities and quasi-couplings}

Let us call the components $\gamma\left(v\right)$ of $\mathbf{Q}$\emph{
}in Theorem \ref{thm: quasi}\emph{ (signed) quasi-probability masses},
and let us call the function $\gamma$ \emph{(signed) quasi-probability
distribution}. Using this terminology, the system of linear equations
$\mathbf{M^{*}Q}=\mathbf{P}^{*}$ (hence also $\mathbf{MQ=P}$) always
produces (generally an infinity of) quasi-probability distributions
of hidden outcomes as solutions for $\mathbf{Q}$.

If the quasi-probability distribution of the hidden outcomes is not
a proper probability distribution, then it does not define a coupling
for the \mbox{c-c} system we are dealing with. However, we can introduce
the notion of a (maximally connected) \emph{quasi-coupling}, by replicating
the definition of a (maximally connected) coupling, but with all references
to probabilities being replaced with quasi-probabilities. 

A \emph{quasi-random variable} $X$ in general is defined as a pair
\begin{equation}
X=\left(\id X,\qdi X\right),
\end{equation}
where $\id X$ is as before and $\qdi X$ is a\emph{ }quasi-probability
distribution function mapping a (finite) set $V_{X}$ of possible
values of $X$ into the set of reals,
\begin{equation}
\qdi X:V_{X}\rightarrow\mathbb{R}.
\end{equation}
The only constraint is that
\begin{equation}
\sum_{v\in V_{X}}\qdi X\left(v\right)=1.
\end{equation}
For any subset $V$ of $V_{X}$ we define quasi-probabilities\footnote{This is a specialization of the measure-theoretic notion of \emph{signed
measure} (or \emph{charge}) to probability spaces with finite sets. }
\[
\qPr\left[X\in V\right]=\sum_{v\in V}\qdi X\left(v\right).
\]
The rest of the conceptual set-up (the class $\E$ generated from
a base set $\R$, the notion of jointly distributed quasi-random variables,
their marginals, functions, etc.) precisely parallels one for the
ordinary, or proper probability distributions and random variables.
We can safely omit details.
\begin{defn}
A \emph{quasi-coupling} $S_{\mathcal{R}}$ of a \mbox{c-c} system
$\mathcal{R}$ is defined as a set of jointly distributed quasi-random
variables in a one-to-one correspondence with the union of the components
of the bunches of $\mathcal{R}$, such that the quasi-probability
distribution of every marginal of $S_{\mathcal{R}}$ that corresponds
to a bunch of the system coincides with the (proper) distribution
of this bunch. A quasi-coupling $S_{\mathcal{R}}$ of $\mathcal{R}$
is \emph{maximally connected} if every marginal of $S_{\mathcal{R}}$
that corresponds to a connection of the system is a maximal coupling
for this connection.
\end{defn}
Let us illustrate this definition on the contextual rank 2 cyclic
system shown in Fig.~\ref{fig: QQ system}. The set of hidden outcomes
here consists of the 16 four-component combinations of 1's and -1's
shown in Fig.~\ref{tab:The-Boolean-matrix-1}. The contextuality
of this system means that these hidden outcomes cannot be assigned
proper probabilities. We can, however, assign real numbers to these
outcomes as shown in Fig.~\ref{tab:The-Boolean-matrix-2}. 
\begin{enumerate}
\item Taking the scalar product of this vector of numbers with the first
row (i.e., summing these numbers), we get 1. This shows that our assignment
of the numbers is a quasi-probability distribution, so we can consider
a quasi-random variable $S$ with this distribution. 
\item Taking the scalar product of the quasi-probability masses with the
subsequent six rows, we get numerical values that are equal to the
corresponding (proper) probabilities characterizing the bunches of
the system in Fig.~\ref{fig: QQ system}. This shows that $S$ is
a quasi-coupling of this system. 
\item Finally, the scalar products of the quasi-probability masses with
the last two rows yield the values of the probabilities characterizing
the maximal couplings for the connections of the system in Fig.~\ref{fig: QQ system}.
This shows that $S$ is a maximally connected quasi-coupling of this
system. 
\end{enumerate}
\begin{figure*}
\begin{footnotesize}

\begin{tabular}{|c|c|c|c|c|c|c|c|c|c|c|c|c|c|c|c|}
\hline 
\multicolumn{16}{|c|}{$\mathbf{M}^{*}$}\tabularnewline
$\!\begin{array}{cc}
+ & +\\
+ & +
\end{array}\!$ & $\!\begin{array}{cc}
+ & +\\
+ & -
\end{array}\!$ & $\!\begin{array}{cc}
+ & +\\
- & +
\end{array}\!$ & $\!\begin{array}{cc}
+ & +\\
- & -
\end{array}\!$ & $\!\begin{array}{cc}
+ & -\\
+ & +
\end{array}\!$ & $\!\begin{array}{cc}
+ & -\\
+ & -
\end{array}\!$ & $\!\begin{array}{cc}
+ & -\\
- & +
\end{array}\!$ & $\!\begin{array}{cc}
+ & -\\
- & -
\end{array}\!$ & $\!\begin{array}{cc}
- & +\\
+ & +
\end{array}\!$ & $\!\begin{array}{cc}
- & +\\
+ & -
\end{array}\!$ & $\!\begin{array}{cc}
- & +\\
- & +
\end{array}\!$ & $\!\begin{array}{cc}
- & +\\
- & -
\end{array}\!$ & $\!\begin{array}{cc}
- & -\\
+ & +
\end{array}\!$ & $\!\begin{array}{cc}
- & -\\
+ & -
\end{array}\!$ & $\!\begin{array}{cc}
- & -\\
- & +
\end{array}\!$ & $\!\begin{array}{cc}
- & -\\
- & -
\end{array}\!$\tabularnewline
\hline 
\hline 
1 & 1 & 1 & 1 & 1 & 1 & 1 & 1 & 1 & 1 & 1 & 1 & 1 & 1 & 1 & 1\tabularnewline
\hline 
1 & 1 & 1 & 1 & 1 & 1 & 1 & 1 &  &  &  &  &  &  &  & \tabularnewline
\hline 
1 & 1 &  &  & 1 & 1 &  &  & 1 & 1 &  &  & 1 & 1 &  & \tabularnewline
\hline 
1 &  & 1 &  & 1 &  & 1 &  & 1 &  & 1 &  & 1 &  & 1 & \tabularnewline
\hline 
1 &  & 1 &  & 1 &  & 1 &  & 1 &  & 1 &  & 1 &  & 1 & \tabularnewline
\hline 
1 & 1 & 1 & 1 &  &  &  &  &  &  &  &  &  &  &  & \tabularnewline
\hline 
1 &  &  &  & 1 &  &  &  & 1 &  &  &  & 1 &  &  & \tabularnewline
\hline 
1 & 1 &  &  & 1 & 1 &  &  &  &  &  &  &  &  &  & \tabularnewline
\hline 
1 &  & 1 &  &  &  &  &  & 1 &  & 1 &  &  &  &  & \tabularnewline
\hline 
\end{tabular}%
\begin{tabular}{|c|}
\hline 
$\mathbf{Q}$\tabularnewline
\hline 
\hline 
$\frac{35}{256}$\tabularnewline
\hline 
$\frac{69}{256}$\tabularnewline
\hline 
$\frac{11}{32}$\tabularnewline
\hline 
$-\frac{1}{4}$\tabularnewline
\hline 
$-\frac{1}{8}$\tabularnewline
\hline 
$\frac{7}{32}$\tabularnewline
\hline 
$-\frac{1}{16}$\tabularnewline
\hline 
$-\frac{1}{32}$\tabularnewline
\hline 
$-\frac{1}{128}$\tabularnewline
\hline 
$-\frac{1}{64}$\tabularnewline
\hline 
$\frac{7}{256}$\tabularnewline
\hline 
$-\frac{1}{256}$\tabularnewline
\hline 
$-\frac{1}{256}$\tabularnewline
\hline 
$\frac{7}{256}$\tabularnewline
\hline 
$\frac{49}{256}$\tabularnewline
\hline 
$\frac{73}{256}$\tabularnewline
\hline 
\end{tabular}%
\begin{tabular}{c}
\multirow{11}{*}{}\tabularnewline
\tabularnewline
\tabularnewline
\tabularnewline
\tabularnewline
\tabularnewline
\tabularnewline
\tabularnewline
\tabularnewline
\tabularnewline
\tabularnewline
\end{tabular}%
\begin{tabular}{|c|}
\multicolumn{1}{c}{}\tabularnewline
\multicolumn{1}{c}{$\mathbf{P}^{*}$}\tabularnewline
\hline 
$1$\tabularnewline
\hline 
$\nicefrac{1}{2}$\tabularnewline
\hline 
$\nicefrac{1}{2}$\tabularnewline
\hline 
$\nicefrac{1}{2}$\tabularnewline
\hline 
$\nicefrac{1}{2}$\tabularnewline
\hline 
$\nicefrac{1}{2}$\tabularnewline
\hline 
$0$\tabularnewline
\hline 
$\nicefrac{1}{2}$\tabularnewline
\hline 
$\nicefrac{1}{2}$\tabularnewline
\hline 
\end{tabular}\end{footnotesize}

\caption{\label{tab:The-Boolean-matrix-3}The same as Fig.~\ref{tab:The-Boolean-matrix-2},
but with the quasi-probabilities that ensure the smallest possible
value of the quasi-coupling's total variation.}
\end{figure*}

\begin{figure*}
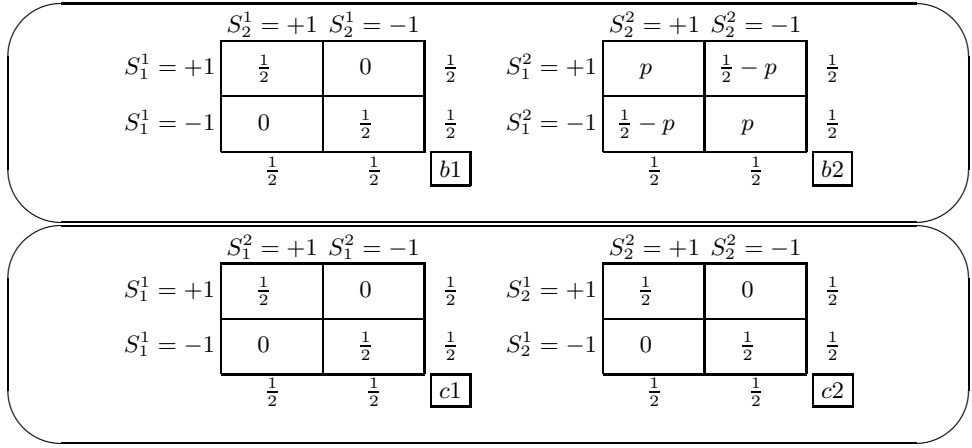

\begin{centering}
\ovalbox{\begin{minipage}[t]{0.7\textwidth}%
\begin{center}
\begin{tabular}{cccc}
 & $S_{2}^{1}=+1$ & $S_{2}^{1}=-1$ & \tabularnewline
\cline{2-3} 
\multicolumn{1}{c|}{$S_{1}^{1}=+1$} & \multicolumn{1}{c|}{$\frac{1}{2}$$\begin{array}{c}
\\
\\
\end{array}$} & \multicolumn{1}{c|}{$0$$\begin{array}{c}
\\
\\
\end{array}$} & $\frac{1}{2}$\tabularnewline
\cline{2-3} 
\multicolumn{1}{c|}{$S_{1}^{1}=-1$} & \multicolumn{1}{c|}{$0$$\begin{array}{c}
\\
\\
\end{array}$} & \multicolumn{1}{c|}{$\frac{1}{2}$$\begin{array}{c}
\\
\\
\end{array}$} & $\frac{1}{2}$\tabularnewline
\cline{2-3} 
 & $\frac{1}{2}$ & $\frac{1}{2}$ & $\boxed{b1}$\tabularnewline
 &  &  & \tabularnewline
\end{tabular}$\quad$%
\begin{tabular}{cccc}
 & $S_{2}^{2}=+1$ & $S_{2}^{2}=-1$ & \tabularnewline
\cline{2-3} 
\multicolumn{1}{c|}{$S_{1}^{2}=+1$} & \multicolumn{1}{c|}{$p$$\begin{array}{c}
\\
\\
\end{array}$} & \multicolumn{1}{c|}{$\frac{1}{2}-p$$\begin{array}{c}
\\
\\
\end{array}$} & $\frac{1}{2}$\tabularnewline
\cline{2-3} 
\multicolumn{1}{c|}{$S_{1}^{2}=-1$} & \multicolumn{1}{c|}{$\frac{1}{2}-p$$\begin{array}{c}
\\
\\
\end{array}$} & \multicolumn{1}{c|}{$p$$\begin{array}{c}
\\
\\
\end{array}$} & $\frac{1}{2}$\tabularnewline
\cline{2-3} 
 & $\frac{1}{2}$ & $\frac{1}{2}$ & $\boxed{b2}$\tabularnewline
 &  &  & \tabularnewline
\end{tabular}
\par\end{center}%
\end{minipage}}
\par\end{centering}

\begin{centering}
\ovalbox{\begin{minipage}[t]{0.7\textwidth}%
\begin{center}
\begin{tabular}{cccc}
 & $S_{1}^{2}=+1$ & $S_{1}^{2}=-1$ & \tabularnewline
\cline{2-3} 
\multicolumn{1}{c|}{$S_{1}^{1}=+1$} & \multicolumn{1}{c|}{$\frac{1}{2}$$\begin{array}{c}
\\
\\
\end{array}$} & \multicolumn{1}{c|}{$0$$\begin{array}{c}
\\
\\
\end{array}$} & $\frac{1}{2}$\tabularnewline
\cline{2-3} 
\multicolumn{1}{c|}{$S_{1}^{1}=-1$} & \multicolumn{1}{c|}{$0$$\begin{array}{c}
\\
\\
\end{array}$} & \multicolumn{1}{c|}{$\frac{1}{2}$$\begin{array}{c}
\\
\\
\end{array}$} & $\frac{1}{2}$\tabularnewline
\cline{2-3} 
 & $\frac{1}{2}$ & $\frac{1}{2}$ & $\boxed{c1}$\tabularnewline
 &  &  & \tabularnewline
\end{tabular}$\quad$%
\begin{tabular}{cccc}
 & $S_{2}^{2}=+1$ & $S_{2}^{2}=-1$ & \tabularnewline
\cline{2-3} 
\multicolumn{1}{c|}{$S_{2}^{1}=+1$} & \multicolumn{1}{c|}{$\frac{1}{2}$$\begin{array}{c}
\\
\\
\end{array}$} & \multicolumn{1}{c|}{$0$$\begin{array}{c}
\\
\\
\end{array}$} & $\frac{1}{2}$\tabularnewline
\cline{2-3} 
\multicolumn{1}{c|}{$S_{2}^{1}=-1$} & \multicolumn{1}{c|}{$0$$\begin{array}{c}
\\
\\
\end{array}$} & \multicolumn{1}{c|}{$\frac{1}{2}$$\begin{array}{c}
\\
\\
\end{array}$} & $\frac{1}{2}$\tabularnewline
\cline{2-3} 
 & $\frac{1}{2}$ & $\frac{1}{2}$ & $\boxed{c2}$\tabularnewline
 &  &  & \tabularnewline
\end{tabular}
\par\end{center}%
\end{minipage}}
\par\end{centering}

\caption{\label{fig: QQ system-1}A cyclic system of rank 2 with the value
of $p$ between 0 and $\frac{1}{2}$. The degree of contextuality
in this system decreases as $p$ increases, and the system becomes
(trivially) noncontextual at $p=\frac{1}{2}$. }
\end{figure*}

\subsection{Contextuality measure based on quasi-couplings}

Proper probability distributions are quasi-probability distributions
with no negative quasi-probability masses. If $X$ is a quasi-random
variable, the value
\begin{equation}
\left\Vert X\right\Vert =\sum_{v\in V_{X}}\left|\qdi X\left(v\right)\right|
\end{equation}
is known in the theory of signed measures as the \emph{total variation}
of $\qdi X$ (or simply of $X$). Its value is 1 if $X$ is a proper
random variable. Otherwise $\left\Vert X\right\Vert >1$, and the
excess of $\left\Vert X\right\Vert $ over 1 can be thought of as
a measure of ``improperness'' of $X$. 

Consider the set of all maximally connected quasi-couplings $S_{\mathcal{R}}$
of a \mbox{c-c} system. We are now interested in the total variations
$\left\Vert S_{\mathcal{R}}\right\Vert $ of (the distributions of)
the quasi-couplings. The set of these values is bounded from below
by 1, therefore it has an infimum, 
\begin{equation}
t=\inf\left\Vert S_{\mathcal{R}}\right\Vert .
\end{equation}
It can be readily seen that this infimum is in fact a minimum of all
$\left\Vert S_{\mathcal{R}}\right\Vert $, i.e., that the set of all
$S_{\mathcal{R}}$ contains a quasi-coupling $S_{\mathcal{R}}^{*}$
with
\begin{equation}
\left\Vert S_{\mathcal{R}}^{*}\right\Vert =t.
\end{equation}
Indeed, if the set of all $S_{\mathcal{R}}$ is finite, the statement
is trivially true; otherwise we choose an infinite sequence of $\left\Vert S_{\mathcal{R}}\right\Vert $
converging to $t$. Without loss of generality, we can assume for
all members of this sequence 
\begin{equation}
\left\Vert S_{\mathcal{R}}\right\Vert -t<\varepsilon.
\end{equation}
It follows that the quasi-probability masses $\qdi S_{\mathcal{R}}\left(v\right)$
for all hidden outcomes $v$ in all these quasi-couplings are confined
within a closed interval $\left[-t-\varepsilon,t+\varepsilon\right]$.
So the quasi-probability distributions $\qdi S_{\mathcal{R}}$ (viewed
as vectors of real numbers)\footnote{Note that $\qdi S_{\mathcal{R}}$ can be viewed as the same entity
as the vector $\mathbf{Q}$ in the matrix equation $\mathbf{MQ}=\mathbf{P}$,
because $\mathbf{Q}$ is a vector of real numbers indexed by the hidden
outcomes. There is a subtlety here (and throughout this paper) related
to distinguishing indexed values and the pairs consisting of indexes
and values, but we will ignore it. } are confined within a cube $\left[-t-\varepsilon,t+\varepsilon\right]^{N}$,
where $N$ is the number of the hidden outcomes. Since the cube is
compact, from the sequence of $S_{\mathcal{R}}$ one can choose a
converging subsequence, with the limit $S_{\mathcal{R}}^{*}$, and
it is easy to see that $\left\Vert S_{\mathcal{R}}^{*}\right\Vert $
cannot exceed $t$ (otherwise the original sequence of $\left\Vert S_{\mathcal{R}}\right\Vert $
would have converged to two distinct limits).

So, the set of quasi-couplings $S_{\mathcal{R}}$ for any \mbox{c-c}
system $\mathcal{R}$ contains a quasi-coupling $S_{\mathcal{R}}^{*}$
with the smallest possible value of the total variation $\left\Vert S_{\mathcal{R}}\right\Vert $.
If $\left\Vert S_{\mathcal{R}}^{*}\right\Vert $ equals 1, then the
system is noncontextual, because then $S_{\mathcal{R}}^{*}$ is a
proper maximally connected coupling. If $\left\Vert S_{\mathcal{R}}^{*}\right\Vert >1$,
then no proper maximally connected coupling for $\mathcal{R}$ exists,
and the quantity $\left\Vert S_{\mathcal{R}}^{*}\right\Vert -1$ can
be taken as a measure (of degree) of contextuality. Note that while
the minimum total variation $\left\Vert S_{\mathcal{R}}^{*}\right\Vert $
is unique, the quasi-coupling $S_{\mathcal{R}}^{*}$ generally is
not.

Consider again Fig.~\ref{tab:The-Boolean-matrix-2}. The value of
$\left\Vert S_{\mathcal{R}}\right\Vert $ in it is $6\cdot\frac{1}{2}=3$.
Is this the smallest possible value? It is not. A direct minimization
of $\left\Vert S_{\mathcal{R}}\right\Vert $ subject to the linear
equations $\mathbf{MQ}=\mathbf{P}$ shows the minimum value in the
case of the system depicted in Fig.~\ref{fig: QQ system} to be 2.
It is reached, e.g., in the quasi-probability distribution shown in
Fig.~\ref{tab:The-Boolean-matrix-3}. This distribution therefore
defines an $S_{\mathcal{R}}^{*}$. There are other quasi-probability
distributions (in fact an infinity of them) with this minimum value
of $\left\Vert S_{\mathcal{R}}\right\Vert $.

It is instructive to see how this total variation measure changes
as we change the value of $p=\Pr\left[S_{1}^{2}=1,S_{2}^{2}=1\right]$
from zero to the maximal possible value $\frac{1}{2}$ while keeping
all other probabilities fixed (see Fig.~\ref{fig: QQ system-1}).
The relationship turns out to be linear:
\begin{equation}
\left\Vert S_{\mathcal{R}}^{*}\right\Vert =2\left(1-p\right).
\end{equation}
The system is maximally contextual at $p=0$, the case we focused
on in our examples. When $p$ reaches $\frac{1}{2}$, the system is
noncontextual: trivially so, because then its two bunches are identical. 

A direct minimization of $\left\Vert S_{\mathcal{R}}\right\Vert $
subject to the linear equations $\mathbf{MQ}=\mathbf{P}$ is a nonlinear
problem, but it can be reduced to a linear programming one, in the
following way:
\begin{enumerate}
\item Create a matrix $\mathbf{M}_{wide}$ by horizontally concatenating
$\mathbf{M}$ and $\mathbf{M'}=\left(-1\right)\cdot\mathbf{M}$,
\begin{equation}
\mathbf{M}_{wide}=\left(\mathbf{M}\,\brokenvert\,\mathbf{M'}\right)
\end{equation}
 Each hidden outcome $v$ labels two columns of $\mathbf{M}_{wide}$
(one in the $\mathbf{M}$ half and one in the $\mathbf{M'}$ half).
\item Create a column vector $\mathbf{Q}_{long}$ whose length is twice
that of $\mathbf{Q}$,
\begin{equation}
\mathbf{Q}_{long}=\left(\begin{array}{c}
\mathbf{Q}_{1}\\
--\\
\mathbf{Q}_{2}
\end{array}\right)
\end{equation}
 Its elements are labelled in the same way as the columns of $\mathbf{M}_{wide}$.
\item Solve the linear programming problem 
\[
\mathbf{M}_{wide}\mathbf{Q}_{long}=\mathbf{P}
\]
 subject to three constraints: (a) nonnegativity of the components
of $\mathbf{Q}_{long}$, and (b) minimality of the sum of the components
of $\mathbf{Q}_{2}$.
\end{enumerate}
To every hidden outcome $v$ there correspond two elements of $\mathbf{Q}_{long}$,
denoted $\gamma^{+}\left(v\right)$ and $\gamma^{-}\left(v\right)$,
and the quasi-probability mass assigned to $v$ is $\gamma^{+}\left(v\right)-\gamma^{-}\left(v\right)$.
The sum of these quasi-probabilities across all $v$ equals 1, and
the sum of their absolute values is the minimal total variation $\left\Vert S_{\mathcal{R}}^{*}\right\Vert $.
The reasoning above (the proof that $S_{\mathcal{R}}^{*}$ always
exists) guarantees that this linear programming problem always has
a solution, generally non-unique.

\section{Conclusion}

In this paper we have described the basic elements of a theory aimed
at analyzing systems of random variables classified in two ways: by
their conteXts and by their conteNts. Irrespective of one's terminological
preferences, the classification of such systems into contextual and
noncontextual ones, as well as into consistently connected and inconsistently
connected ones, is meaningful and, at least in some applications of
the theory, fundamentally important. We would like to conclude this
paper by recapitulating a few points made in this paper and by offering
a general observation.

\paragraph{(1)}

Inconsistency of connectedness should be distinguished from contextuality.
One may call inconsistency ``a kind of contextuality,'' but it is
contextuality of a different kind. Inconsistency of connectedness
is about direct influences of certain elements of conteXts upon random
variables. Such influences are revealed on the level of stochastically
unrelated random variables sharing the same conteNt (i.e., within
connections). Direct influences cannot act, e.g., from the future
to the past or from one event to a spatially separated but simultaneous
one. Contextuality, by contrast, is revealed in joint distributions
of random variables, and it is not constrained by considerations of
causality. As an example, consider a physical realization of the Suppes-Zanotti-Legett-Garg
system (Leggett and Garg, 1985) that consists in three measurements
made at three moments in time with respect to some zero point, as
shown below:
\[
\xymatrix{\ar@{-}+<0ex,0ex>;[rrrr] & \overset{}{\underset{\begin{array}{c}
t_{1}\end{array}}{\bullet}} & \overset{}{\underset{\begin{array}{c}
t_{2}\end{array}}{\bullet}} & \overset{}{\underset{\begin{array}{c}
t_{3}\end{array}}{\bullet}} & \,}
\]
The measurements are always made in pairs: at moments $t_{1}$ and
$t_{2}$ or at moments $t_{1}$ and $t_{3}$ or at moments $t_{2}$
and $t_{3}$. Each pair of times moments defines a conteXt, and each
moment defines a conteNt (because the measurement in this analysis
are distinguished only by the time moments at which they are made).
Now, it is perfectly possible that the distribution of $R_{t_{2}}^{\left(t_{1},t_{2}\right)}$
differs from the distribution of $R_{t_{2}}^{\left(t_{2},t_{3}\right)}$,
because in the former case the measurement at moment $t_{2}$ can
be directly influenced by the fact that a measurement was made at
some moment in the past, $t_{1}$ (if the system is a quantum one,
its quantum state, prepared at moment zero, can be changed by a previous
measurement); but a measurement cannot be influenced by another measurement
yet to be made at a future moment, $t_{3}$. By contrast, in any joint
distributions of the variables, such as $\left(R_{t_{1}}^{\left(t_{1},t_{2}\right)},R_{t_{2}}^{\left(t_{1},t_{2}\right)}\right)$,
the future random variable stochastically depends on the past one
exactly whenever the past one stochastically depends on the future
one. Contextuality is only revealed by looking at such joint distributions
within bunches and comparing them across bunches.

\paragraph{(2) }

A distinguishing feature of Contextuality-by-Default, and the main
reason for the ``by-default'' in its name, is that it treats random
variables in different conteXts as different random variables, even
if they have conteNts in common. As a result, the bunches of a system
never overlap, and the problem of contextuality therefore is not posed
as a problem of compatibility of different overlapping groups of random
variables. Rather it is posed as a problem of compatibility between
the bunches on the one hand and maximal couplings for the connections
on the other. In this respect Contextuality-by-Default is distinct
from other approaches to contextuality, e.g., the prominent line of
contextuality research by Abramsky and his colleagues (Abramsky \&
Brandenburger, 2011; Abramsky et al, 2015).

\paragraph{(3) }

Treating random variables in different conteXts as different, however,
in no way means that conteXts are fused (or confused) with conteNts.
On the contrary, Contextuality-by-Default is based on a strict differentiation
of these entities, although, being an abstract mathematical theory,
it cannot determine what constitutes conteNts and conteXts in a given
empirical situation. This determination is made before the theory
applies. If one changes one's double-classification of the random
variables (by the conteNts and by the conteXts), the contextuality
of the system changes too.

\paragraph{(4) }

Contextuality-by-Default is not a model for empirical phenomena. As
any abstract mathematical theory, it has no predictive power as a
result of having no predictive intent. It is a theoretical language,
on a par with, say, real analysis or probability theory. In fact,
if presented in full generality to include arbitrary systems of arbitrary
random variables (the presentation in this paper was confined to finite
sets of categorical variables only), Contextuality-by-Default is essentially
co-extensive with Kolmogorovian theory of random variables. The main
difference from Kolmogorovian probability theory is that the Contextuality-by-Default
theory may (but does not have to) be constructed without sample spaces,
that it prominently uses the notion of stochastic unrelatedness (implicit
or underemphasized in Kolmogorovian probability theory), and that
the theory of couplings (rather peripheral to the mainstream Kolmogorovian
theory) is at the very heart of the Contextuality-by-Default theory.

\paragraph{(5) }

In dealing with contextuality, the Contextuality-by-Default theory
is about compatibility (or lack thereof) of the observed bunches of
a system with maximal couplings of the separately taken connections.
Maximality is not, however, the only possible constraint imposable
on the couplings of the connections. Contextuality-by-Default can
be expanded or modified in various ways by replacing it with other
constraints, and any new constraint replacing maximality would tackle
a new meaning of contextuality. Using the same logic as in Section
\ref{sub: The-intuition-of} and in Definition \ref{def: cntx}, if
separately taken connections can be coupled subject to some constraint
$\mathsf{C}$, then the system is ``$\mathsf{C}$-noncontextual''
if it can be coupled so that all subcouplings corresponding to its
connections satisfy $\mathsf{C}$; otherwise the system is ``$\mathsf{C}$-contextual.''\textquotedbl{}
It is remarkable that the representation of the contextuality problem
as a linear programming task (Section \ref{sec: Contextuality-as-LP})
and the construction of the measure of contextuality based on the
quasi-couplings (Section \ref{sec: How-to-measure}) apply with no
modifications to any choice of $\mathsf{C}$ such that a coupling
satisfying $\mathsf{C}$ exists for any connection taken separately.
Indeed, the only property of connection probabilities required for
the construction of the matrix-vector pair $\mathbf{M^{*}}$-$\mathbf{P}^{*}$
(hence also $\mathbf{M}$-$\mathbf{P}$) is that these probabilities
exist, not the way they are computed. (The choices of $\mathsf{C}$
for which a coupling satisfying $\mathsf{C}$ may not exist for some
connections taken separately requires a modification in the definition
of contextuality, but can be handled too: any system possessing these
connections can be treated as ``automatically'' contextual.) In
choosing a constraint $\mathsf{C}$ to replace maximality, one can
be guided by certain reasonable desiderata, one of them being that
$\mathsf{C}$ should be reduced to the identity constraint when a
system is consistently connected. Another reasonable desideratum could
be that the ``$\mathsf{C}$-theory'' reduces to the one described
in Section \ref{sec: Cyclic-c-c-systems} when specialized to cyclic
systems with binary random variables. We will elaborate elsewhere.\footnote{\label{fn: As-the-reviewing}As the reviewing of this paper was nearing
completing and no substantive changes could be made, we proposed a
new version of CbD, with $\mathsf{C}$-couplings being ``multimaximal''
ones (Dzhafarov and Kujala, 2016a,b): in such a coupling of a connection
any of its subcouplings is a maximal coupling of the corresponding
subset of the connection. If random variables in a connection are
binary, their multimaximal coupling always exists and is unique. For
connections with more-than-binary categorical variables, one possible
approach is to replace them with all their possible dichotomizations;
in each conteXt, these dichotomizations are jointly distributed and
form a sub-bunch of the bunch corresponding to the context. The replacement
of maximality with multimaximality affects classification of the systems
into contextual and noncontextual, but it does not affect the validity
of our theorems related to the measure of contextuality, as the maximality
constraint was not used in their proofs. The specializations of CbD
to consistently connected systems and to cyclic systems with binary
random variables remain unchanged.}

\subsection*{Acknowledgments.}

This research has been supported by NSF grant SES-1155956, AFOSR grant
FA9550-14-1-0318, and A. von Humboldt Foundation. We are grateful
to Victor H. Cervantes of Purdue University for his insights on maximal
couplings that helped us in the linear programming treatment of contextuality
and the construction of a measure of contextuality. The latter was
also inspired by the use of quasi-probability distributions (``negative
probabilities'') in dealing with contextuality by Samson Abramsky
of Oxford University and J. Acacio de Barros of San Francisco State
University. We greatly benefited from numerous discussions with them
and their colleagues. We are grateful to Matt Jones of the University
of Colorado whose critical analysis of our treatment of contextuality
helped us to improve the motivation and argumentation for our approach
to contextuality. Victor H. Cervantes and Farzin Shamloo of Purdue
University were most helpful in discussing and finding imprecisions
and typos in earlier versions of the paper.

\end{document}